\documentclass[final,3p,times,nonatbib]{elsarticle}

\usepackage[english]{babel}
\usepackage{graphicx}
\usepackage{siunitx}
\usepackage[utf8]{inputenc}
\usepackage[T1]{fontenc}
\usepackage[labelfont=bf,format=plain,font=small]{caption}
\usepackage[labelformat=simple]{subcaption}

\usepackage{algorithmic}
\usepackage{algorithm}
\renewcommand{\algorithmiccomment}[1]{\bgroup\hfill//~#1\egroup}
\usepackage{tikz}
\usetikzlibrary{arrows.meta,arrows}
\usetikzlibrary{patterns,positioning,shapes.geometric}
\usepackage{pgfplots}
\pgfplotsset{compat=1.13}
\tikzstyle{every node}=[font=\small]
\pgfplotsset{grid style={gray!30,thin}}
\pgfplotsset{axis line style={thin}}
\pgfplotsset{every tick label/.append style={font=\scriptsize}}
\tikzset{every picture/.style={line width=1.2pt}}        
\tikzset{every node/.style={font=\small}}    
\tikzset {_0x7o5ee8j/.code = {\pgfsetadditionalshadetransform{ \pgftransformshift{\pgfpoint{0 bp } { 0 bp }  }  \pgftransformrotate{0 }  \pgftransformscale{2 }  }}}
\pgfdeclarehorizontalshading{_sfax6qxjp}{150bp}{rgb(0bp)=(0.82,0.82,0.82);
rgb(37.5bp)=(0.82,0.82,0.82);
rgb(62.5bp)=(0.47,0.47,0.47);
rgb(100bp)=(0.47,0.47,0.47)}

\usepackage[
    bookmarks,
    bookmarksopen=false,
    colorlinks=true,
    linkcolor=orange,
    anchorcolor=brown,
    citecolor=blue,
    filecolor=magenta,
    menucolor=red,
    urlcolor=cyan,
    pagebackref=false,
    plainpages=false,
    pdfpagelabels,
    hypertexnames=true,
    final
]{hyperref}

\usepackage{amsmath,amssymb,amsthm}
\usepackage{dsfont}
\usepackage{bm}
\usepackage{mathtools}

\cslet{blx@noerroretextools}\empty
\usepackage[backend=biber,style=apa,natbib,maxcitenames=2,mincitenames=1,maxbibnames=99,minbibnames=98,sorting=ynt,doi=false,isbn=false,url=false,eprint=false]{biblatex}
\usepackage{booktabs}

\usepackage[nolist]{acronym}

\usepackage{my_macros}
\makeatletter
\renewcommand*{\AC@verridelabel}[1]{}
\makeatother

\usepackage{aliascnt}
\usepackage[noabbrev,capitalize]{cleveref}

\usepackage{autonum}
\usepackage{csquotes}

\newtheoremstyle{kursiv}{}{}{\itshape}{}{\bfseries}{}{ }{\thmname{#1}\thmnumber{ #2}{\normalfont\thmnote{ (#3)}}:}
\newtheoremstyle{normal}{}{}{}{}{\bfseries}{}{ }{\thmname{#1}\thmnumber{ #2}{\normalfont\thmnote{ (#3)}}:}
\theoremstyle{kursiv}
\newtheorem{theorem}{Theorem}
\newaliascnt{corollary}{theorem}
\newtheorem{corollary}[corollary]{Corollary}
\aliascntresetthe{corollary}
\newaliascnt{lemma}{theorem}
\newtheorem{lemma}[lemma]{Lemma}
\aliascntresetthe{lemma}
\newaliascnt{proposition}{theorem}
\newtheorem{proposition}[proposition]{Proposition}
\aliascntresetthe{proposition}
\newaliascnt{definition}{theorem}
\newtheorem{definition}[definition]{Definition}
\aliascntresetthe{definition}
\theoremstyle{normal}
\newaliascnt{remark}{theorem}
\newtheorem{remark}[remark]{Remark}
\aliascntresetthe{remark}
\crefname{theorem}{Theorem}{Theorems}
\crefname{corollary}{Corollary}{Corollaries}
\crefname{lemma}{Lemma}{Lemmas}
\crefname{proposition}{Proposition}{Propositions}
\crefname{definition}{Definition}{Definitions}
\crefname{remark}{Remark}{Remarks}

\journal{Journal of Computational Physics}

\begin{document}
\begin{frontmatter}
\title{Material Optimization of Optical Properties of Particulate Products by Transition Matrix Method and Sequential Global Programming}

\author[a]{Nico Nees\corref{cor1}}
\ead{nico.nees@fau.de}
\author[a,b]{Lukas Pflug}
\author[a]{Michael Stingl}
\cortext[cor1]{Corresponding author}
\affiliation[a]{organization={Friedrich-Alexander-Universit\"at Erlangen-N\"urnberg (FAU), Department of Mathematics, Chair of Applied Mathematics (Continuous Optimization)},
                addressline={Cauerstra\ss{}e 11},
                postcode={91058},
                city={Erlangen},
                country={Germany}}
\affiliation[b]{organization={Friedrich-Alexander-Universit\"at Erlangen-N\"urnberg (FAU), FAU Competence Center Scientific Computing (FAU CSC)},
                addressline={Martensstra\ss{}e 5a},
                postcode={91058},
                city={Erlangen},
                country={Germany}}

\begin{abstract}
Material optimization of nanoparticle assemblies, such as structurally coloured photonic crystals, attracts increasing interest. We combine the transition matrix method (TMM), which avoids spatial discretization of particles and surrounding medium, with sequential global programming (SGP) to optimize optical properties of particulate products by iteratively assigning predefined catalogue materials to individual particles.

We derive exact separable models, reproducing the fully coupled objective under single-particle material changes, and computationally accessible first-order correct approximations for general continuously differentiable objectives. Unlike conventional gradient-based models, these approximations aim to capture the nonlinear, potentially nonconvex dependence on individual material variables beyond their local first-order behaviour. Within SGP, globally minimizing the particle-wise subproblems on a prescribed candidate grid can yield substantial improvements in relatively few outer iterations, including steps away from unfavourable stationary points, as previously demonstrated for discrete dipole approximation\nbhy based SGP. The construction exploits the decomposition of the TMM system matrix into a material-dependent block-diagonal and a material-independent coupling matrix, while the state and adjoint problems remain fully coupled.

Applications include extinction, scattering quantities, and colour objectives. In a \num{35}\nbhy particle extinction minimization problem, the objective decreases by \qty{31.0}{\percent} in \num{13} iterations. A fine-grid evaluation of the fully coupled objective shows that no tested single-particle material change improves the final design, providing a stronger numerical optimality test than first-order stationarity alone. For a \num{627}\nbhy particle photonic crystal, chroma increases by \qtyrange{58.4}{122.1}{\percent} for different particle diameters in \numrange{17}{42} iterations. The optimized designs favour absorbing material in layers away from the incident light, consistent with experimental findings.
\end{abstract}

\end{frontmatter}

\section{Introduction}\label{section:introduction}
Structural colours result from the interaction of visible light with organized nanostructures and can be modified by introducing absorbing material into a particulate assembly \citep{kinoshita2005,vignolini2012,goerlitzer2018,carina2024,bleyer2025}. In photonic crystals, the spatial distribution of absorber can suppress unwanted scattering and alter the observed colour. This leads to the following material optimization question: given the positions of the particles, which material should be assigned to each particle to improve a prescribed optical property? We study this as a numerical optimization problem governed by fully coupled electromagnetic scattering, without imposing fabrication constraints on individual particle assignments.

The forward model must account for multiple scattering among the particles. The \acf{tmm} uses vector spherical wave functions to describe the response of each particle through its transition matrix and couples these responses through translation operators \citep{waterman1965,Waterman1971,Mishchenko2002,egel2017,BEUTEL2024109076}. It does not require a volume discretization of the particles or of the surrounding medium. By comparison, the \acf{fem} usually discretizes a bounded computational domain and an absorbing boundary treatment \citep{Monk2003}, while the \acf{dda} represents the scatterer by interacting dipoles \citep{draine1994}. FFT\nbhy based \ac{dda} implementations may use regular dipole grids whose inactive sites also affect computational cost. The \ac{tmm} is therefore attractive for assemblies of separated particles with available individual T\nbhy matrices, subject to the usual dependence of cost and accuracy on particle geometry and multipole truncation. \ac{tmm}\nbhy related methods have also been used to optimize particle geometry and arrangement \citep{zhan2018,zhelyeznyakov2020,Semmler2021,asadova2026}.

The \acf{sgp} method replaces the original material optimization problem by a sequence of separable subproblems \citep{semmler2018,Nees2022,vu2023}. Separability alone is shared with ordinary first-order Taylor models. The distinguishing feature of \ac{sgp} is that it seeks computationally tractable approximations to exact separable models, which reproduce the objective function whenever only one design coordinate is varied. Such approximations can preserve nonlinear and potentially nonconvex dependence on individual material variables that the local gradient does not capture. Their coordinate-wise accuracy also motivates a stronger optimality criterion than first-order stationarity, which also applies to finite, discrete material catalogues: no change of the material assigned to a single particle should improve the true objective. When accurate enough, their globally minimized scalar subproblems can propose substantial improvements and lead to convergence in comparatively few outer iterations. In particular, the published \ac{dda}\nbhy \ac{sgp} study by \citet[Section~6.1, Figure~11]{Nees2022} documents progress through several unfavourable Karush--Kuhn--Tucker (KKT) points, in contrast to the behaviour of a local convex-model method in that example, although escaping stationary points is not guaranteed in general. \ac{sgp} was originally developed also for discrete material optimization with several admissible materials \citep{semmler2018}. Here we use a continuously relaxed bi-material catalogue.

A particle-wise \ac{dda}\nbhy \ac{sgp} formulation was already considered in \citet{Nees2022}. There, dipoles belonging to one particle are combined into a single material design variable, although the surrogate originates from dipole-level updates. This grouping can diminish the model's ability to represent particle-wise material changes. In contrast, the \ac{tmm} system admits one material-dependent diagonal block per particle, with the interparticle coupling represented by material-independent off-diagonal blocks. A material change in one particle therefore gives a low-rank matrix update and leads directly to accurate particle-wise separable models. Only an auxiliary inverse of the system matrix is approximated in the computationally tractable surrogate, while the state and adjoint equations remain fully coupled.

The main contribution of this work is the explicit construction and analysis of exact separable \ac{tmm} models and of computationally accessible first-order correct approximations that retain the nonlinear single-particle dependence in the fully coupled \ac{tmm} setting. The construction applies to general continuously differentiable functions of finitely many linear scattering observables and encompasses several wavelengths and incident polarizations. An additional linearization of the monitor scattering intensity reduces the number of adjoint solves needed for colour objectives.

The numerical study first examines extinction minimization for \num{35} particles. Besides demonstrating accurate surrogates and convergence in \num{13} \ac{sgp} iterations, it tests, by fully coupled evaluations on \num{201} admissible material values per particle, whether any sampled single-particle change can improve the final design. This coordinate-wise optimality check is stronger than merely testing first-order KKT stationarity, though it does not certify optimality on the entire continuous material interval. A second example maximizes chroma for a \num{627}\nbhy particle photonic crystal with wavelength-dependent materials, obtaining substantial improvements in \numrange{17}{42} iterations. The prescribed optical material data are informed by earlier experimental studies \citep{carina2024,bleyer2025}, but exact fabrication of the optimized particle-by-particle distributions is not assumed.

\Cref{section:general_scattering_problem,section:TMM} summarize the scattering observables and the \ac{tmm} model. \Cref{section:material_optimization} introduces the design variables and \ac{sgp} framework. The derivation of the general separable models, followed by their specialization to optical and colour objectives, is given in \cref{section:TMM_SGP_models_general}. The numerical results follow in \cref{section:numerical_experiments}. Standard \ac{tmm} relations and detailed algebraic steps are collected in the appendices. Part of the \ac{tmm} background and related preliminary numerical material appeared in the doctoral dissertation by \citet{nees2025Thesis}. The present paper is self-contained with respect to the model construction and its validation.

\section{Electromagnetic Scattering Theory}\label{section:general_scattering_problem}
In the following, we introduce the general electromagnetic scattering problem and the optical properties of particulate products that are considered throughout this work.

\subsection{General Electromagnetic Scattering Problem}\label{section:general_electromagnetic_scattering_problem}
We start from the curl-curl formulation of the time-harmonic Maxwell's equations, which is the standard setting for the description of electromagnetic scattering by bounded objects (see, e.g., \citet{Monk2003,Angell2004,colton2013b,stratton2015,Hohenester2020}{}).

In this work, we employ the standard Cartesian coordinate system in $\RR^3$ with the origin $\mathbf{0} = (0,0,0)^\top$. The $i$\nbhy th standard unit vector is written as~$\eb^i$, and the entries of a vector~$\mathbf{u}$ and a matrix~$\mathbf{K}$ are denoted by $\mathbf{u}_i = \eb^i\cdot\mathbf{u}$ and $\mathbf{K}_{ij}$, respectively. For complex-valued quantities, $\overline{z}$ denotes the complex conjugate of $z\in\CC$, $\mathbf{K}^\dagger = \overline{\mathbf{K}}^\top$ the conjugate transpose of a matrix~$\mathbf{K}$, and $\Re$ and $\Im$ the operators that map a complex number to its real and imaginary part, respectively. Furthermore, $\mathbf{u}\cdot\mathbf{v} = \sum_i \mathbf{u}_i\mathbf{v}_i$ denotes the dot product, which is taken without complex conjugation, $\mathbf{u}\times\mathbf{v}$ the cross product, $\Vert\cdot\Vert_p$ the $p$\nbhy norm for a parameter $p\geq 1$ with the Euclidean norm $\Vert\cdot\Vert_2$ as a special case, and $\nabla\cdot$ and $\nabla\times$ the divergence and the curl operator. Finally, the imaginary unit is represented by $\ci = \sqrt{-1}$, and the vacuum wavenumber associated with a wavelength $\lambda\in\RR_{>0}$ is defined by
\begin{align}\label{eq:wavenumber}
    k_\lambda = \frac{2\pi}{\lambda}.
\end{align}

The scatterer is a bounded domain $\Omega = \bigcup_{i=1}^\Npar \Omega_i \subset \RR^3$ consisting of $\Npar$ particles~$\Omega_i\subset\RR^3$ that do not overlap each other, i.e., $\Omega_i\cap\Omega_j = \emptyset$ for $i\neq j$. Without loss of generality, the object is assumed to be centred at the origin. It is characterized by its relative permittivity $\epsilon_\mathrm{r}\colon\Omega\to\CC$, and therefore by the complex refractive index~$\nbar$, and it is embedded in a homogeneous background medium with relative permittivity $\epsilon_\mathrm{b} = \nbar_\mathrm{b}^2$. Throughout this work, the background medium is assumed to be non-absorbing, i.e., $\nbar_\mathrm{b}\in\RR_{>0}$. Moreover, all materials are assumed to be non-magnetic, i.e., the relative permeability satisfies $\mu_\mathrm{r}\equiv 1$ on~$\RR^3$. The propagation of both the incident and the scattered wave outside the scatterer is then governed by the wavenumber of the background medium,
\begin{align}\label{eq:wavenumber_background}
    k_\mathrm{b} = \sqrt{\epsilon_\mathrm{b}}\,k_\lambda = \nbar_\mathrm{b}\frac{2\pi}{\lambda},
\end{align}
which is real and positive by the above assumption.

As incident fields we consider plane electromagnetic waves. For a propagation direction $\db\in\RR^3$ and a polarization $\pb\in\RR^3$ with $\Vert\db\Vert_2 = \Vert\pb\Vert_2 = 1$ and $\pb^\top\db = 0$, the incident electric field $\boldsymbol{E}_\mathrm{inc}\colon\RR^3\to\CC^3$ of unit amplitude is given at a position $\xb\in\RR^3$ by
\begin{align}\label{eq:plane_wave}
    \boldsymbol{E}_\mathrm{inc}(\xb) = \pb \exp\bigl(\ci k_\mathrm{b} \xb^\top \db\bigr),
\end{align}
i.e., a wave oscillating in the direction~$\pb$ and travelling in the direction~$\db$. The associated magnetic field follows from Maxwell's equations as $\boldsymbol{H}_\mathrm{inc} = \frac{1}{\ci k_\lambda}\nabla\times\boldsymbol{E}_\mathrm{inc} = \nbar_\mathrm{b}\,(\db\times\pb)\exp(\ci k_\mathrm{b}\xb^\top\db)$. Consequently, the incident irradiance~$I_\mathrm{inc}$, i.e., the modulus of the time-averaged energy flux carried by the incident wave, is the constant $I_\mathrm{inc} = \frac{1}{2}\nbar_\mathrm{b}$. In the absence of the scatterer, this plane wave solves Maxwell's equations in the background medium~$\epsilon_\mathrm{b}$, see, e.g., \citet{Monk2003}{}.

\begin{figure}[!htb]
    \centering
    \input{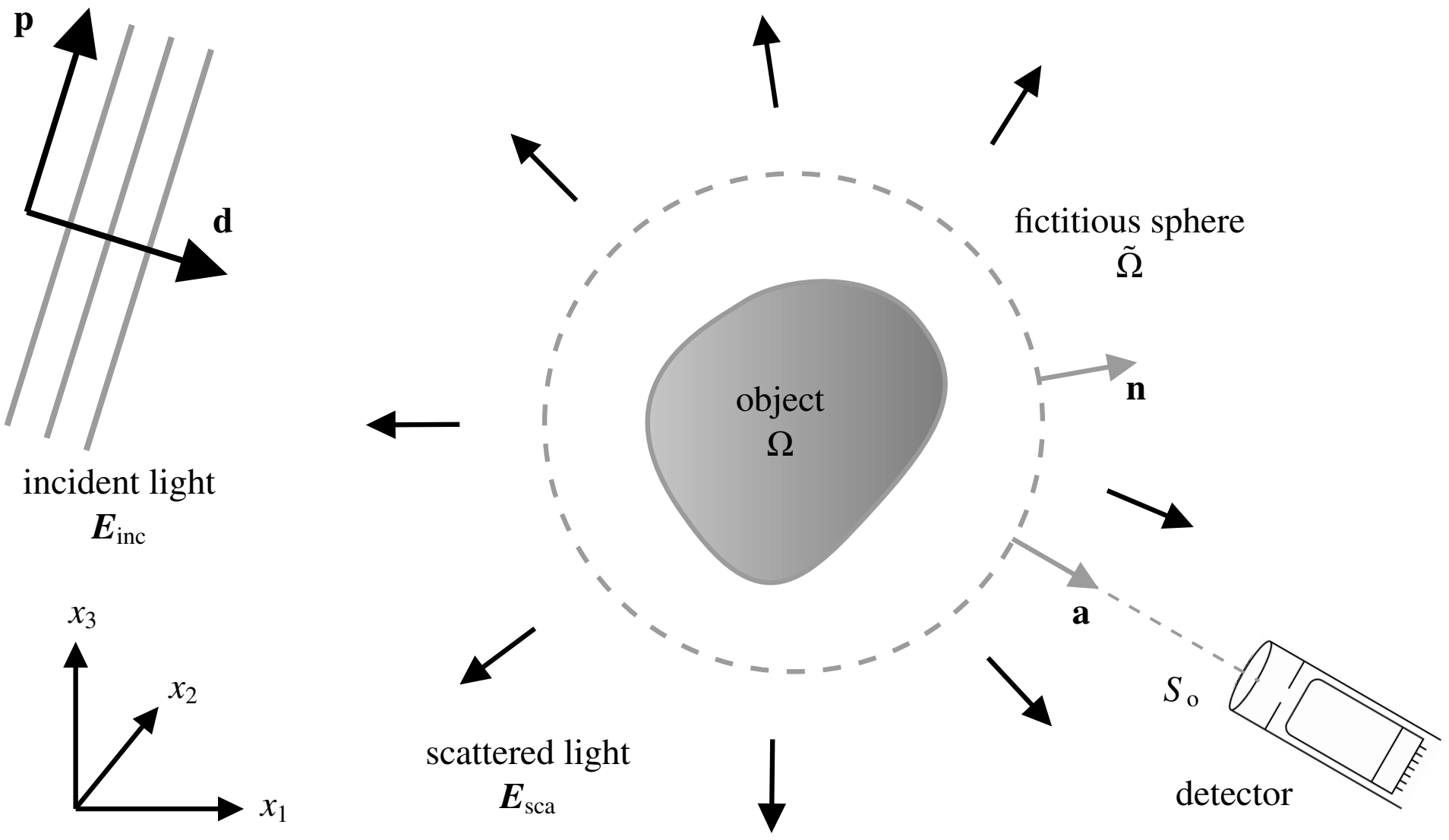}
    \caption{General electromagnetic scattering problem in which the object~$\Omega$ is illuminated by an incident light~$\boldsymbol{E}_\mathrm{inc}$ with propagation direction~$\db$ and polarization~$\pb$. When the incident field interacts with the scatterer, some light may be scattered in different directions, which is indicated by the scattered field~$\boldsymbol{E}_\mathrm{sca}$. The scatterer's orientation as well as the direction and polarization of the incident light are characterized by a reference coordinate system. Furthermore, the fictitious sphere $\tilde\Omega$ is visualized together with an exemplary local outward unit vector~$\mathbf{n}$. A well-collimated detector in the far-field zone, indicated by the grey dashed line, with the detection surface~$S_\mathrm{o}$ is aligned with the scattering direction~$\ab$.}
    \label{fig:scattering_problem}
\end{figure}

On impact with the object, some energy of the incident wave may be absorbed and some deflected. The physical process of deflecting energy is called scattering and is described by the scattered field~$\boldsymbol{E}_\mathrm{sca}$, which propagates in different directions depending on the material properties as well as the size and shape of the object. The total electric field is given by the superposition $\boldsymbol{E}_\mathrm{tot} = \boldsymbol{E}_\mathrm{inc}+\boldsymbol{E}_\mathrm{sca}$ of the incident and the scattered field, and analogously for the magnetic fields. This setting is visualized schematically in \cref{fig:scattering_problem}, together with a fictitious sphere~$\tilde\Omega$ centred on and surrounding the scatterer, whose surface~$S_{\tilde\Omega}$ serves as observation sphere in \cref{section:monitor_scattering}, and a well-collimated detector in the far-field zone.

The scattered field is the solution of the so-called electromagnetic scattering problem. Inside and outside the scatterer it satisfies the curl-curl equation with the respective permittivity, where the interaction of the incident field with the object acts as a source term. At the boundary $\Gamma = \partial\Omega$ with outer unit normal~$\nub$, appropriate interface and transmission conditions have to be imposed, as outlined by \citet[Section~1.2.2]{Monk2003}{}. Finally, the Silver--Müller radiation condition \citep{mueller1969} guarantees that the scattered field is outgoing, i.e., that it originates at the scatterer and propagates outward. Following \citet[Section~1.4.3]{Monk2003}{}, the scattered electric field $\boldsymbol{E}_\mathrm{sca}\colon\RR^3\to\CC^3$ thus solves
\begin{subequations}\label{eq:ME_scat_problem}
    \begin{align+}
        \nabla \times \bigl(\nabla \times \boldsymbol{E}_\mathrm{sca}(\xb)\bigr) - k_\lambda^2\epsilon_\mathrm{r}(\xb) \boldsymbol{E}_\mathrm{sca}(\xb) &= -k_\lambda^2\bigl(\epsilon_\mathrm{b}-\epsilon_\mathrm{r}(\xb)\bigr) \boldsymbol{E}_\mathrm{inc}(\xb) && \text{in } \Omega,\label{eq:ME_scat_problem_1}\\
        \nabla \times \bigl(\nabla \times \boldsymbol{E}_\mathrm{sca}(\xb)\bigr) - k_\lambda^2\epsilon_\mathrm{b} \boldsymbol{E}_\mathrm{sca}(\xb) &= \mathbf{0} && \text{in } \RR^3 \setminus \Omega, \label{eq:ME_scat_problem_2}\\
        \bigl[\boldsymbol{E}_\mathrm{sca}(\xb) \times \nub\bigr] = \bigl[\bigl(\nabla \times \boldsymbol{E}_\mathrm{sca}(\xb)\bigr) \times \nub\bigr] &= \mathbf{0} && \text{on } \Gamma, \label{eq:ME_scat_problem_3}
    \end{align+}
\end{subequations}
together with the Silver--Müller radiation condition
\begin{align+}\label{eq:ME_silver_mueller}
    \lim_{\Vert\xb\Vert_2\to\infty} \Bigl(\bigl(\nabla \times \boldsymbol{E}_\mathrm{sca}(\xb)\bigr)\times \xb -\ci k_\mathrm{b} \Vert\xb\Vert_2\boldsymbol{E}_\mathrm{sca}(\xb)\Bigr) = \mathbf{0}.
\end{align+}

\subsection{Optical Properties of Particulate Products}\label{section:optical_properties}
The solution of the scattering problem gives access to a wide range of optical properties of a particulate product, see, e.g., \citet{Mishchenko2002,Hohenester2020}{}. In this manuscript, we focus on three properties, which are introduced in the following sections. The first one is the extinction cross-section, a single scalar quantity that measures how much energy the object removes from the incident beam and that is obtained directly from the solution of the scattering problem. The second one is the monitor scattering intensity, i.e., the light energy captured by a detector of finite aperture. Its evaluation requires the scattered field in a whole set of scattering directions and therefore builds on the directional scattering intensity introduced below. The third one is the colour of the object, which in addition requires the monitor scattering intensity for a whole spectrum of wavelengths as well as its conversion into colour values.

\subsubsection{Extinction Cross-Section}\label{section:cross_section}
Let us begin with one of the most prominent optical properties, the so-called extinction cross-section. When the incident wave interacts with the object, energy is removed from the incident beam, either because it is absorbed by the particles and subsequently transformed into another form of energy, such as heat, or because it is deflected into other directions. The associated net rates of energy are the extinction energy~$W_\mathrm{ext}$, the absorbed energy~$W_\mathrm{abs}$, and the scattered energy~$W_\mathrm{sca}$. They are obtained by integrating the corresponding time-averaged energy flux, given by the respective Poynting vector, over a closed surface surrounding the scatterer, see \citet[Equation~(1.31)]{Mishchenko2002} or \citet[Section~2]{Frezza2018}{}. Since the background medium is assumed to be non-absorbing, the three quantities are linked by the energy balance $W_\mathrm{ext} = W_\mathrm{abs} + W_\mathrm{sca}$, i.e., the energy that is subtracted from the incident field is either absorbed by the object or scattered by it.

Normalizing these energies by the incident irradiance~$I_\mathrm{inc}$ yields the extinction, absorption, and scattering cross-section,
\begin{align}\label{eq:scattering_cross_section}
    \sigma_\mathrm{ext} = \frac{1}{I_\mathrm{inc}} W_\mathrm{ext},\qquad
    \sigma_\mathrm{abs} = \frac{1}{I_\mathrm{inc}} W_\mathrm{abs},\qquad
    \sigma_\mathrm{sca} = \frac{1}{I_\mathrm{inc}} W_\mathrm{sca},
\end{align}
respectively. All three have the dimension of an area and inherit the relation $\sigma_\mathrm{ext} = \sigma_\mathrm{abs}+\sigma_\mathrm{sca}$.

\subsubsection{Directional Scattering}\label{section:directional_scattering}
While the cross-sections of the previous section are integral quantities, many applications require information about the angular distribution of the scattered light. Let us therefore consider a so-called well-collimated detector \citep[Section~11.4]{Mishchenko2014}, which only measures light from a very narrow range of angles, in the far-field zone. For an illustration, see \cref{fig:scattering_problem}. The scattered light intensity detected by such an external detector or observer in a specific direction is referred to as directional scattering or angle-dependent scattering intensity.

Let $\ab\in\RR^3$ be a unit vector specifying the scattering direction of interest. Following the idea of \citet[Section~2.5]{Mishchenko2002}, the observer's detection surface~$S_\mathrm{o}$ is a plane centred on and perpendicular to~$\ab$, placed at a distance~$R$ from the scatterer. If the detector is far away from the scatterer and its surface is very small, the electromagnetic power it receives is essentially the flux of the scattered field through this surface. This power still depends on the size of the detector and on its distance from the scatterer. Normalizing it by the incident irradiance~$I_\mathrm{inc}$ and by the solid angle $A(S_\mathrm{o})/R^2$ subtended by the detector removes this dependence and leaves a quantity that characterizes the scatterer alone, the differential scattering cross-section
\begin{align}\label{eq:directional_scattering_general}
    s_\ab = \frac{\mathrm{d}\sigma_\mathrm{sca}}{\mathrm{d}\Omega}(\ab) = \lim_{R\to\infty} \frac{R^2}{I_\mathrm{inc}}\,\boldsymbol{S}_\mathrm{sca}(R\,\ab)\cdot\ab ,
\end{align}
which is well defined because the scattered energy flux decays like $R^{-2}$ in the far field. Throughout this work, the directional scattering intensity~$s_\ab$ denotes this differential scattering cross-section with the dimension of an area per solid angle. Here, $\boldsymbol{S}_\mathrm{sca}\colon\RR^3\to\RR^3$ denotes the time-averaged Poynting vector of the scattered wave, i.e., the energy flux carried by the scattered field~$(\boldsymbol{E}_\mathrm{sca},\boldsymbol{H}_\mathrm{sca})$ alone. It is defined at a position $\xb\in\RR^3$ by (\citeauthor{poynting1884}, \citeyear{poynting1884}; \citeauthor{Mishchenko2014}, \citeyear{Mishchenko2014}, Section~2.4)
\begin{align}\label{eq:poynting_vector_scattered}
    \boldsymbol{S}_\mathrm{sca}(\xb) = \frac{1}{2} \Re\bigl(\boldsymbol{E}_\mathrm{sca}(\xb) \times \overline{\boldsymbol{H}_\mathrm{sca}(\xb)}\bigr).
\end{align}

\subsubsection{Monitor Scattering}\label{section:monitor_scattering}
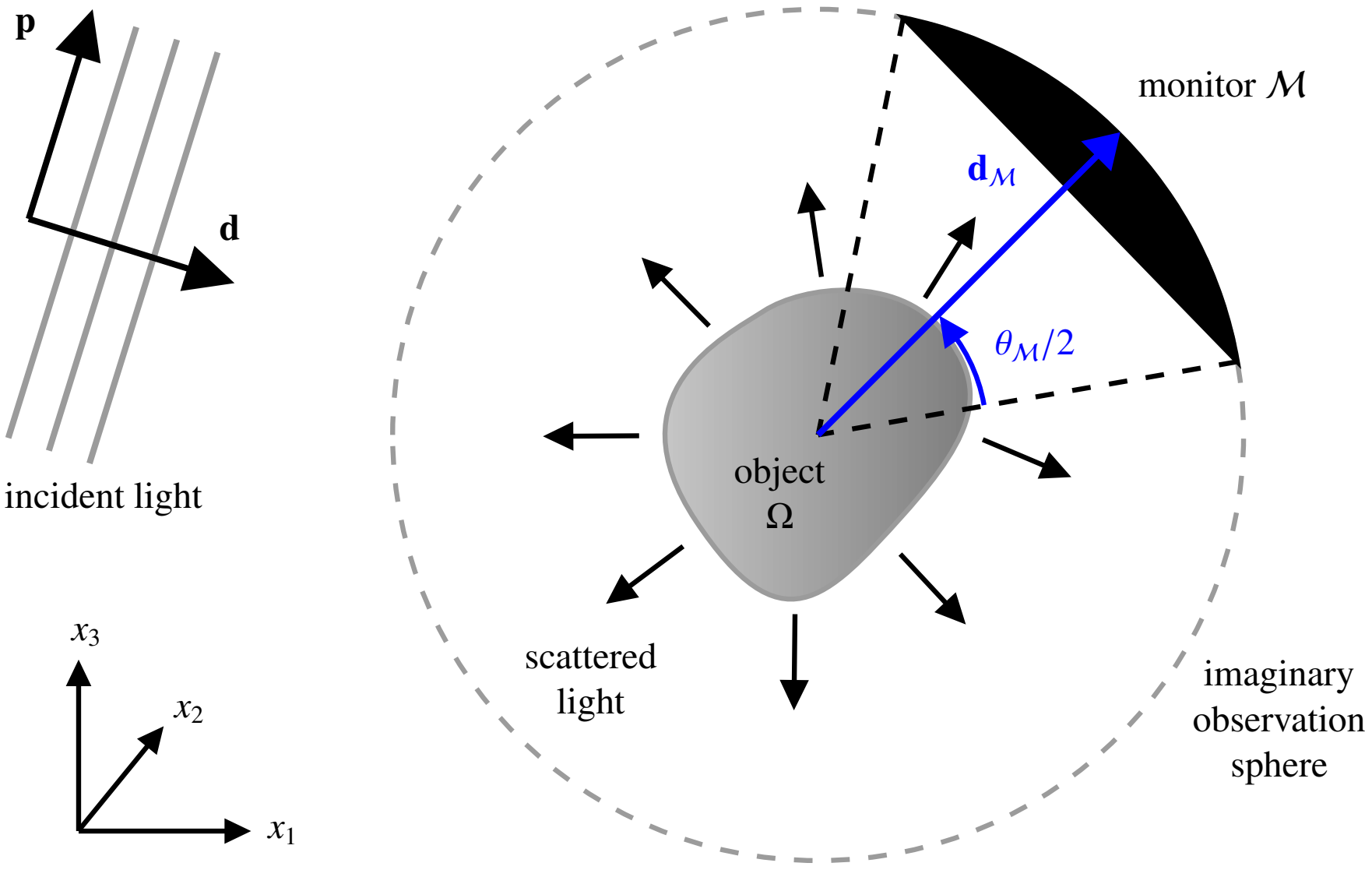
\begin{figure}[!htb]
    \centering
    \begin{tikzpicture}[scale=0.0029357\textwidth,x=0.75pt,y=0.75pt,yscale=-1,xscale=1]
\draw  [color={rgb, 255:red, 155; green, 155; blue, 155 }  ,draw opacity=1 ][dash pattern={on 4.5pt off 4.5pt}] (250,180) .. controls (250,124.77) and (294.77,80) .. (350,80) .. controls (405.23,80) and (450,124.77) .. (450,180) .. controls (450,235.23) and (405.23,280) .. (350,280) .. controls (294.77,280) and (250,235.23) .. (250,180) -- cycle ;
\draw [color={rgb, 255:red, 155; green, 155; blue, 155 }  ,draw opacity=1 ][line width=1.5]    (169.21,183.7) -- (199.32,87.15) ;
\draw [color={rgb, 255:red, 155; green, 155; blue, 155 }  ,draw opacity=1 ][line width=1.5]    (159.67,180.72) -- (189.77,84.18) ;
\draw [color={rgb, 255:red, 155; green, 155; blue, 155 }  ,draw opacity=1 ][line width=1.5]    (178.76,186.68) -- (208.87,90.13) ;
\draw [line width=1.5]    (164.47,129.25) -- (209.02,143.15) ;
\draw [line width=1.5]    (164.47,129.25) -- (178.37,84.7) ;
\draw [shift={(212.84,144.34)}, rotate = 197.32] [fill={rgb, 255:red, 0; green, 0; blue, 0 }  ][line width=0.08]  [draw opacity=0] (11.61,-5.58) -- (0,0) -- (11.61,5.58) -- cycle    ;
\draw [shift={(179.5,80)}, rotate = 107.32] [fill={rgb, 255:red, 0; green, 0; blue, 0 }  ][line width=0.08]  [draw opacity=0] (11.61,-5.58) -- (0,0) -- (11.61,5.58) -- cycle    ;

\path  [shading=_sfax6qxjp,_0x7o5ee8j] (337.27,150.87) .. controls (346.6,144.87) and (368.67,142.33) .. (379.27,155.53) .. controls (389.87,168.73) and (387.93,178.87) .. (369.93,198.87) .. controls (351.93,218.87) and (341.27,230.2) .. (321.27,200.2) .. controls (301.27,170.2) and (327.93,156.87) .. (337.27,150.87) -- cycle ;
\draw  [color={rgb, 255:red, 155; green, 155; blue, 155 }  ,draw opacity=1 ] (337.27,150.87) .. controls (346.6,144.87) and (368.67,142.33) .. (379.27,155.53) .. controls (389.87,168.73) and (387.93,178.87) .. (369.93,198.87) .. controls (351.93,218.87) and (341.27,230.2) .. (321.27,200.2) .. controls (301.27,170.2) and (327.93,156.87) .. (337.27,150.87) -- cycle ;
\draw    (307.93,180.2) -- (288.33,180.32) ;
\draw [shift={(285.33,180.33)}, rotate = 359.66] [fill={rgb, 255:red, 0; green, 0; blue, 0 }  ][line width=0.08]  [draw opacity=0] (7.14,-3.43) -- (0,0) -- (7.14,3.43) -- cycle    ;
\draw    (318.26,206.18) -- (302.55,217.9) ;
\draw [shift={(300.15,219.69)}, rotate = 323.27] [fill={rgb, 255:red, 0; green, 0; blue, 0 }  ][line width=0.08]  [draw opacity=0] (7.14,-3.43) -- (0,0) -- (7.14,3.43) -- cycle    ;
\draw    (375.11,147.8) -- (386.05,130.39) ;
\draw [shift={(387.05,128.69)}, rotate = 122.39] [fill={rgb, 255:red, 0; green, 0; blue, 0 }  ][line width=0.08]  [draw opacity=0] (7.14,-3.43) -- (0,0) -- (7.14,3.43) -- cycle    ;
\draw    (388.7,181.1) -- (406.75,188.74) ;
\draw [shift={(409.52,189.91)}, rotate = 202.95] [fill={rgb, 255:red, 0; green, 0; blue, 0 }  ][line width=0.08]  [draw opacity=0] (7.14,-3.43) -- (0,0) -- (7.14,3.43) -- cycle    ;
\draw    (344.57,222.11) -- (344.35,241.71) ;
\draw [shift={(344.32,244.71)}, rotate = 270.64] [fill={rgb, 255:red, 0; green, 0; blue, 0 }  ][line width=0.08]  [draw opacity=0] (7.14,-3.43) -- (0,0) -- (7.14,3.43) -- cycle    ;
\draw    (350.78,142.87) -- (347.89,123.49) ;
\draw [shift={(347.44,120.52)}, rotate = 81.52] [fill={rgb, 255:red, 0; green, 0; blue, 0 }  ][line width=0.08]  [draw opacity=0] (7.14,-3.43) -- (0,0) -- (7.14,3.43) -- cycle    ;
\draw    (324.38,154.41) -- (310.62,140.45) ;
\draw [shift={(308.51,138.31)}, rotate = 45.41] [fill={rgb, 255:red, 0; green, 0; blue, 0 }  ][line width=0.08]  [draw opacity=0] (7.14,-3.43) -- (0,0) -- (7.14,3.43) -- cycle    ;
\draw    (369.31,208) -- (382.68,222.34) ;
\draw [shift={(384.72,224.53)}, rotate = 226.99] [fill={rgb, 255:red, 0; green, 0; blue, 0 }  ][line width=0.08]  [draw opacity=0] (7.14,-3.43) -- (0,0) -- (7.14,3.43) -- cycle    ;
\draw    (176.2,273.07) -- (213.74,273.07) ;
\draw [shift={(216.74,273.07)}, rotate = 180] [fill={rgb, 255:red, 0; green, 0; blue, 0 }  ][line width=0.08]  [draw opacity=0] (6.25,-3) -- (0,0) -- (6.25,3) -- cycle    ;
\draw    (176.2,273.07) -- (176.2,235.78) ;
\draw [shift={(176.2,232.78)}, rotate = 90] [fill={rgb, 255:red, 0; green, 0; blue, 0 }  ][line width=0.08]  [draw opacity=0] (6.25,-3) -- (0,0) -- (6.25,3) -- cycle    ;
\draw    (176.2,272.94) -- (194.3,250.76) ;
\draw [shift={(196.2,248.44)}, rotate = 129.23] [fill={rgb, 255:red, 0; green, 0; blue, 0 }  ][line width=0.08]  [draw opacity=0] (6.25,-3) -- (0,0) -- (6.25,3) -- cycle    ;
\draw  [color={rgb, 255:red, 0; green, 0; blue, 0 }  ,draw opacity=1 ][fill={rgb, 255:red, 0; green, 0; blue, 0 }  ,fill opacity=1 ] (370.03,82.01) .. controls (410.13,90.16) and (441.51,122.3) .. (448.53,162.8) -- cycle ;
\draw  [dash pattern={on 4.5pt off 4.5pt}]  (370.03,82.01) -- (350,180) ;
\draw  [dash pattern={on 4.5pt off 4.5pt}]  (448.53,162.8) -- (350,180) ;
\draw [color={rgb, 255:red, 0; green, 0; blue, 255 }, line width=1.5 , shift={(350,180)}, rotate = -45 ]  (0,0) -- (99,0) ;
\draw [shift={(350.3+100*0.707106781186547524,79.6+100*(1-0.707106781186547524))}, rotate = 135] [fill={rgb, 255:red, 0; green, 0; blue, 255 }  ][line width=0.3]  [draw opacity=0] (8,-5) -- (0,0) -- (8,5) -- cycle    ;
\draw  [color={rgb, 255:red, 0; green, 0; blue, 255 },shift={(389,173)}] (0,0) arc (-10:-43:39) ; 
\draw  [shift={(378.5,152)}, rotate = 55] [fill={rgb, 255:red, 0; green, 0; blue, 255 }  ][line width=0.08]  [draw opacity=0] (7.14,-3.43) -- (0,0) -- (7.14,3.43) -- cycle    ;

\draw (208,127) node [anchor=north west][inner sep=0.75pt]   [align=left] {$\db$};
\draw (160,80) node [anchor=north west][inner sep=0.75pt]   [align=left] {$\pb$};
\draw (432.67,231.67) node [anchor=north west][inner sep=0.75pt]   [align=left] {\begin{minipage}[lt]{50.17pt}\setlength\topsep{0pt}
\begin{center}{imaginary}\\{observation}\\{sphere}\end{center}\end{minipage}};
\draw (157.5,190) node [anchor=north west][inner sep=0.75pt]   [align=left] {incident light};
\draw (275.58,227.58) node [anchor=north west][inner sep=0.75pt]   [align=left] {\begin{minipage}[lt]{40.48pt}\setlength\topsep{0pt}
\begin{center}{scattered}\\{light}\end{center}\end{minipage}};
\draw (320,184) node [anchor=north west][inner sep=0.75pt]   [align=left] {\begin{minipage}[lt]{40.48pt}\setlength\topsep{0pt}
\begin{center}{object}\\{$\Omega$}\end{center}\end{minipage}};
\draw (219,273.5) node [anchor=west][inner sep=0.75pt] [align=left] {$\displaystyle x_1$};
\draw (197,249) node [anchor=south west][inner sep=0.75pt] [align=left] {$\displaystyle x_2$};
\draw (178,231) node [anchor=south][inner sep=0.75pt] [align=left] {$\displaystyle x_3$};
\node[align=center,color={rgb, 255:red, 0; green, 0; blue, 255 }] at (391,118) {$\mathbf{d}_\mathcal{M}$};
\node[align=center,color={rgb, 255:red, 0; green, 0; blue, 255 }] at (401,159) {$\theta_{\mathcal{M}}/2$};
\draw (420,93) node [anchor=north west][inner sep=0.75pt] [align=left] {\begin{minipage}[lt]{48.53pt}\setlength\topsep{0pt}
\begin{center}{monitor $\displaystyle \mathcal{M}$}\end{center}\end{minipage}};
\end{tikzpicture}
    \caption{General monitor scattering setting. The light energy captured by a monitor~$\Mcal$, a spherical cap, is of interest. The monitor is defined by the position vector of its centre~$\db_\Mcal$ and an opening angle~$\theta_\Mcal$.}
    \label{fig:monitor_scat_setting}
\end{figure}

We want to compute the light intensity captured by a segment of the surface~$S_{\tilde\Omega}$ of the imaginary observation sphere~$\tilde\Omega$, called the monitor $\Mcal\subseteq S_{\tilde\Omega}$. The general setting is illustrated in \cref{fig:monitor_scat_setting}. In this work, we restrict our focus to spherical caps that represent the monitor, i.e., 
\begin{align}
    \Mcal= \Bigl\{\mathbf{x}\in S_{\tilde\Omega}\ \vert\ \langle\mathbf{x},\db_\Mcal \rangle \leq \frac{\theta_\Mcal}{2}\Bigr\},
\end{align}
where $\db_\Mcal\in\RR^3$ denotes the position vector of the centre of the monitor, $\theta_\Mcal\in(0,2\pi]$ the opening angle, and $\langle \mathbf{x},\db_\Mcal \rangle$ the angle between $\mathbf{x}\in\RR^3$ and $\db_\Mcal$. The opening angle~$\theta_\Mcal$ is defined in reference to the configuration shown in \cref{fig:monitor_scat_setting}. In general, it is possible to consider monitors of any shape.

In order to compute the scattering $s_\Mcal$ captured by the monitor~$\Mcal$, which is characterized by the predefined position vector~$\db_\Mcal$ and the opening angle~$\theta_\Mcal$, it is necessary to integrate the directional scattering intensity~$s_\ab$ from \cref{eq:directional_scattering_general} over the solid angle subtended by the monitor. In the following, we identify the monitor~$\Mcal$ with the corresponding set of scattering directions $\{\xb/\Vert\xb\Vert_2\ \vert\ \xb\in\Mcal\}$ on the unit sphere~$S^2$. Therefore, the monitor scattering intensity is given by
\begin{align}\label{eq:monitor_scat_integral}
    s_{\Mcal} = \int_{\Mcal} s_\ab\, \mathrm d \Omega(\ab),
\end{align}
where $\mathrm d \Omega$ denotes the surface measure of the unit sphere, i.e., the solid angle element.
Since $s_\ab$ is a differential scattering cross-section, $s_\Mcal$ is the partial scattering cross-section associated with the monitor, i.e., it has the dimension of an area. The integral can be well approximated by a tensor-product Romberg quadrature rule \citep[Section~6.3]{davis1975} in the polar and the azimuthal angle of a local spherical coordinate system aligned with~$\db_\Mcal$, resulting in a weighted sum of the directional scattering intensity:
\begin{align}\label{eq:directional_scat_Romberg}
    s_{\Mcal} \approx \sum_{t=1}^{\tilde{n}} \omega_t s_{\ab^t}.
\end{align}
Thereby, $\tilde n\in\NN$ denotes the number of integration points or scattering directions~$\ab^t\in\RR^3$ on the monitor~$\Mcal$, and $\omega_t\in\RR$ represent the corresponding integration weights for $t\in\{1,\ldots,\tilde n\}$. Note that the weights~$\omega_t$ comprise both the Romberg weights of the two one-dimensional rules and the surface element~$\sin\theta$ of the local spherical coordinate system.

\begin{remark}[Scattering Cross-Section]
    Choosing $\theta_\Mcal = 2\pi$ results in the integration of the differential scattering cross-section over the whole observation sphere (cf. \cref{fig:monitor_scat_setting}). This yields the scattering cross-section defined by \cref{eq:scattering_cross_section}, i.e., $s_\Mcal = \sigma_\mathrm{sca}$ in this case.
\end{remark}

\begin{remark}[Reflection and Transmission]\label{rem:reflection_transmission}
    Let the angle between the position vector of the monitor~$\db_\Mcal$ and the propagation direction of the incident field~$\db$ be denoted by $\langle \db_\Mcal,\db \rangle$. The light intensity scattered in a monitor placed in forward direction, i.e., for $\langle \db_\Mcal,\db \rangle < \frac{\pi}{2}$ with $\theta_\Mcal \leq \pi - 2\langle \db_\Mcal,\db \rangle$, is called transmission. Conversely, the scattered light captured by a monitor placed in backward direction, so for $\langle \db_\Mcal,\db \rangle > \frac{\pi}{2}$ with $\theta_\Mcal \leq 2\langle \db_\Mcal,\db \rangle - \pi$, is designated as reflection \citep{Bohren1998}.
\end{remark}

In this work, we will focus on the so-called reflectance, which can be derived from the reflection (see \cref{rem:reflection_transmission}). Therefore, let $A_\mathrm{p}(\Omega,\db)$ be the projection area of the scatterer~$\Omega$ onto a plane perpendicular to the incident propagation direction~$\db$. Then, the reflectance~$r_{\Mcal}$ captured by the monitor~$\Mcal$, which is characterized for the reflection by $\langle \db_\Mcal,\db \rangle > \frac{\pi}{2}$ with $\theta_\Mcal \leq 2\langle \db_\Mcal,\db \rangle - \pi$, is computed by
\begin{align}\label{eq:reflectance}
    r_{\Mcal} = \frac{1}{A_\mathrm{p}(\Omega,\db)} s_{\Mcal}.
\end{align}
The reflectance is thus a dimensionless partial scattering efficiency. If $r_{\Mcal} = 0$, no light is reflected or scattered in the direction of the monitor $\Mcal$. 

\subsubsection{Colour Functions}\label{section:Lab_colour}
The following section addresses the computation of an object's colour. Based on a spectrum of a scattering property, we derive formulas for describing colour values in an appropriate colour space. We deal with the so-called CIE\nbhy $\Lcolour\acolour\bcolour$ space defined by the International Commission on Illumination (CIE) \citep{CIE2006}. In this colour space, a colour is defined by a tuple~$(\Lcolour,\acolour,\bcolour)$ with the following meanings:
\begin{itemize}
    \item $\Lcolour$: perceptual lightness of the colour ($\Lcolour = 0 \rightarrow$ black; $\Lcolour = 100 \rightarrow$ diffuse white)
    \item $\acolour$: position between red and green ($\acolour < 0 \rightarrow$ green; $\acolour > 0 \rightarrow$ red)
    \item $\bcolour$: position between yellow and blue ($\bcolour < 0 \rightarrow$ blue; $\bcolour > 0 \rightarrow$ yellow)
\end{itemize}

\begin{figure}[!htb]
    \centering
    \includegraphics[trim = 3.2cm 20.3cm 8.2cm 3.5cm, clip, width=0.61\textwidth]{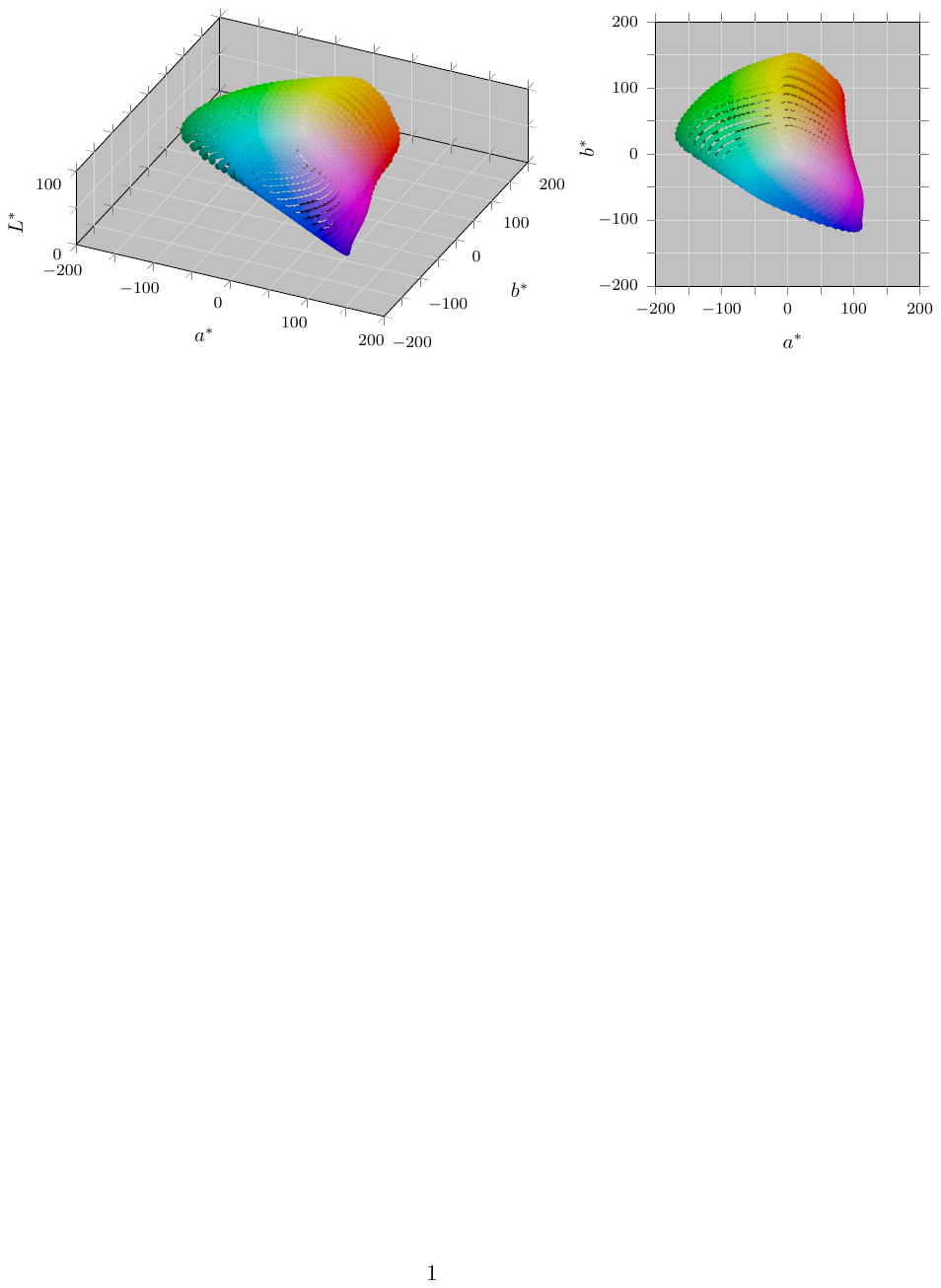}
    \includegraphics[trim = 13cm 20.2cm 1.9cm 3.6cm, clip, width=0.38\textwidth]{pics/Lab_gamut.pdf}
    \caption{Visualization of the CIE\nbhy $\Lcolour\acolour\bcolour$ colour space. In the centre of this space, the colours are low in intensity, whereas at the outer boundary, i.e., for high $\acolour$ or $\bcolour$ values, they are highly intense and vivid. The right plot corresponds to a top view of the colour space. Note that only the outer boundary of the visible gamut of the colour space is illustrated here.}
    \label{fig:CIElab_colour_space}
\end{figure}

For an illustration of the prescribed colour space, see \cref{fig:CIElab_colour_space}.

In order to compute colour values, it is necessary to calculate some light scattering properties for the whole spectrum of visible light, which ranges from approximately \qtyrange{380}{780}{\nm}. Our focus will be on the reflectance~$r_{\Mcal}$, defined in the previous section, where $\Mcal$ defines the monitor. We denote the reflectance at a given wavelength $\lambda \in \RR_{>0}$ by $r_\Mcal(\lambda)$. The CIE\nbhy $\Lcolour\acolour\bcolour$ colour values are then computed according to the CIE standards \citep{klein2010,IOS2019}. Let $\Lambda$ define an entire spectrum. In this work, we use $\Lambda = [\qty{400}{\nm},\qty{800}{\nm}]$, which is contained in the tabulation range $[\qty{390}{\nm},\qty{830}{\nm}]$ of the CIE~2006 colour matching functions employed below. Furthermore, let the functions $x,y,z\colon\Lambda \to \RR$ represent predefined colour matching functions, and $I\colon\Lambda \to \RR$ define a standard illuminant. For the conversion from a given spectrum to the CIE\nbhy $\Lcolour\acolour\bcolour$ values, it is first necessary to compute the so-called CIE\nbhy $\mathit{XYZ}$ values by
\begin{align}
    X &= \frac{1}{\bar{y}} \int_\Lambda x(\lambda) r_\Mcal(\lambda) I(\lambda)\, \mathrm{d}\lambda,\label{eq:CIE_XYZ_1} \\
    Y &= \frac{1}{\bar{y}} \int_\Lambda y(\lambda) r_\Mcal(\lambda) I(\lambda)\, \mathrm{d}\lambda,\label{eq:CIE_XYZ_2} \\
    \intertext{and}
    Z &= \frac{1}{\bar{y}} \int_\Lambda z(\lambda) r_\Mcal(\lambda) I(\lambda)\, \mathrm{d}\lambda\label{eq:CIE_XYZ_3}
\end{align}
with
\begin{align}
    \bar{y} &= \int_\Lambda y(\lambda) I(\lambda)\, \mathrm{d}\lambda.\label{eq:CIE_XYZ_4}
\end{align}
Moreover, let the tuple $(X_\mathrm{w},Y_\mathrm{w},Z_\mathrm{w})$ define a reference white point. In our case, these values are provided by
\begin{align}
    X_\mathrm{w} &= \num{0.9472528492},\label{eq:CIE_white_point_1}\\
    Y_\mathrm{w} &= \num{1.0},\label{eq:CIE_white_point_2} \\
    \intertext{and}
    Z_\mathrm{w} &= \num{1.0713012997},\label{eq:CIE_white_point_3}
\end{align}
and are computed using \cref{eq:CIE_XYZ_1,eq:CIE_XYZ_2,eq:CIE_XYZ_3} with $r_\Mcal(\lambda) = 1$ for all $\lambda \in \Lambda$. Additionally, we define the function $f\colon\RR \to \RR$ by
\begin{align}\label{eq:Lab_colour_f}
    f(x) = \begin{cases}
        \sqrt[3]{x}, & \text{if } x > \varepsilon, \\
        \frac{\kappa x +16}{116}, & \text{otherwise},
    \end{cases}
\end{align}
with CIE parameters $\varepsilon = \frac{216}{24389}$ and $\kappa = \frac{24389}{27}$. Then, the CIE\nbhy $\Lcolour\acolour\bcolour$ colour values are calculated via
\begin{align}
    L^\ast &= 116 f\Bigl(\frac{Y}{Y_\mathrm{w}}\Bigr) - 16,\label{eq:CIE_Lab_1} \\
    a^\ast &= 500 \biggl(f\Bigl(\frac{X}{X_\mathrm{w}}\Bigr) - f\Bigl(\frac{Y}{Y_\mathrm{w}}\Bigr)\biggr),\label{eq:CIE_Lab_2}\\
    \intertext{and}
    b^\ast &= 200 \biggl(f\Bigl(\frac{Y}{Y_\mathrm{w}}\Bigr) - f\Bigl(\frac{Z}{Z_\mathrm{w}}\Bigr)\biggr).\label{eq:CIE_Lab_3}
\end{align}

In our studies, the CIE\nbhy $\Lcolour\acolour\bcolour$ values are determined using the so-called CIE D65 standard illuminant \citep{CIE2022}, commonly referred to as the ``Daylight 6500K''. This illuminant is used to represent typical daylight conditions at noon in western Europe. Furthermore, there are multiple options available for the colour matching functions that are employed to measure how the human eye perceives colours \citep{CIE2006,IOS2019}. We consider the colour matching functions of the CIE~2006 convention at a step size of \qty{5}{\nm} for the \ang{10} standard observer, i.e., for a field of view of~\ang{10} \citep{CIE2006}. 

Once the colour values $\Lcolour$, $\acolour$ and $\bcolour$ are determined, different colour properties can be computed. In this work, we focus on the perceived colour saturation, measured here by CIE\nbhy $\Lcolour\acolour\bcolour$ chroma~$\Ccolour$, defined as
\begin{align}\label{eq:chroma_function}
    \Ccolour = \sqrt{(\acolour)^2 + (\bcolour)^2}.
\end{align}
In other words, the chroma value determines the distance of a given colour from the $\Lcolour$\nbhy axis in the CIE\nbhy $\Lcolour\acolour\bcolour$ colour space with a fixed lightness value~$\Lcolour$.
\section{Transition Matrix Method (TMM)}\label{section:TMM}
The \ac{tmm}, originally proposed by \citet[]{waterman1965,Waterman1971}, generalizes Mie theory to systems of multiple particles of rather arbitrary shape and is therefore also referred to as the generalized multi-particle Mie method. In contrast to volume discretization methods such as \ac{dda} or \ac{fem}, the particles themselves do not have to be discretized. Instead, the incident, internal and scattered fields are expanded into vector spherical wave functions, and the entire scattering behaviour of a particle is encoded in its transition matrix (T\nbhy matrix), which relates the expansion coefficients of the field exciting the particle to those of its scattered field. Coupling all particles yields a single linear system for the scattering coefficients, from which all optical properties of interest can be computed. In the following, we state this system directly in terms of the particle materials and positions, since these are the quantities that enter the optimization in \cref{section:material_optimization}. The underlying derivation, which is mainly based on the formulations of \citet[]{egel2017,BEUTEL2024109076}, is summarized in \ref{appendix:TMM}.

Let us consider a system of $\Npar\in\NN$ non-overlapping particles. We first fix a wavelength~$\lambda$. Throughout this section, material symbols denote refractive-index values at that wavelength, and the wavelength argument is suppressed. In spectral optimization, each particle instead carries a single material function $m_i\colon\Lambda\to\CC$ with value $m_i(\lambda)$ at the current wavelength. The material of particle~$i\in\{1,\ldots,\Npar\}$ is characterized by its complex refractive index $m_i\in\CC$, and its position by its centre $\rb^i\in\RR^3$. We collect these quantities in the material vector and the position vector
\begin{align}
    \mb = \bigl(m_1,\ldots,m_\Npar\bigr)^\top \in \CC^\Npar
    \quad\text{and}\quad
    \Rb = \bigl(\rb^1,\ldots,\rb^\Npar\bigr)^\top \in \RR^{\Npar\times 3},
\end{align}
i.e., $\Rb_i=(\rb^i)^\top\in\RR^{1\times 3}$. The field expansions associated with particle~$i$ are truncated at the multipolar order $N_\mathrm{t}^i\in\NN$. Each expansion coefficient is addressed by a multi-index $\nu=(\ell,n,m)$, which subsumes the polarization $\ell\in\{1,2\}$ (transverse electric, TE, for $\ell=1$; transverse magnetic, TM, for $\ell=2$), the multipolar index $n\in\{1,\ldots,N_\mathrm{t}^i\}$ and the mode $m\in\{-n,\ldots,n\}$. The ordering of the coefficients is specified in \ref{appendix:TMM_expansions}. Consequently, particle~$i$ carries
\begin{align}\label{eq:TMM_number_coeffs_i}
    M^i = 2N_\mathrm{t}^i\bigl(N_\mathrm{t}^i+2\bigr)
\end{align}
expansion coefficients, and the total number of unknowns of the particulate system is
\begin{align}\label{eq:TMM_number_coeffs}
    \bar{M} = \sum_{i=1}^\Npar M^i.
\end{align}

Three quantities enter the \ac{tmm} system. First, the T\nbhy matrix $\Tf^i(m_i)\in\CC^{M^i\times M^i}$ of particle~$i$ depends only on the particle itself, namely on its size, shape, orientation and material~$m_i$, as well as on the background medium and the wavelength. It is independent of the incident field and of all other particles \citep[Section~5.1]{Mishchenko2002}. For an isotropic sphere, it is diagonal with the Mie coefficients as entries \citep[]{Mishchenko1996}. Second, the coupling matrix $\Cf(\rb^i-\rb^j)\in\CC^{M^i\times M^j}$, $j\neq i$, translates the outgoing waves scattered by particle~$j$ into regular waves about the centre of particle~$i$. It depends only on the relative position $\rb^i-\rb^j$ (and on the background medium and the wavelength) and is defined in \cref{eq:coupling_matrix}. As short form, we write $\Cf^{i\leftarrow j}=\Cf(\rb^i-\rb^j)$. Third, the incident field coefficient vector $\pf^i_\mathrm{inc}(\rb^i)\in\CC^{M^i}$ contains the expansion coefficients of the incident plane wave with respect to the centre~$\rb^i$ and is given analytically by \cref{eq:incident_coefficients_t}. The unknowns are the scattering coefficients $\af^i\in\CC^{M^i}$ of the scattered field of each particle.

\begin{remark}[Invertibility of the T\nbhy Matrix]\label{rem:invertibility_Tmatrix}
    The formulation below uses the inverses $\Tf^i(m_i)^{-1}$. We assume that the truncated single-particle T\nbhy matrices are invertible for every considered material and wavelength. In particular, equality between particle and background refractive indices is excluded, since it gives $\Tf^i=\mathbf{0}$. However, unequal indices alone do not guarantee invertibility because individual Mie coefficients may also vanish. The diagonal entries typically decay at high multipole orders, so that unnecessarily large truncation orders can produce ill-conditioning. We therefore choose only the orders required for convergence.
\end{remark}

The field exciting particle~$i$ is the superposition of the incident field and the fields scattered by all other particles. Together with the T\nbhy matrix relation, this yields the state equation of particle~$i$ (see \ref{appendix:TMM_state_equation}):
\begin{align}\label{eq:TMM_state_equation_i}
    \Tf^i(m_i)^{-1}\af^i - \sum_{\substack{j=1\\j\neq i}}^\Npar \Cf(\rb^i-\rb^j)\,\af^j = \pf^i_\mathrm{inc}(\rb^i),
    \quad i\in\{1,\ldots,\Npar\}.
\end{align}
Combining these equations for all particles leads to the linear system of equations of the \ac{tmm}
\begin{align}\label{eq:TMM_system}
    \Tcalb(\mb,\Rb)\,\Acalb(\mb,\Rb) = \Pcalb_\mathrm{inc}(\Rb)
\end{align}
with the system matrix $\Tcalb(\mb,\Rb)\in\CC^{\bar{M}\times\bar{M}}$ given by
\begin{align}\label{eq:TMM_system_all}
    \Tcalb(\mb,\Rb) = \begin{pmatrix}
        \Tf^1(m_1)^{-1} & -\Cf(\rb^1-\rb^2) & \cdots & -\Cf(\rb^1-\rb^\Npar) \\
        -\Cf(\rb^2-\rb^1) & \Tf^2(m_2)^{-1} & \ddots & \vdots \\
        \vdots & \ddots & \ddots & -\Cf(\rb^{\Npar-1}-\rb^\Npar) \\
        -\Cf(\rb^\Npar-\rb^1) & \cdots & -\Cf(\rb^\Npar-\rb^{\Npar-1}) & \Tf^\Npar(m_\Npar)^{-1}
    \end{pmatrix},
\end{align}
the state vector $\Acalb(\mb,\Rb)\in\CC^{\bar{M}}$ and the right-hand side $\Pcalb_\mathrm{inc}(\Rb)\in\CC^{\bar{M}}$ given by
\begin{align}\label{eq:TMM_sol_vector}
    \Acalb(\mb,\Rb) = \begin{pmatrix} \af^1(\mb,\Rb)\\ \vdots\\ \af^\Npar(\mb,\Rb) \end{pmatrix}
    \quad\text{and}\quad
    \Pcalb_\mathrm{inc}(\Rb) = \begin{pmatrix} \pf^1_\mathrm{inc}(\rb^1)\\ \vdots\\ \pf^\Npar_\mathrm{inc}(\rb^\Npar) \end{pmatrix}.
\end{align}
Throughout this work, we assume that $\Tcalb(\mb,\Rb)$ is invertible at every material assignment and wavelength under consideration, so that $\Acalb(\mb,\Rb)$ is the unique solution of \cref{eq:TMM_system} with material assignment~$\mb$ and position assignment~$\Rb$. In block notation, we write $\Acalb_i=\af^i$, $\Pcalb_{\mathrm{inc},i}=\pf^i_\mathrm{inc}(\rb^i)$ and
\begin{align}
    \Tcalb_{ij}(\mb,\Rb) = \begin{cases}
        \Tf^i(m_i)^{-1}, & \text{if } j=i,\\
        -\Cf(\rb^i-\rb^j), & \text{if } j\neq i,
    \end{cases}
    \qquad i,j\in\{1,\ldots,\Npar\}.
\end{align}
Hence, the diagonal blocks depend only on the material properties and characterize the individual particles, whereas the off-diagonal blocks depend only on the positions and describe the interaction between the particles. Since every particle interacts with every other particle, the system matrix~$\Tcalb$ is in general dense, i.e., the system is fully coupled. For the sake of readability, we omit the dependence of $\Acalb$, $\Pcalb_\mathrm{inc}$ and all optical properties derived from them on the wavelength and on the direction and polarization of the incident field.

\subsection{On the Structure of the TMM System Matrix}\label{section:TMM_structure}
The system matrix~$\Tcalb$ from \cref{eq:TMM_system_all} can be decomposed into a material-dependent block-diagonal and a position-dependent off-block-diagonal part via
\begin{align}\label{eq:TMM_split_matrix}
    \Tcalb(\mb,\Rb) = \Tcalb^\mathrm{mat}(\mb) + \Tcalb^\mathrm{pos}(\Rb).
\end{align}
Thereby, the matrix functions $\Tcalb^\mathrm{mat}\colon \CC^\Npar \to \CC^{\bar{M} \times \bar{M}}$ and $\Tcalb^{\mathrm{pos}}\colon\RR^{\Npar\times 3}\to\CC^{\bar{M} \times \bar{M}}$ are defined by
\begin{align}
    \Tcalb^{\mathrm{mat}}(\mb) &= \sum_{i=1}^\Npar \Bf_i \Tf^i(m_i)^{-1} \Bf_i^\top,\label{eq:TMM_diagonal_part}\\
    \intertext{and}
    \Tcalb^\mathrm{pos}(\Rb) &= -\sum_{i=1}^\Npar \sum_{\substack{j=1\\j\neq i}}^\Npar \Bf_i \Cf(\rb^i-\rb^j) \Bf_j^\top,\label{eq:TMM_offdiagonal_part}
\end{align}
where the matrices $\Bf_i\in\RR^{\bar{M}\times M^i}$ for $i\in\{1,\ldots,\Npar\}$ are given by
\begin{align}\label{eq:B_matrix_TMM}
    \Bf_i = \bigl(\eb^{\bar{M}^{i-1}+1},\eb^{\bar{M}^{i-1}+2},\ldots,\eb^{\bar{M}^i-1},\eb^{\bar{M}^i}\bigr),
\end{align}
where $\eb^i\in\RR^{\bar{M}}$ denotes the $i$\nbhy th unit vector in $\RR^{\bar{M}}$ and
\begin{align}\label{eq:bar_M_i}
    \bar{M}^i = \sum_{j=1}^i M^j, \qquad \bar{M}^0 = 0.
\end{align}
The matrix~$\Bf_i$ thus selects the block of coefficients belonging to particle~$i$. In particular, the columns of the $\Bf_i$ are pairwise distinct unit vectors, so that
\begin{align}\label{eq:B_orthogonality}
    \Bf_j^\top \Bf_i = \begin{cases}
        \one_{M^i}, & \text{if } j=i,\\
        \mathbf{0}\in\RR^{M^j\times M^i}, & \text{if } j\neq i,
    \end{cases}
    \qquad i,j\in\{1,\ldots,\Npar\},
\end{align}
which will be used repeatedly in the following derivations.

The material-dependent component of the \ac{tmm} system matrix $\Tcalb^{\mathrm{mat}}$ is a block-diagonal matrix comprising inverted T\nbhy matrices~$\Tf^i$, $i\in\{1,\ldots,\Npar\}$. In contrast, the off-block-diagonal part $\Tcalb^{\mathrm{pos}}$ does not depend on the particle properties itself but on their positions and has blocks of zeros on the block-diagonal.

\subsection{Optical Properties Based on TMM}\label{section:TMM_optical_properties}
The state vector~$\Acalb(\mb,\Rb)$ can be employed to compute a range of optical properties. Recall from \cref{section:general_electromagnetic_scattering_problem} that the background medium is assumed to be non-absorbing, so that $k_\mathrm{b}\in\RR_{>0}$ and the optical properties defined below are real-valued.

\subsubsection{Extinction Cross-Section}
In the \ac{tmm} framework, the extinction cross-section of an individual constituent in our system of particles can be calculated using the corresponding scattering and incident expansion coefficients (cf. \citet[Section~2.3]{mackowski2011} and \citet[Section~3.6.1]{schebarchov2022}).

\begin{definition}[Individual Extinction Cross-Section]
    Let $\Npar\in\NN$ be the number of particles within the particulate system, and $M^i$ the number of expansion coefficients of the $i$\nbhy th particle defined by \cref{eq:TMM_number_coeffs_i}, $i\in\{1,\ldots,\Npar\}$. Furthermore, let $\mb\in \CC^\Npar$ and $\Rb\in\RR^{\Npar\times 3}$ be given and let $\af^i(\mb,\Rb)\in\CC^{M^i}$ be the vector of scattering expansion coefficients with material assignment~$\mb$ and position assignment~$\Rb$, and $\pf_{\mathrm{inc}}^i(\rb^i)\in\CC^{M^i}$ the corresponding incident field expansion coefficient vector of particle~$i$. Moreover, let $k_\mathrm{b}$ be the wavenumber of the background medium. Then, the \emph{individual extinction cross-section} $\sigma_{\mathrm{ext},i}\colon \CC^\Npar\times \RR^{\Npar\times 3} \to \RR$ of particle~$i$ in the \ac{tmm} setting is defined as 
    \begin{align}
        \sigma_{\mathrm{ext},i}(\mb,\Rb) = -\frac{\pi}{k_\mathrm{b}^2} \Re\biggl(\sum_{\nu=1}^{M^i} \overline{\pf_{\mathrm{inc},\nu}^i(\rb^i)} \af^i_\nu(\mb,\Rb)\biggr),
    \end{align}
    which can be written in compact form as
    \begin{align}
        \sigma_{\mathrm{ext},i}(\mb,\Rb) = -\frac{\pi}{k_\mathrm{b}^2} \Re\Bigl(\bigl(\pf_{\mathrm{inc}}^i(\rb^i)\bigr)^\dagger \af^i(\mb,\Rb)\Bigr).
    \end{align}
\end{definition}
    
The extinction cross-section of the entire cluster can be computed by summing the extinction cross-sections of the individual particles within the system.

\begin{definition}[Cluster Extinction Cross-Section]\label{def:extinction_TMM}
    Let $\Npar\in\NN$ be the number of particles within the particulate system, $\bar{M}$ the corresponding total number of expansion coefficients defined by \cref{eq:TMM_number_coeffs}, and let $\mb\in \CC^\Npar$ and $\Rb\in\RR^{\Npar\times 3}$ be given. Furthermore, let $\Acalb(\mb,\Rb) = \Bigl(\bigl(\af^1(\mb,\Rb)\bigr)^\top,\ldots,\bigl(\af^\Npar(\mb,\Rb)\bigr)^\top\Bigr)^\top\in\CC^{\bar{M}}$ be the state vector comprising the scattering expansion coefficients uniquely solving the \ac{tmm} state equation given by {\crefname{equation}{System}{Systems}\cref{eq:TMM_system}} with material assignment~$\mb$ and position assignment~$\Rb$, and $\Pcalb_\mathrm{inc}(\Rb) = \Bigl(\bigl(\pf_{\mathrm{inc}}^1(\rb^1)\bigr)^\top,\ldots,\bigl(\pf_{\mathrm{inc}}^\Npar(\rb^\Npar)\bigr)^\top\Bigr)^\top\in\CC^{\bar{M}}$ the vector of incident field expansion coefficients. Moreover, let $k_\mathrm{b}$ be the wavenumber of the background medium. Then, the \emph{extinction cross-section} $\sigma_{\mathrm{ext}}\colon \CC^\Npar\times \RR^{\Npar\times 3}\to\RR$ of the entire cluster in the \ac{tmm} setting is defined as 
    \begin{align}
        \sigma_\mathrm{ext}(\mb,\Rb) = \sum_{i=1}^\Npar \sigma_{\mathrm{ext},i}(\mb,\Rb) = -\frac{\pi}{k_\mathrm{b}^2}    \sum_{i=1}^\Npar \Re\Bigl(\bigl(\pf_{\mathrm{inc}}^i(\rb^i)\bigr)^\dagger \af^i(\mb,\Rb)\Bigr),    
    \end{align}
    or, alternatively, in compact notation:
    \begin{align}
        \sigma_\mathrm{ext}(\mb,\Rb) = -\frac{\pi}{k_\mathrm{b}^2}  \Re\bigl(\Pcalb_\mathrm{inc}(\Rb)^\dagger \Acalb(\mb,\Rb)\bigr).
    \end{align}
\end{definition}

\subsubsection{Directional Scattering}
The light intensity scattered in a specific direction~$\ab\in\RR^3$, generally defined in \cref{section:directional_scattering}, can be calculated in the \ac{tmm} framework from the so-called plane wave pattern expansion coefficients (see, e.g., \citet{egel2017}), using the transformation operator~$B_{j,\nu}$ defined in \cref{eq:transformation_operator}.

Therefore, let us denote by the tuple $(a,\theta_\ab,\varphi_\ab)$ the spherical coordinates of the direction $\ab$, where $a$ describes the radial distance, $\varphi_\ab$ the azimuthal angle and $\theta_\ab$ the polar angle. Let the corresponding radial basis vector $\hat{\rb}_\ab\in\RR^3$ be defined by
\begin{align}
    \hat{\rb}_\ab = \begin{pmatrix}
        \sin{\theta_\ab} \cos{\varphi_\ab} \\
        \sin{\theta_\ab} \sin{\varphi_\ab} \\
        \cos{\theta_\ab}
    \end{pmatrix}.
\end{align}

\begin{definition}[Directional Scattering Intensity]\label{def:directional_scattering}
    Let $\Npar\in\NN$ be the number of particles within the particulate system, $\bar{M}$ the corresponding total number of expansion coefficients defined by \cref{eq:TMM_number_coeffs}, and let $\mb\in \CC^\Npar$ and $\Rb\in\RR^{\Npar\times 3}$ be given. Furthermore, let $\Acalb(\mb,\Rb) = \Bigl(\bigl(\af^1(\mb,\Rb)\bigr)^\top,\ldots,\bigl(\af^\Npar(\mb,\Rb)\bigr)^\top\Bigr)^\top\in\CC^{\bar{M}}$ be the state vector comprising the scattering expansion coefficients uniquely solving the \ac{tmm} state equation given by {\crefname{equation}{System}{Systems}\cref{eq:TMM_system}} with material assignment~$\mb$ and position assignment~$\Rb$. Moreover, let $k_\mathrm{b}$ be the wavenumber of the background medium and $\ab\in\RR^3$ the scattering direction of interest. Then, the \emph{directional scattering intensity} $s_\ab\colon \CC^\Npar\times \RR^{\Npar\times 3}\to\RR$ of the entire cluster in the \ac{tmm} setting is defined as 
    \begin{align}
        s_\ab(\mb,\Rb) & = \frac{1}{k_\mathrm{b}^2}\sum_{j=1}^2\Biggl\vert\sum_{i=1}^\Npar \sum_{\ell=1}^2 \sum_{n=1}^{N_\mathrm{t}^i} \sum_{m=-n}^n \exp(\ci m \varphi_\ab) B_{j,\ell n m}(\theta_\ab) \exp\bigl(-\ci k_\mathrm{b} (\hat{\rb}_\ab)^\top \rb^i\bigr) \af^i_{\ell n m}(\mb,\Rb)\Biggr\vert^2\\
        & = \frac{1}{k_\mathrm{b}^2}\sum_{j=1}^2\Biggl\vert\sum_{i=1}^\Npar\sum_{\nu=1}^{M^i} \exp(\ci m \varphi_\ab) B_{j,\nu}(\theta_\ab) \exp\bigl(-\ci k_\mathrm{b} (\hat{\rb}_\ab)^\top \rb^i\bigr) \af^i_\nu(\mb,\Rb)\Biggr\vert^2,
    \end{align}
    where in the second line $m$ denotes the third component of the multi-index $\nu=(\ell,n,m)$.
    Define the matrix $\Ecalb(\ab,\Rb)\in\CC^{\bar{M}\times 2}$ by
    \begin{align}
        \Ecalb(\ab,\Rb) &= \bigl(\Ecalb_1(\ab,\Rb),\Ecalb_2(\ab,\Rb)\bigr)
    \end{align}
    with
    \begin{align}
        \Ecalb_j(\ab,\Rb) &= \bigl(\tilde{\Ef}^1_j(\ab,\rb^1)^\top,\ldots,\tilde{\Ef}^\Npar_j(\ab,\rb^\Npar)^\top\bigr)^\top,\quad j\in\{1,2\}.
    \end{align}
    The maps $\tilde{\Ef}^i_j\colon \RR^3\times\RR^3\to\CC^{M^i}$, $i\in\{1,\ldots,\Npar\}$, are defined component-wise by
    \begin{align}
        \bigl(\tilde{\Ef}^i_{j}(\ab,\rb^i)\bigr)_\nu = \bigl(\tilde{\Ef}^i_{j}(\ab,\rb^i)\bigr)_{\ell n m} = \exp(\ci m \varphi_\ab) B_{j,\ell n m}(\theta_\ab) \exp\bigl(-\ci k_\mathrm{b} (\hat{\rb}_\ab)^\top \rb^i\bigr).
    \end{align}
    Then, the directional scattering intensity can be written in a more compact way as
    \begin{align}
        s_\ab(\mb,\Rb) & = \frac{1}{k_\mathrm{b}^2}\sum_{j=1}^2\Biggl\vert\sum_{i=1}^\Npar \tilde{\Ef}^i_{j}(\ab,\rb^i)^\top\af^i(\mb,\Rb)\Biggr\vert^2\\
        & = \frac{1}{k_\mathrm{b}^2}\sum_{j=1}^2\Bigl\vert \Ecalb_{j}(\ab,\Rb)^\top\Acalb(\mb,\Rb)\Bigr\vert^2\\
        & = \frac{1}{k_\mathrm{b}^2} \bigl\Vert\Ecalb(\ab,\Rb)^\top\Acalb(\mb,\Rb)\bigr\Vert_2^2.\label{eq:directional_scattering}
    \end{align}
\end{definition}

\subsubsection{Monitor Scattering}
The monitor scattering intensity denoted by $s_{\Mcal}$ can be expressed as a weighted sum of the directional scattering intensity, see \cref{section:monitor_scattering}. To this end, we insert the definition of the directional scattering intensity from \cref{eq:directional_scattering} into \cref{eq:directional_scat_Romberg}:
\begin{align}
    s_{\Mcal}(\mb,\Rb) & = \sum_{t=1}^{\tilde{n}} \omega_t s_{\ab^t}(\mb,\Rb) = \frac{1}{k_\mathrm{b}^2} \sum_{t=1}^{\tilde n} \omega_t \bigl\Vert \Ecalb(\ab^t,\Rb)^\top \Acalb(\mb,\Rb)\bigr\Vert^2_2 \\
    & = \frac{1}{k_\mathrm{b}^2} \sum_{t=1}^{\tilde n} \omega_t \bigl(\Ecalb(\ab^t,\Rb)^\top \Acalb(\mb,\Rb)\bigr)^\top\ \overline{\Ecalb(\ab^t,\Rb)^\top \Acalb(\mb,\Rb)}\\
    & = \frac{1}{k_\mathrm{b}^2} \sum_{t=1}^{\tilde n} \omega_t \Acalb(\mb,\Rb)^\top \Ecalb(\ab^t,\Rb) \Ecalb(\ab^t,\Rb)^\dagger \overline{\Acalb(\mb,\Rb)},
\end{align}
where $\mb\in \CC^\Npar$ is the material vector, $\Rb\in\RR^{\Npar\times 3}$ the position vector, and $\omega_t$ and $\ab^t$, $t\in\{1,\ldots,\tilde n\}$, are the Romberg quadrature weights and scattering directions, respectively. Now, we can write the monitor scattering intensity as
\begin{align}
    s_{\Mcal}(\mb,\Rb) = \frac{1}{k_\mathrm{b}^2} \Acalb(\mb,\Rb)^\top \biggl( \sum_{t=1}^{\tilde n} \omega_t  \Ecalb(\ab^t,\Rb) \Ecalb(\ab^t,\Rb)^\dagger \biggr) \overline{\Acalb(\mb,\Rb)}.
\end{align}

Let us define the matrix-valued function $\Fcalb\colon\RR^{3\times \tilde n} \times \RR^{\Npar \times 3} \to \CC^{\bar{M} \times \bar{M}}$ by 
\begin{align}\label{eq:matrix_F}
    \Fcalb(\ab^1,\ldots,\ab^{\tilde n}, \Rb) = \sum_{t=1}^{\tilde n} \omega_t  \Ecalb(\ab^t,\Rb) \Ecalb(\ab^t,\Rb)^\dagger.
\end{align}
Note that the function~$\Fcalb$ is a sum of Hermitian matrices and therefore itself Hermitian, i.e., $\Fcalb^\dagger = \Fcalb$. Using the above derivations, we can define the monitor scattering intensity in the \ac{tmm} setting.

\begin{definition}[Monitor Scattering Intensity]\label{def:monitor_scat}
    Let $\Npar\in\NN$ be the number of particles within the particulate system, $\bar{M}$ the corresponding total number of expansion coefficients defined by \cref{eq:TMM_number_coeffs}, and let $\mb\in \CC^\Npar$ and $\Rb\in\RR^{\Npar\times 3}$ be given. Furthermore, let $\Acalb(\mb,\Rb) = \Bigl(\bigl(\af^1(\mb,\Rb)\bigr)^\top,\ldots,\bigl(\af^\Npar(\mb,\Rb)\bigr)^\top\Bigr)^\top\in\CC^{\bar{M}}$ be the state vector comprising the scattering expansion coefficients uniquely solving the \ac{tmm} state equation given by {\crefname{equation}{System}{Systems}\cref{eq:TMM_system}} with material assignment~$\mb$ and position assignment~$\Rb$. Moreover, let $k_\mathrm{b}$ be the wavenumber of the background medium, $\db_\Mcal$ the position vector of the centre of the monitor, and $\theta_\Mcal\in(0,2\pi]$ the opening angle. Then, the \emph{monitor scattering intensity} $s_{\Mcal}\colon \CC^\Npar\times \RR^{\Npar\times 3}\to\RR$ in the \ac{tmm} framework is defined by
    \begin{align}\label{eq:monitor_scat_function_rewrite}
        s_{\Mcal}(\mb,\Rb) = \frac{1}{k_\mathrm{b}^2} \Acalb(\mb,\Rb)^\top \Fcalb(\ab^1,\ldots,\ab^{\tilde n},\Rb) \overline{\Acalb(\mb,\Rb)}
    \end{align}
    with
    \begin{align}
        \Fcalb(\ab^1,\ldots,\ab^{\tilde n},\Rb) = \sum_{t=1}^{\tilde n} \omega_t  \Ecalb(\ab^t,\Rb) \Ecalb(\ab^t,\Rb)^\dagger
    \end{align}
    where $\Ecalb$ is defined according to \cref{def:directional_scattering}, $\ab^t$ denotes the $t$\nbhy th scattering direction and $\omega_t$ the corresponding Romberg integration weight, $t\in\{1,\ldots,\tilde n\}$. The monitor enters this definition through the quadrature nodes and weights.
\end{definition}

Note that both the directional and monitor scattering intensity are quadratic in the state variable, i.e., the scattering expansion coefficients, $\Acalb$.

Once the \ac{tmm} system has been solved, i.e., once the scattering expansion coefficient vector $\Acalb$ has been determined, it is possible to compute the scattering in any direction and in any monitor of interest without solving the system again.

\section{Material Optimization by Sequential Global Programming}\label{section:material_optimization}
We optimize the material assignment of a fixed assembly of particles using the \ac{sgp} method \citep{semmler2018,Nees2022}. The \ac{tmm} supplies the fully coupled scattering state for the current design. The material updates are determined by separable models of the resulting optical objective. We recall the necessary published \ac{sgp} background before presenting the \ac{tmm}\nbhy specific construction in \cref{section:TMM_SGP_models_general}.

\subsection{Design Optimization Problem}\label{section:general_problem}
Let $\lambda_1,\ldots,\lambda_Q\in\Lambda$ be the $Q$ wavelengths used for the prescribed spectral quadrature. Note that a monochromatic problem corresponds to $Q=1$. We represent a particle material by its sampled complex refractive indices,
\begin{align}\label{eq:spectral_material_vector}
    m_i=\bigl(m_i(\lambda_1),\ldots,m_i(\lambda_Q)\bigr)\in\MatSpace:=\CC^Q,
    \qquad \mb=(m_1,\ldots,m_\Npar)\in\MatSpace^\Npar.
\end{align}
Thus one particle is one material coordinate even though it has $Q$ spectral components. The incident fields and fully coupled \ac{tmm} state equations are evaluated at these nodes (and, where necessary, for both polarizations).

For the bi-material optimization problems considered here, the prescribed catalogue consists of two distinct sampled material vectors,
\begin{align}\label{eq:material_catalogue}
    \mathcal{U}=\{\nbar_1,\nbar_2\}\subset\MatSpace,\qquad \nbar_1\neq\nbar_2.
\end{align}
Here $\nbar_j=(\nbar_j(\lambda_q))_{q=1}^{Q}$ is the sampled vector of the prescribed material dispersion, $j=1,2$, and the vectors are assumed distinct at the chosen nodes. The discrete assignment set is $\mathcal U_\mathrm D=\mathcal U^\Npar$. At each $\lambda_q$, the fully coupled \ac{tmm} state satisfies \cref{eq:TMM_system}. The material assignment problem can therefore be written as
\begin{subequations}\label{eq:general_mat_problem_TMM}
\begin{align}
    \min_{\mb\in\mathcal{U}_\mathrm D}\ & J(\mb)\label{eq:general_mat_problem_TMM_obj_func}\\
    \text{s.t.}\ & \Tcalb(\mb,\lambda_q)\Acalb(\mb,\lambda_q)=\Pcalb_\mathrm{inc}(\lambda_q),\quad q=1,\ldots,Q.\label{eq:general_mat_problem_TMM_system}
\end{align}
\end{subequations}
Here $J$ denotes a real-valued optical objective, with any spectral integral represented by the chosen fixed quadrature. The particle positions are fixed and suppressed in the notation.

To obtain a continuous relaxation, we parameterize the catalogue by one scalar design coordinate per particle, following \citet{Nees2022}.
\begin{definition}[Parameterization of the Material Catalogue]\label{def:parameterization_of_catalogue}
    A mapping $P\colon[0,1]\to\MatSpace$ is a parameterization of $\mathcal U$ if it is twice continuously differentiable with respect to the real scalar $\rho$, satisfies $P(0)=\nbar_1$ and $P(1)=\nbar_2$, and is injective.
\end{definition}
The corresponding admissible set is
\begin{align}\label{eq:material_admissible_set}
    \Uadmat=\bigl\{P(\rho)\mid\rho\in[0,1]\bigr\}^\Npar.
\end{align}
We use the linear refractive-index interpolation
\begin{align}\label{eq:linear_parameterization}
    P(\rho)=\bigl((1-\rho)\nbar_1(\lambda_q)+\rho\nbar_2(\lambda_q)\bigr)_{q=1}^Q.
\end{align}
The same scalar $\rho_i$ determines all spectral components of particle $i$: there are exactly $\Npar$ scalar design parameters. Because the sampled endpoint vectors are distinct, this interpolation is injective. Intermediate refractive indices are a numerical interpolation between the prescribed materials. Interpreting $\rho$ as a volume fraction would require a physically justified mixing model, see, e.g., the absorber-material discussion in \citet{bleyer2025}. Intermediate (grey) indices can be discouraged by adding a greyness-penalty term, as in \citet{Nees2022}. This serves a purpose similar to the suppression of intermediate materials in solid isotropic material with penalization (SIMP) approaches \citep{bendsoe2004}, although the penalization mechanisms are distinct. No greyness penalty is used in the present computations.

The relaxed optimization problem becomes
\begin{subequations}\label{eq:general_relaxed_mat_problem_TMM}
\begin{align+}
    \min_{\mb\in\Uadmat}\ & J(\mb)\label{eq:general_relaxed_mat_problem_TMM_obj_func}\\
    \text{s.t.}\ & \Tcalb(\mb,\lambda_q)\Acalb(\mb,\lambda_q)=\Pcalb_\mathrm{inc}(\lambda_q),\quad q=1,\ldots,Q.\label{eq:general_relaxed_mat_problem_TMM_system}
\end{align+}
\end{subequations}
For later reference, the single-particle admissible set is
\begin{align}\label{eq:admissible_set_single_element}
    \mathcal V_\mathrm{ad}=\{P(\rho)\mid\rho\in[0,1]\},\qquad \Uadmat=\mathcal V_\mathrm{ad}^\Npar.
\end{align}
We assume that the individual truncated T\nbhy matrices and the fully coupled \ac{tmm} systems are invertible for the admissible material assignments at every sampled wavelength. Thus the \ac{tmm} state and its observables are well defined on the product set $\Uadmat$ and extend smoothly to a regular open neighbourhood $\mathcal D\subset\MatSpace^\Npar$ of that set, under the local T\nbhy matrix regularity specified in \cref{section:TMM_SGP_models_general}. Smoothness of a subsequently chosen objective is assumed only at the reference values where derivatives are used (in particular, chroma requires nonzero chroma at the reference).

\subsection{Sequential Global Programming}\label{section:SGP}
We minimize $J$ over $\Uadmat$ using \ac{sgp}, which successively replaces $J$ by separable approximations and solves independent one-particle subproblems. Each block coordinate is the sampled material vector of one particle. The bi-material problems are parameterized by $\rho_i\in[0,1]$.

\subsubsection{General Separable Approximations}\label{section:general_separable_approximations}
We recall the following notions from \citet{Nees2022,Gangl2024}, now stated for finite-dimensional particle-wise material vectors. Functions are understood on the admissible product set, and derivatives on a regular open neighbourhood of the reference design.
\begin{definition}[Separable Function]\label{def:separable_function}
    A function $S\colon\Uadmat\to\RR$ is \emph{separable} if there are functions $S_i\colon\mathcal V_\mathrm{ad}\to\RR$ and a constant $S_0\in\RR$ such that
    \begin{align}
        S(\mb)=S_0+\sum_{i=1}^\Npar S_i(m_i).
    \end{align}
\end{definition}
\begin{definition}[Exact Separable Model]\label{def:exact_separable_model}
    Given $f$ defined on a regular neighbourhood of $\Uadmat$ and $\tilde{\mb}\in\Uadmat$, a separable $S[f,\tilde{\mb}]$ is an \emph{exact separable model} of $f$ at $\tilde\mb$, i.e., it is exact coordinate-wise, if
    \begin{align}\label{eq:prop_exact_sep_model}
      S[f,\tilde\mb](\tilde\mb+h_i\eb^i)=f(\tilde\mb+h_i\eb^i)
    \end{align}
    for every particle $i$ and every $h_i\in\MatSpace$ with $\tilde m_i+h_i\in\mathcal V_\mathrm{ad}$. Here $h_i\eb^i$ modifies only the sampled material vector of particle $i$.
\end{definition}
\begin{definition}[First-Order Separable Approximation]\label{def:first_order_separable_model}
    Suppose $f$ and $S[f,\tilde\mb]$ admit continuously \emph{real differentiable extensions} to an open neighbourhood of $\tilde\mb\in\Uadmat$. The separable model is \emph{first-order correct} at $\tilde\mb$ if
    \begin{align}
      S[f,\tilde\mb](\tilde\mb)=f(\tilde\mb),\qquad
      DS[f,\tilde\mb](\tilde\mb)[\mathbf{h}]=Df(\tilde\mb)[\mathbf{h}]
    \end{align}
    for all $\mathbf{h}\in\MatSpace^\Npar\simeq\RR^{2Q\Npar}$.
\end{definition}
For differentiable objectives, exactness along admissible particle-wise variations implies first-order agreement in the corresponding admissible directions. The particular exact models constructed below admit smooth extensions for which first-order agreement holds in all ambient directions.
\begin{theorem}\label{thm:exact_separable_model}
    For $f$ defined on a regular neighbourhood of $\Uadmat$ and $\tilde\mb\in\Uadmat$, set
    \begin{align}
      S_i[f,\tilde\mb](u)=f\bigl(\tilde\mb+(u-\tilde m_i)\eb^i\bigr),\qquad u\in\mathcal V_\mathrm{ad}.
    \end{align}
    Then
    \begin{align}
      S[f,\tilde\mb](\mb)=(1-\Npar)f(\tilde\mb)+\sum_{i=1}^\Npar S_i[f,\tilde\mb](m_i)
    \end{align}
    is an exact separable model of $f$ at $\tilde\mb$.
\end{theorem}
\begin{proof}
    If only coordinate $i$ changes, the other $\Npar-1$ terms equal $f(\tilde\mb)$, cancelling the prefactor and the remaining term is $f(\tilde\mb+h_i\eb^i)$.
\end{proof}
\begin{definition}[Coordinate-wise Minimum]\label{def:coordinatewise_minimum}
    A design $\mb^*\in\Uadmat$ is a \emph{coordinate-wise minimum} of $J$ if for every $i\in\{1,\ldots,\Npar\}$ and $u\in\mathcal V_\mathrm{ad}$
    \begin{align}\label{eq:separable_local_minimum}
        J(\mb^*)\leq J\bigl(\mb^*+(u-m_i^*)\eb^i\bigr).
    \end{align}
    
\end{definition}
This is also called a \emph{one-particle} or \emph{one-coordinate minimum}. Unlike first-order stationarity, it requires neither derivatives nor continuous design variables and therefore also applies to finite catalogues. Equivalently, it is a pure-strategy Nash equilibrium of the identical-interest game in which each particle chooses its material and all particles minimize $J$. For the differentiable reduced objective $j(\boldsymbol\rho)=J\bigl(P(\rho_1),\ldots,P(\rho_\Npar)\bigr)$, \cref{eq:separable_local_minimum} implies the box-constrained KKT conditions, since $\partial_i j(\boldsymbol\rho^*)(t-\rho_i^*)\geq0$ for every $t\in[0,1]$. The converse is false in general: coordinate-wise minimality also need not imply a local minimum with respect to simultaneous changes of several particles. If $\mb^*$ globally minimizes an exact separable model constructed at $\mb^*$, then \cref{eq:prop_exact_sep_model} immediately gives \cref{eq:separable_local_minimum}.

Constructing this exact model by separate full-state solves at every trial material would be expensive. The \ac{tmm} structure instead permits inexpensive approximate coordinate-wise evaluations.

We use the Sherman--Morrison--Woodbury identity \citep[Section~2.1.4]{golub2013matrix}: for an invertible matrix $\Ab$, matrices $\Ub$ and $\Vb$ of compatible size with inner dimension $k\in\NN$, and invertible $\one_k+\Vb\Ab^{-1}\Ub$,
\begin{align}\label{eq:SMW_formula}
    (\Ab+\Ub\Vb)^{-1}=\Ab^{-1}-\Ab^{-1}\Ub(\one_k+\Vb\Ab^{-1}\Ub)^{-1}\Vb\Ab^{-1}.
\end{align}

\subsubsection{SGP Algorithm}\label{section:SGP_algorithm}
Let us consider the optimization problem presented by {\crefname{equation}{Problem}{Problems}\cref{eq:general_relaxed_mat_problem_TMM}} with a real-valued objective function~$J$ and an admissible set~$\Uadmat$. Furthermore, let $S[J,\tilde\mb](\cdot)$ be a general separable approximation of the objective function~$J$ at the point of approximation $\tilde{\mb} \in \Uadmat$. Then, the \ac{sgp} algorithm in the context of this work is given by \cref{alg:SGP_alg} (cf. also \citet[Section~4]{semmler2018}, \citet[Section~3]{Nees2022} and \citet[Section~3.3]{vu2023}).

The model and the proximal term in \cref{alg:SGP_alg} are both separable across particles. The inner loop increases the proximal parameter until the full coupled objective satisfies the descent test. The parameter is then reduced for the next outer iteration. Each particle-wise model is minimized over the same admissible material curve, and the numerical implementation selects the best point on a prescribed finite grid of $\rho\in[0,1]$. Thus ``global'' solution of a scalar subproblem here refers to its sampled candidate grid, not a certificate for the entire continuous interval.

\begin{algorithm}[!htb]
    \caption{Sequential Global Programming (SGP) Algorithm}
    \label{alg:SGP_alg}
    \begin{algorithmic}[1]
        \REQUIRE $n=0,\ \mb^0\in\Uadmat,\ \tau>0,\ p\geq 1,\ \overline{\vartheta} > 1,\ \underline{\vartheta} \in (0,1),\ \varepsilon_1 \geq 0,\ \varepsilon_2 \geq 0$
        \STATE $j^0 \leftarrow J(\mb^0)$
        \REPEAT
        \REPEAT
        \STATE $\displaystyle{\mb^\ast \leftarrow \argmin_{\mb\in\Uadmat}\ S[J,\mb^n](\mb) + \tau \Vert\mb^n-\mb\Vert_p^p}$ \algorithmiccomment{Solution of the subproblem}
        \IF{$\mb^\ast=\mb^n$} \STATE \textbf{terminate} \ENDIF
        \STATE $\tau \leftarrow \overline{\vartheta} \tau$
        \STATE $j^{n+1} \leftarrow J(\mb^\ast)$
        \UNTIL $j^{n+1} < j^n - \varepsilon_2 \Vert\mb^\ast - \mb^n\Vert_2^2$
        \STATE $\tau \leftarrow \underline{\vartheta} \tau$
        \STATE $\mb^{n+1} \leftarrow \mb^\ast$
        \STATE $n \leftarrow n+1$
        \UNTIL $\Vert\mb^n - \mb^{n-1}\Vert_2^2 < \varepsilon_1$
    \end{algorithmic}
\end{algorithm}

The zero-update termination in \cref{alg:SGP_alg} avoids applying a strict descent test to an unchanged design. The proximal and descent terms are evaluated in the scalar design parameters $\rho_i\in[0,1]$. Write $\boldsymbol{\rho}=(\rho_1,\ldots,\rho_\Npar)$ and $\mb=(P(\rho_1),\ldots,P(\rho_\Npar))$ for the corresponding material assignment. In particular, the proximal term used for the numerical computations is
\begin{align}\label{eq:separable_prox_point_term}
   \|\mb^n-\mb\|_p^p\;\mathrel{:=}\;\|\boldsymbol{\rho}^n-\boldsymbol{\rho}\|_p^p\;=\sum_{i=1}^\Npar|\rho_i^n-\rho_i|^p,
\end{align}
and analogously $\|\mb^n-\mb\|_2^2:=\sum_i|\rho_i^n-\rho_i|^2$ in the acceptance and stopping tests.

In the bi-material setting considered here, the subproblem in the inner loop of \cref{alg:SGP_alg} decouples into $\Npar$ univariate minimization problems over the single admissible set~$\normalfont \mathcal{V}_\mathrm{ad}$ from \cref{eq:admissible_set_single_element}, i.e., over the scalar parameter $\rho\in[0,1]$ in \cref{eq:linear_parameterization}. In the idealized algorithm, each of them is solved globally on the continuous interval. In our implementation, the best candidate on a fine finite grid of $\rho$\nbhy values is selected for each particle. The joint combinatorial problem over all particles is not solved. For the general multi-material treatment, see \citet{Nees2022}.

Convergence and globalization results for \ac{sgp} are discussed in \citet{semmler2018}. They do not by themselves guarantee that the present iteration with approximate models converges to a coordinate-wise minimum in the sense of \cref{def:coordinatewise_minimum}.

\begin{remark}[Simultaneous Updates and Coordinate-wise Optimality]\label{rem:coordinatewise_globalization}
    For exact separable models and exact particle-wise minimization, coordinate-wise optimality can be made compatible with simultaneous \ac{sgp} updates through a possible additional globalization mechanism: first test the full set of proposed material changes. If insufficient descent is observed, rank the changes by their predicted individual improvements and successively keep only the better half. In the limiting case only the best single-particle change remains, for which the exact separable model agrees with the full objective. A short convergence argument for this hypothetical reduction strategy, without a proximal term, is given in \ref{appendix:coordinatewise_globalization}. The strategy is not implemented in \cref{alg:SGP_alg}.
\end{remark}

In the following, a general separable model of the objective function~$J$ from {\crefname{equation}{Problem}{Problems}\cref{eq:general_relaxed_mat_problem_TMM}} at the approximation point $\tilde\mb\in\Uadmat$ is represented by the function~$S[J,\tilde\mb]$. In this context, an exact separable model according to \cref{def:exact_separable_model} is denoted by $\Sexact[J,\tilde\mb]$. Furthermore, a general first-order separable approximation with respect to \cref{def:first_order_separable_model} is denoted by $\Sapprox[J,\tilde\mb]$, and the corresponding computationally accessible choice based on the block-diagonal approximation of the inverse system matrix by $\Sdiag[J,\tilde\mb]$.

\section{SGP Separable Models Based on TMM}\label{section:TMM_SGP_models_general}
We derive exact separable models and computationally accessible first-order approximations for objectives depending on finitely many linear observables of the fully coupled \ac{tmm} state. The resulting formulas include polynomial objectives as special cases and also apply to more general differentiable optical properties.

Fix a sampled wavelength $\lambda_q$ and suppress the index $q$ in the pointwise derivation. The particle material $m_i\in\MatSpace=\CC^Q$ enters this state equation only through its $q$\nbhy th component $m_{i,q}=m_i(\lambda_q)$. We write $\Tf^i(m_i)$ for $\Tf^i(m_{i,q},\lambda_q)$, and likewise suppress $\lambda_q$ in the system matrix and state. The extension to the jointly affected states is given in \cref{cor:several_state_equations}.

For $g\in\mathcal C^1(\CC^L,\RR)$ and $\Lcalb\in\CC^{\bar M\times L}$ consider
\begin{align}\label{eq:general_obj_func_TMM}
    J(\mb)=g\bigl(\Lcalb^\top\Acalb(\mb)\bigr),\qquad \mb\in\mathcal D,
\end{align}
where $\Acalb(\mb)$ uniquely solves the fully coupled state equation
\begin{align}\label{eq:TMM_system_SGP}
    \Tcalb(\mb)\Acalb(\mb)=\Pcalb_\mathrm{inc}.
\end{align}
The $L$ entries of $\Lcalb^\top\Acalb(\mb)$ are complex linear scattering observables.

Note that the function~$g$ is real-valued and, therefore, not complex differentiable unless it is constant. Consequently, the differentiability of~$g$ has to be understood in the real sense, i.e., by identifying $\CC^L$ with $\RR^{2L}$ via $\mathbf{w}_j = x_j+\ci y_j$, $j\in\{1,\ldots,L\}$. Let us collect the partial derivatives of~$g$ in the complex vector $\nabla g(\mathbf{w})\in\CC^L$ with the entries
\begin{align}\label{eq:complex_gradient}
    \bigl(\nabla g(\mathbf{w})\bigr)_j = \frac{\partial g}{\partial x_j}(\mathbf{w}) + \ci\,\frac{\partial g}{\partial y_j}(\mathbf{w}), \qquad j\in\{1,\ldots,L\},
\end{align}
for $\mathbf{w}\in\CC^L$, so that the directional derivative of~$g$ in a direction $\boldsymbol{\delta}\in\CC^L$ is given by
\begin{align}
    \lim_{\substack{\epsilon\to 0\\ \epsilon\in\RR}} \frac{1}{\epsilon}\bigl(g(\mathbf{w}+\epsilon\boldsymbol{\delta})-g(\mathbf{w})\bigr) = \Re\bigl(\nabla g(\mathbf{w})^\dagger\boldsymbol{\delta}\bigr).
\end{align}
We assume that the inverse single-particle T\nbhy matrices $z\mapsto\Tf^i(z,\lambda_q)^{-1}$ are holomorphic on neighbourhoods of the relevant sampled refractive indices and that the coupled matrix is invertible on the regular neighbourhood $\mathcal D$ introduced in \cref{section:general_problem}. Let $\partial_z\Tf^i(z,\lambda_q)^{-1}$ denote the entrywise complex derivative. The dependence on the finite-dimensional material vector is therefore (real) differentiable, and the real chain rule for a variation $h_i\in\MatSpace$ reads
\begin{align}\label{eq:chain_rule_g}
 D_{m_i}J(\mb)[h_i]=\Re\Bigl(\nabla g(\Lcalb^\top\Acalb(\mb))^\dagger\,D_{m_i}(\Lcalb^\top\Acalb)(\mb)[h_i]\Bigr).
\end{align}
For the fixed state equation, the direction $h_i$ enters the T\nbhy matrix derivative only through $h_{i,q}$. The required state derivative is stated explicitly in \cref{lem:single_element_update}. For first-order agreement at a reference material, it is sufficient that $g$ be $\mathcal C^1$ in a neighbourhood of the corresponding observables. Its global $\mathcal C^1$ regularity is a convenient stronger assumption.

Before deriving an exact separable model of the aforementioned objective function, let us introduce the notation
\begin{align}\label{eq:delta_T}
    \Delta \Tf^i(z_1,z_2) = \Tf^i(z_1)^{-1} - \Tf^i(z_2)^{-1}
\end{align}
for $z_1,z_2\in\CC$, with the wavelength fixed. When vector-valued materials are inserted, both arguments are evaluated at the fixed component $q$. Thus $\Delta\Tf^i(\tilde m_i,m_i)$ abbreviates $\Delta\Tf^i(\tilde m_{i,q},m_{i,q})$ for the present state.

Let $\tilde\mb,\mb\in\Uadmat$ and $\Kcalb\in\CC^{\bar M\times\bar M}$. Whenever the indicated inverse exists, set
\begin{align}\label{eq:X_matrix}
 \Xf^i[\tilde\mb,\Kcalb](m_i)=\bigl(\one_{M^i}-\Delta\Tf^i(\tilde m_i,m_i)\Bf_i^\top\Kcalb\Bf_i\bigr)^{-1}\Delta\Tf^i(\tilde m_i,m_i).
\end{align}
The auxiliary matrix $\Kcalb$ approximates $\Tcalb(\tilde\mb)^{-1}$ only in this surrogate update. By \cref{thm:exact_separable_model}, an exact separable model has the generic form $(1-\Npar)J(\tilde\mb)+\sum_i J\bigl(\tilde\mb+(m_i-\tilde m_i)\eb^i\bigr)$. \cref{lem:single_element_update} computes the coordinate-wise observables without a new full state solve for every candidate.

\begin{lemma}[Single-Element Update of the Linear Functionals]\label{lem:single_element_update}
    Let $\Lcalb\in\CC^{\bar{M}\times L}$ and $\tilde{\mb}\in \Uadmat$ be given, and let $\Acalb(\tilde{\mb})\in\CC^{\bar{M}}$ be the unique solution of the \ac{tmm} state system given by \cref{eq:TMM_system_SGP} with material assignment $\tilde \mb$. Moreover, let $\Qcalb(\tilde{\mb})\in\CC^{\bar{M}\times L}$ uniquely solve the adjoint equation
    \begin{align}\label{eq:general_adjoint_equation_TMM}
        \Tcalb(\tilde{\mb})^\top\Qcalb = \Lcalb
    \end{align}
    with material assignment~$\tilde{\mb}$. Furthermore, for $i\in\{1,\ldots,\Npar\}$ and $\Kcalb\in\CC^{\bar{M}\times\bar{M}}$ such that the matrix $\Xf^i[\tilde{\mb},\Kcalb](m_i)$ from \cref{eq:X_matrix} is well defined, let
    \begin{align}\label{eq:zeta_scalar}
        \boldsymbol{\zeta}^i[\tilde{\mb},\Kcalb](m_i) = \Qcalb(\tilde{\mb})^\top\Bf_i\Xf^i[\tilde{\mb},\Kcalb](m_i)\Bf_i^\top\Acalb(\tilde{\mb}) \in\CC^L.
    \end{align}
    Then, the following statements hold:
    \begin{enumerate}
        \item\label{lem:single_element_update_exact} For $\Kcalb = \Tcalb(\tilde{\mb})^{-1}$, the update of the linear functionals is exact, i.e.,
        \begin{align}\label{eq:single_element_update_exact}
            \Lcalb^\top\Acalb\bigl(\tilde{\mb}+(m_i - \tilde m_i)\eb^i\bigr) = \Lcalb^\top\Acalb(\tilde{\mb}) + \boldsymbol{\zeta}^i\bigl[\tilde{\mb},\Tcalb(\tilde{\mb})^{-1}\bigr](m_i).
        \end{align}
        \item\label{lem:single_element_update_first_order} For an arbitrary matrix~$\Kcalb$, the update is first-order correct, i.e., $\boldsymbol{\zeta}^i[\tilde{\mb},\Kcalb](\tilde m_i) = \mathbf{0}$ and
        \begin{align}\label{eq:single_element_update_first_order}
            D_{m_i}\boldsymbol{\zeta}^i[\tilde\mb,\Kcalb](\tilde m_i)[h_i] =D_{m_i}(\Lcalb^\top\Acalb)(\tilde\mb)[h_i]
            =-\Qcalb(\tilde\mb)^\top\Bf_i\,\partial_z\Tf^i(z,\lambda_q)^{-1}\big|_{z=\tilde m_{i,q}}\,\Bf_i^\top\Acalb(\tilde\mb)\,h_{i,q},\quad h_i\in\MatSpace.
        \end{align}
    \end{enumerate}
\end{lemma}

\begin{proof}
For a change of the $i$\nbhy th sampled material vector, the fully coupled matrix at a fixed wavelength changes only in the $i$\nbhy th diagonal block:
\begin{align}\label{eq:single_update_low_rank}
  \Tcalb\bigl(\tilde\mb+(m_i-\tilde m_i)\eb^i\bigr)
  =\Tcalb(\tilde\mb)-\Bf_i\Delta\Tf^i(\tilde m_i,m_i)\Bf_i^\top.
\end{align}
Applying \cref{eq:SMW_formula} and using $\Lcalb^\top\Tcalb(\tilde\mb)^{-1}=\Qcalb(\tilde\mb)^\top$ gives \cref{eq:single_element_update_exact}.
At the reference material, $\Delta\Tf^i(\tilde m_i,\tilde m_i)=0$, so $\Xf^i[\tilde\mb,\Kcalb](\tilde m_i)=0$ for every admissible $\Kcalb$. Differentiating the state equation and the update at this point gives the same directional derivative, since all terms depending on $\Kcalb$ multiply $\Delta\Tf^i$ and vanish at the reference point. Thus \cref{eq:single_element_update_first_order} follows. Note that the full algebra is given in \ref{appendix:single_particle_update}.
\end{proof}

By inserting the exact update from \cref{lem:single_element_update} into the function~$g$, it is possible to state the desired exact separable model.

\begin{theorem}[Exact Separable Model for \cref{eq:general_obj_func_TMM}]\label{thm:exact_sep_model_TMM}
    Let $\Lcalb\in\CC^{\bar{M}\times L}$ and $\tilde{\mb}\in \Uadmat$ be given, and let $\Acalb(\tilde{\mb})$ and $\Qcalb(\tilde{\mb})$ be the unique solutions of the state and the adjoint equation given by \cref{eq:TMM_system_SGP,eq:general_adjoint_equation_TMM}, respectively. Furthermore, let the matrix $\Xf^i\bigl[\tilde{\mb},\Tcalb(\tilde{\mb})^{-1}\bigr](m_i)$ from \cref{eq:X_matrix} be well defined for all $i\in\{1,\ldots,\Npar\}$ and all considered material assignments $\mb\in \MatSpace^\Npar$. Then, the function $\Sexact[J,\tilde{\mb}](\cdot)\colon\Uadmat \to \RR$ defined by
    \begin{align}\label{eq:general_exact_sep_model_TMM}
        \Sexact[J,\tilde{\mb}](\mb) = J(\tilde{\mb})(1-\Npar) + \sum_{i=1}^\Npar g\Bigl(\Lcalb^\top \Acalb(\tilde{\mb}) + \boldsymbol{\zeta}^i\bigl[\tilde{\mb},\Tcalb(\tilde{\mb})^{-1}\bigr](m_i)\Bigr),
    \end{align}
    where $\boldsymbol{\zeta}^i$ is given by \cref{eq:zeta_scalar}, is an exact separable model according to \cref{def:exact_separable_model} of the objective function~$J$ from \cref{eq:general_obj_func_TMM}.
\end{theorem}

\begin{proof}
    We have to show that the function~$\Sexact[J,\tilde{\mb}]$ defined in \cref{eq:general_exact_sep_model_TMM} coincides with the exact separable model from \cref{thm:exact_separable_model}. According to \cref{eq:single_element_update_exact}, the $i$\nbhy th summand in \cref{eq:general_exact_sep_model_TMM} is given by
    \begin{align}
        g\Bigl(\Lcalb^\top \Acalb(\tilde{\mb}) + \boldsymbol{\zeta}^i\bigl[\tilde{\mb},\Tcalb(\tilde{\mb})^{-1}\bigr](m_i)\Bigr) = g\Bigl(\Lcalb^\top\Acalb\bigl(\tilde{\mb}+(m_i - \tilde m_i)\eb^i\bigr)\Bigr) = J\bigl(\tilde{\mb}+(m_i - \tilde m_i)\eb^i\bigr).
    \end{align}
    Therefore, both functions coincide, which concludes the proof.
\end{proof}

Similarly to the \ac{dda} framework in \citet[]{Nees2022}, the computation of the inverse \ac{tmm} system matrix, $\Tcalb(\tilde{\mb})^{-1}$, and consequently the exact separable model can be computationally prohibitive. Nevertheless, by choosing an arbitrary $\bar{M}\times\bar{M}$~matrix for the inverse system matrix in the expression for the exact separable model from \cref{eq:general_exact_sep_model_TMM}, a first-order correct separable approximation can be obtained.

\begin{theorem}[First-Order Separable Approximation for \cref{eq:general_obj_func_TMM}]\label{thm:general_first_order_separable_approx_TMM}
    Let $\Lcalb\in\CC^{\bar{M}\times L}$ and $\tilde{\mb}\in \Uadmat$ be given, and let $\Acalb(\tilde{\mb})\in\CC^{\bar{M}}$ and $\Qcalb(\tilde{\mb})\in\CC^{\bar{M}\times L}$ be the unique solutions of the \ac{tmm} state and adjoint equation given by \cref{eq:TMM_system_SGP,eq:general_adjoint_equation_TMM}, respectively, with material assignment $\tilde{\mb}$. Furthermore, let $\Kcalb\in\CC^{\bar{M}\times\bar{M}}$ be arbitrary such that the matrix $\Xf^i[\tilde{\mb},\Kcalb](m_i)$ from \cref{eq:X_matrix} is well defined for all $i\in\{1,\ldots,\Npar\}$ and all considered material assignments $\mb\in \MatSpace^\Npar$. Then, the function $\Sapprox[J,\tilde{\mb}](\cdot)\colon \Uadmat \to \RR$ defined by
    \begin{align}\label{eq:general_approx_sep_model_TMM}
        \Sapprox[J,\tilde{\mb}](\mb) = J(\tilde{\mb})(1-\Npar) + \sum_{i=1}^\Npar g\bigl(\Lcalb^\top \Acalb(\tilde{\mb}) + \boldsymbol{\zeta}^i[\tilde{\mb},\Kcalb](m_i)\bigr),
    \end{align}
    where $\boldsymbol{\zeta}^i$ is given by \cref{eq:zeta_scalar}, is a first-order separable approximation according to \cref{def:first_order_separable_model} of the objective function~$J$ from \cref{eq:general_obj_func_TMM} in the point of approximation~$\tilde{\mb}$.
\end{theorem}

\begin{proof}
By \cref{lem:single_element_update}, all update terms vanish at the approximation point, giving $\Sapprox[J,\tilde\mb](\tilde\mb)=J(\tilde\mb)$. For each particle variation $h_i\in\MatSpace$, the same lemma and the real chain rule \cref{eq:chain_rule_g} give
\begin{align}
 D_{m_i}\Sapprox[J,\tilde\mb](\tilde\mb)[h_i]
 =\Re\Bigl(\nabla g(\Lcalb^\top\Acalb(\tilde\mb))^\dagger D_{m_i}(\Lcalb^\top\Acalb)(\tilde\mb)[h_i]\Bigr)
 =D_{m_i}J(\tilde\mb)[h_i].
\end{align}
Since the model is separable, these equalities for every coordinate imply equality of the full real Fr\'echet derivatives.
\end{proof}

Thus any matrix~$\Kcalb$ yields a first-order separable approximation. As in the \ac{dda} framework of \citet[]{Nees2022}, $\Kcalb$ should nevertheless be chosen such that the exactness property is retained approximately. We therefore use the following block-diagonal approximation of the inverse \ac{tmm} system matrix.

\begin{definition}[Block-Diagonal Approximation of Inverse \ac{tmm} System Matrix]\label{def:diagonal_approximation_TMM}
    Let $\mb\in \Uadmat$ be a given design vector, and let $\Tcalb(\mb)\in\CC^{\bar{M}\times\bar{M}}$ be the \ac{tmm} system matrix defined in \cref{section:TMM} for the material assignment~$\mb$. Furthermore, let $\Tcalb^\mathrm{mat}(\mb)\in\CC^{\bar{M}\times\bar{M}}$ be the corresponding material-dependent block-diagonal matrix (see \cref{eq:TMM_diagonal_part}). Then, we call $\Tcalb^\mathrm{mat}(\mb)^{-1}$ the \emph{block-diagonal approximation} of the original inverse \ac{tmm} system matrix, $\Tcalb(\mb)^{-1}$, i.e.,
    \begin{align}
        \Tcalb(\mb)^{-1} \approx \Tcalb^\mathrm{mat}(\mb)^{-1}.
    \end{align} 
\end{definition}

Due to the block-diagonal structure of $\Tcalb^{\mathrm{mat}}$, it holds
\begin{align}
    \Tcalb^{\mathrm{mat}}(\mb)^{-1} &= \sum_{i=1}^\Npar \Bf_i \Tf^i(m_i) \Bf_i^\top
\end{align}
with material assignment $\mb\in \Uadmat$. Accordingly, the above approximation is computationally accessible, given that the entire \ac{tmm} approach relies on the individual T\nbhy matrices~$\Tf^i$, $i\in\{1,\ldots,\Npar\}$, i.e., these matrices are typically either known or can be easily computed, as in the case of radially symmetric particles.

\begin{remark}
The block-diagonal choice approximates the inverse matrix used to construct the surrogate and  does not remove particle interactions from the \ac{tmm} state or adjoint solves. Its accuracy may depend on the strength of inter-particle coupling. We do not derive an error bound, but assess model quality for a densely packed configuration in \cref{section:model_validation}.
\end{remark}

Using the block-diagonal approximation of the original \ac{tmm} system matrix allows for the derivation of a first-order separable approximation according to \cref{def:first_order_separable_model} of the objective function defined in \cref{eq:general_obj_func_TMM}, which is computationally accessible.

\begin{corollary}\label{cor:first_order_separable_approx_TMM}
    Let $\Lcalb\in\CC^{\bar{M}\times L}$ and $\tilde{\mb}\in \Uadmat$ be given, and let $\Acalb(\tilde{\mb})\in\CC^{\bar{M}}$ and $\Qcalb(\tilde{\mb})\in\CC^{\bar{M}\times L}$ be the unique solutions of the \ac{tmm} state and the adjoint equation given by \cref{eq:TMM_system_SGP,eq:general_adjoint_equation_TMM}, respectively, with material assignment~$\tilde{\mb}$. Furthermore, let $\Tcalb^\mathrm{mat}(\tilde{\mb})^{-1}\in\CC^{\bar{M}\times\bar{M}}$ be the block-diagonal approximation of the inverse \ac{tmm} system matrix defined in \cref{def:diagonal_approximation_TMM}. Then, according to \cref{def:first_order_separable_model}, the function $\Sdiag[J,\tilde{\mb}](\cdot)\colon \Uadmat \to \RR$ defined by
    \begin{align}\label{eq:diag_approx_sep_model_TMM}
    \begin{split}
        &\Sdiag[J,\tilde{\mb}](\mb) = J(\tilde{\mb})(1-\Npar) + \sum_{i=1}^\Npar g\Bigl(\Lcalb^\top \Acalb(\tilde{\mb}) + \boldsymbol{\zeta}^i\bigl[\tilde{\mb},\Tcalb^\mathrm{mat}(\tilde{\mb})^{-1}\bigr](m_i)\Bigr)
    \end{split}
    \end{align}
    is a first-order separable approximation of the objective function~$J$ from \cref{eq:general_obj_func_TMM} in the point of approximation~$\tilde{\mb}$. Moreover, the corresponding building block is available in the closed form
    \begin{align}\label{eq:X_matrix_diag_closed}
        \Xf^i\bigl[\tilde{\mb},\Tcalb^\mathrm{mat}(\tilde{\mb})^{-1}\bigr](m_i) = \Tf^i(\tilde m_i)^{-1} \bigl(\Tf^i(m_i)-\Tf^i(\tilde m_i)\bigr) \Tf^i(\tilde m_i)^{-1},
    \end{align}
    which requires no inverse of the \ac{tmm} system matrix.
\end{corollary}

\begin{proof}
    In \cref{thm:general_first_order_separable_approx_TMM}, we showed that we can achieve a first-order separable model by arbitrarily choosing a matrix~$\Kcalb$ as an approximation for the inverse \ac{tmm} system matrix $\Tcalb(\tilde{\mb})^{-1}$. Therefore, we just have to show that inserting the block-diagonal approximation~$\Tcalb^\mathrm{mat}(\tilde{\mb})^{-1}$ for the matrix~$\Kcalb$ in \cref{eq:general_approx_sep_model_TMM} yields the desired model~$\Sdiag[J,\tilde{\mb}](\mb)$. With the definition of $\Tcalb^\mathrm{mat}(\tilde{\mb})$, we get that
    \begin{align}
        \bigl(\one_{M^i} - \Delta\Tf^i(\tilde m_i,m_i) \Bf_i^\top \Tcalb^\mathrm{mat}(\tilde{\mb})^{-1} \Bf_i\bigr)^{-1} &= \biggl(\one_{M^i} - \Delta\Tf^i(\tilde m_i,m_i) \Bf_i^\top \sum_{j=1}^\Npar \Bf_j \Tf^j(\tilde m_j) \Bf_j^\top \Bf_i\biggr)^{-1}\\
        &= \bigl(\one_{M^i} - \Delta\Tf^i(\tilde m_i,m_i) \Tf^i(\tilde m_i)\bigr)^{-1}\\
        &= \Tf^i(\tilde m_i)^{-1}\Tf^i(m_i),
    \end{align}
    where $\Delta\Tf^i$ is defined by \cref{eq:delta_T}. Therefore, by \cref{eq:X_matrix},
    \begin{align}
        \Xf^i\bigl[\tilde{\mb},\Tcalb^\mathrm{mat}(\tilde{\mb})^{-1}\bigr](m_i) &= \bigl(\one_{M^i} -\Delta\Tf^i(\tilde m_i,m_i) \Bf_i^\top \Tcalb^\mathrm{mat}(\tilde{\mb})^{-1} \Bf_i\bigr)^{-1} \Delta\Tf^i(\tilde m_i,m_i)\\
        &= \Tf^i(\tilde m_i)^{-1} \bigl(\Tf^i(m_i) - \Tf^i(\tilde m_i)\bigr) \Tf^i(\tilde m_i)^{-1},
    \end{align}
    which is the closed form stated in \cref{eq:X_matrix_diag_closed} and finishes the proof.
\end{proof}

\begin{remark}[On the Invertibility Assumption]\label{rem:invertibility_assumption}
    The assumption of \cref{thm:exact_sep_model_TMM} that the matrix $\Xf^i[\tilde{\mb},\Tcalb(\tilde{\mb})^{-1}](m_i)$ be well defined is not an additional restriction. Indeed, the matrix determinant lemma applied to the rank\nbhy $M^i$ update computed in the proof of \cref{lem:single_element_update} gives
    \begin{align}\label{eq:determinant_lemma}
        \det\bigl(\Tcalb(\tilde{\mb}) - \Bf_i \Delta\Tf^i(\tilde m_i,m_i) \Bf_i^\top\bigr) = \det\bigl(\Tcalb(\tilde{\mb})\bigr)\cdot\det\bigl(\one_{M^i}-\Delta\Tf^i(\tilde m_i,m_i) \Bf_i^\top \Tcalb(\tilde{\mb})^{-1} \Bf_i\bigr),
    \end{align}
    so that, for an invertible matrix~$\Tcalb(\tilde{\mb})$, the matrix $\one_{M^i}-\Delta\Tf^i(\tilde m_i,m_i) \Bf_i^\top \Tcalb(\tilde{\mb})^{-1} \Bf_i$ inverted in \cref{eq:X_matrix} is invertible if and only if $\Tcalb\bigl(\tilde{\mb}+(m_i - \tilde m_i)\eb^i\bigr)$ is invertible, i.e., if and only if the \ac{tmm} state equation is uniquely solvable for the modified material assignment. Since every material assignment considered during the optimization has to admit a unique state, this is satisfied throughout. For the block-diagonal choice of \cref{cor:first_order_separable_approx_TMM} the assumption is even automatic, as the corresponding inverse is computed there in closed form as $\Tf^i(\tilde m_i)^{-1}\Tf^i(m_i)$, which exists by \cref{rem:invertibility_Tmatrix}, so that $\Xf^i\bigl[\tilde{\mb},\Tcalb^\mathrm{mat}(\tilde{\mb})^{-1}\bigr](m_i)$ is always well defined. For a general matrix~$\Kcalb$ in \cref{thm:general_first_order_separable_approx_TMM}, however, invertibility has to be required explicitly, since otherwise $\Sapprox[J,\tilde{\mb}]$ is not defined for every admissible material assignment.
\end{remark}

\begin{remark}[Computational Aspects]\label{rem:cost_sep_models}
    The above exact separable model requires the blocks $\Bf_i^\top \Tcalb(\tilde{\mb})^{-1}\Bf_i$ for all $i\in\{1,\ldots,\Npar\}$, i.e., all diagonal blocks of the inverse system matrix. Obtaining them amounts to $\bar{M}$ linear solves with $\Tcalb(\tilde{\mb})$, which is the reason why the exact model becomes prohibitive for larger particulate systems. The block-diagonal model~$\Sdiag$ avoids this entirely. It still requires one state solve for $\Acalb(\tilde{\mb})$ and $L$ adjoint solves, one for each column of~$\Lcalb$, per approximation point, where $L$ depends on the optical property under consideration. Once state and adjoint states are available, evaluating the model for a candidate material~$m_i$ only involves operations on $M^i\times M^i$ matrices.

    In the \ac{sgp} subproblem, the continuous material curve is sampled on a finite set of candidates, and the model is evaluated for each particle and each candidate independently. Consequently, the individual T\nbhy matrices $\Tf^i(m_i)$ have to be available for all candidate materials -- and, in the spectral setting of \cref{section:Lab_colour}, for all considered wavelengths -- which is unproblematic for radially symmetric particles, where they are given by the Mie coefficients.
\end{remark}

We have now derived an exact separable model as well as a corresponding computationally accessible first-order separable approximation of the objective function from \cref{eq:general_obj_func_TMM}. However, depending on the number~$L$ of linear functionals, it can be prohibitively expensive to compute the separable models directly, since $L$ adjoint problems have to be solved in every iteration of the \ac{sgp} algorithm. In this case, the following remark applies.

\begin{remark}[Linearization of the Objective Function]\label{rem:general_functions}
    Let $\tilde{\mb}\in\Uadmat$ be given, and let $\tilde{\mathbf{w}} = \Lcalb^\top\Acalb(\tilde{\mb})\in\CC^L$. Then, the objective function~$J$ from \cref{eq:general_obj_func_TMM} can be approximated first by linearizing the function~$g$ about the point~$\tilde{\mathbf{w}}$. With the gradient~$\nabla g$ from \cref{eq:complex_gradient}, this leads to
    \begin{align}
        \tilde{J}[\tilde{\mb}](\mb) = g(\tilde{\mathbf{w}}) + \Re\Bigl(\nabla g(\tilde{\mathbf{w}})^\dagger\Lcalb^\top\bigl(\Acalb(\mb) - \Acalb(\tilde{\mb})\bigr)\Bigr) = g(\tilde{\mathbf{w}}) + \Re\Bigl(\tilde{\Lcalb}^\top\bigl(\Acalb(\mb) - \Acalb(\tilde{\mb})\bigr)\Bigr)
    \end{align}
    with $\tilde{\Lcalb} = \Lcalb\,\overline{\nabla g(\tilde{\mathbf{w}})}\in\CC^{\bar{M}}$. Note that the linearized function is of the form of \cref{eq:general_obj_func_TMM} with $L=1$, the single vector~$\tilde{\Lcalb}$ and an affine function~$g$, so that only one adjoint problem has to be solved. Moreover, due to \cref{eq:chain_rule_g}, the linearized function~$\tilde{J}[\tilde{\mb}]$ reproduces both the value and the gradient of~$J$ in~$\tilde{\mb}$. Therefore, every first-order separable approximation of $\tilde{J}[\tilde{\mb}]$ in~$\tilde{\mb}$ is also a first-order separable approximation of~$J$.
\end{remark}

So far, we have considered objective functions that depend on the solution of a single \ac{tmm} state equation. However, some optical properties require the solution of the scattering problem for several incident fields, e.g., for several wavelengths or polarizations. In this case, the aforementioned results can be extended as follows. Let $K\in\NN$, and let us consider the $K$ state equations
\begin{align}\label{eq:TMM_systems_multiple}
    \Tcalb_k(\mb)\Acalb_k(\mb) = \Pcalb_{\mathrm{inc},k}, \qquad k\in\{1,\ldots,K\},
\end{align}
where each system matrix $\Tcalb_k(\mb)\in\CC^{\bar{M}\times\bar{M}}$ is of the structure given by \cref{eq:TMM_split_matrix,eq:TMM_diagonal_part} with the individual T\nbhy matrices $\Tf^i_k(m_i)$, $i\in\{1,\ldots,\Npar\}$. Accordingly, $\Delta\Tf^i_k$, $\Xf^i_k$ and $\Tcalb^\mathrm{mat}_k$ denote the quantities from \cref{eq:delta_T,eq:X_matrix,eq:TMM_diagonal_part} for the $k$\nbhy th state equation. Furthermore, let $L_1,\ldots,L_K\in\NN$, $\Lcalb_k\in\CC^{\bar{M}\times L_k}$ for all $k\in\{1,\ldots,K\}$, and $g\in\mathcal{C}^1(\CC^{L_1}\times\cdots\times\CC^{L_K},\RR)$ be given. Then, let us consider the objective function $J\colon\mathcal D\to\RR$ defined by
\begin{align}\label{eq:general_obj_func_TMM_multiple}
    J(\mb) = g\bigl(\Lcalb_1^\top\Acalb_1(\mb),\ldots,\Lcalb_K^\top\Acalb_K(\mb)\bigr).
\end{align}
Note that the function~$J$ from \cref{eq:general_obj_func_TMM} is included as the special case $K=1$.

\begin{corollary}[Separable Models for Several State Equations]\label{cor:several_state_equations}
    Let $\tilde{\mb}\in\Uadmat$ be given, and for all $k\in\{1,\ldots,K\}$, let $\Acalb_k(\tilde{\mb})\in\CC^{\bar{M}}$ and $\Qcalb_k(\tilde{\mb})\in\CC^{\bar{M}\times L_k}$ be the unique solutions of the state equation from \cref{eq:TMM_systems_multiple} and of the adjoint equation
    \begin{align}\label{eq:adjoint_equation_multiple}
        \Tcalb_k(\tilde{\mb})^\top\Qcalb_k = \Lcalb_k,
    \end{align}
    respectively, with material assignment~$\tilde{\mb}$. Moreover, for $i\in\{1,\ldots,\Npar\}$, $k\in\{1,\ldots,K\}$ and $\Kcalb_k\in\CC^{\bar{M}\times\bar{M}}$ such that the matrix $\Xf^i_k[\tilde{\mb},\Kcalb_k](m_i)$ is well defined, let
    \begin{align}\label{eq:zeta_multiple}
        \boldsymbol{\zeta}^i_k[\tilde{\mb},\Kcalb_k](m_i) = \Qcalb_k(\tilde{\mb})^\top\Bf_i\Xf^i_k[\tilde{\mb},\Kcalb_k](m_i)\Bf_i^\top\Acalb_k(\tilde{\mb}) \in\CC^{L_k}.
    \end{align}
    Then, the following statements hold:
    \begin{enumerate}
        \item\label{cor:several_state_equations_exact} The function $\Sexact[J,\tilde{\mb}](\cdot)\colon\Uadmat\to\RR$ defined by
        \begin{align}\label{eq:exact_sep_model_multiple}
            & \Sexact[J,\tilde{\mb}](\mb) \\
            &= J(\tilde{\mb})(1-\Npar) + \sum_{i=1}^\Npar g\Bigl(\Lcalb_1^\top\Acalb_1(\tilde{\mb}) + \boldsymbol{\zeta}^i_1\bigl[\tilde{\mb},\Tcalb_1(\tilde{\mb})^{-1}\bigr](m_i),\ldots,\Lcalb_K^\top\Acalb_K(\tilde{\mb}) + \boldsymbol{\zeta}^i_K\bigl[\tilde{\mb},\Tcalb_K(\tilde{\mb})^{-1}\bigr](m_i)\Bigr)
        \end{align}
        is an exact separable model according to \cref{def:exact_separable_model} of the objective function~$J$ from \cref{eq:general_obj_func_TMM_multiple}.
        \item\label{cor:several_state_equations_diag} The function $\Sdiag[J,\tilde{\mb}](\cdot)\colon\Uadmat\to\RR$, which is obtained by replacing the inverse system matrix $\Tcalb_k(\tilde{\mb})^{-1}$ in \cref{eq:exact_sep_model_multiple} by its block-diagonal approximation $\Tcalb^\mathrm{mat}_k(\tilde{\mb})^{-1}$ from \cref{def:diagonal_approximation_TMM} for all $k\in\{1,\ldots,K\}$, is a first-order separable approximation according to \cref{def:first_order_separable_model} of the objective function~$J$ from \cref{eq:general_obj_func_TMM_multiple} in the point of approximation~$\tilde{\mb}$. Moreover, the corresponding building blocks $\Xf^i_k\bigl[\tilde{\mb},\Tcalb^\mathrm{mat}_k(\tilde{\mb})^{-1}\bigr](m_i)$ are given by \cref{eq:X_matrix_diag_closed} with $\Tf^i$ replaced by $\Tf^i_k$.
    \end{enumerate}
\end{corollary}

\begin{proof}
Apply \cref{lem:single_element_update} to each state $k$. For a change of one particle material vector, every component $\Lcalb_k^\top\Acalb_k$ is recovered exactly when $\Kcalb_k=\Tcalb_k(\tilde\mb)^{-1}$, so \cref{eq:exact_sep_model_multiple} is the exact separable model. With the block-diagonal choices, each component update has the correct value and directional derivative at $\tilde\mb$. The real chain rule for the function $g$ therefore gives first-order agreement of the full objective in every direction $h_i\in\MatSpace$. Summing these partial derivatives proves the assertion. The closed form follows from \cref{eq:X_matrix_diag_closed} separately for each $k$.
\end{proof}

In spectral computations the state index $k$ may run over sampled wavelengths and incident polarizations. Therefore, the same wavelength may occur more than once. Each state depends on the appropriate component $m_{i,q(k)}$ of the same particle vector $m_i$. Under \cref{eq:linear_parameterization}, a scalar variation $\delta\rho_i$ produces the coherent changes $h_{i,q}=(\nbar_2(\lambda_q)-\nbar_1(\lambda_q))\,\delta\rho_i$ at all sampled wavelengths. Thus the finite-sum objective and surrogate have identical derivatives with respect to every scalar $\rho_i$ by the ordinary chain rule.

Note that, in every iteration of the \ac{sgp} algorithm, one state equation and $L_k$ adjoint problems have to be solved for each $k\in\{1,\ldots,K\}$. If two state equations differ only in the incident field, e.g., for two polarizations with the same wavelength, the system matrices coincide, and only the right-hand sides differ. Furthermore, \cref{rem:general_functions} applies analogously, i.e., linearizing the function~$g$ with respect to all its arguments leads to $L_k=1$ for all $k\in\{1,\ldots,K\}$.

\begin{remark}[Continuous Spectral Setting]\label{rem:continuous_spectral_extension}
    The preceding analysis uses finitely many sampled wavelengths, as in the numerical experiments. The same low-rank construction may be formulated for material functions on a continuous wavelength interval, provided the single-particle and coupled \ac{tmm} systems are uniformly invertible on the relevant admissible materials, their material dependence has suitable uniform regularity, and differentiation can be interchanged with the spectral integral. These additional function-space assumptions are not needed for the fixed quadrature considered here.
\end{remark}

Let us now apply the aforementioned results to the optical properties introduced in \cref{section:TMM_optical_properties}. The extinction cross-section as well as the directional and the monitor scattering intensity are of the form of \cref{eq:general_obj_func_TMM}. In these cases, the function~$g$ is either an affine function of the real part with $L=1$ or a weighted sum of squared absolute values, i.e.,
\begin{align}
    g(w) &= c_0 + c_1\Re(w) &&\text{with}\quad \nabla g(w) = c_1,\\
    g(\mathbf{w}) &= \sum_{j=1}^L \gamma_j\vert \mathbf{w}_j\vert^2 &&\text{with}\quad \nabla g(\mathbf{w}) = \bigl(2\gamma_1\mathbf{w}_1,\ldots,2\gamma_L\mathbf{w}_L\bigr)^\top,
\end{align}
with constants $c_0,c_1,\gamma_1,\ldots,\gamma_L\in\RR$. In particular, the polynomial objective functions $g = p\circ\Re$ and $g = p\circ\Im$ with $L=1$ and a real polynomial~$p$, are included as special cases with $\nabla g(w) = p'\bigl(\Re(w)\bigr)$ and $\nabla g(w) = \ci\,p'\bigl(\Im(w)\bigr)$, respectively. Moreover, the colour functions from \cref{section:Lab_colour} depend on the states at all sampled wavelengths and are therefore treated based on \cref{cor:several_state_equations}.

\subsection{Extinction Cross-Section}
In the \ac{tmm} framework, according to \cref{def:extinction_TMM}, the extinction cross-section is defined by
\begin{align}
    \sigma_\mathrm{ext}(\mb) = -\frac{\pi}{k_\mathrm{b}^2}  \Re\bigl(\Pcalb_\mathrm{inc}^\dagger \Acalb(\mb)\bigr),
\end{align}
where $\mb\in \MatSpace^\Npar$ denotes the material assignment. This optical property is of the form of the generally considered objective function~$J$ from \cref{eq:general_obj_func_TMM} with $L=1$, $g(w) = -\frac{\pi}{k_\mathrm{b}^2}\Re(w)$ for all $w\in\CC$, and $\Lcalb = \overline{\Pcalb_\mathrm{inc}}$. Employing \cref{thm:exact_sep_model_TMM}, an exact separable model of the aforementioned extinction cross-section can be expressed as follows:
\begin{align}
    &\Sexact[\sigma_\mathrm{ext},\tilde{\mb}](\mb) \\
    & = \sigma_\mathrm{ext}(\tilde{\mb})(1-\Npar) + \sum_{i=1}^\Npar  - \frac{\pi}{k_\mathrm{b}^2} \Re\Bigl(\Pcalb_\mathrm{inc}^\dagger \Acalb(\tilde{\mb}) + \Qcalb(\tilde{\mb})^\top\Bf_i \Xf^i\bigl[\tilde{\mb},\Tcalb(\tilde{\mb})^{-1}\bigr](m_i)  \Bf_i^\top\Acalb(\tilde{\mb})\Bigr) \\
    & = \sigma_\mathrm{ext}(\tilde{\mb}) - \frac{\pi}{k_\mathrm{b}^2} \sum_{i=1}^\Npar \Re\Bigl(\Qcalb(\tilde{\mb})^\top \Bf_i \Xf^i\bigl[\tilde{\mb},\Tcalb(\tilde{\mb})^{-1}\bigr](m_i)  \Bf_i^\top\Acalb(\tilde{\mb})\Bigr),\label{eq:exact_sep_model_extinction_TMM}
\end{align}
where $\Qcalb(\tilde{\mb})\in \CC^{\bar M}$ is the unique solution of the adjoint equation
\begin{align}
    \Tcalb(\tilde{\mb})^\top \Qcalb = \overline{\Pcalb_\mathrm{inc}}.
\end{align}
Hence, in addition to the state equation, only one adjoint problem has to be solved in every iteration.

A corresponding first-order separable approximation of the extinction cross-section in the \ac{tmm} setting is given by applying \cref{cor:first_order_separable_approx_TMM} as 
\begin{align}\label{eq:sep_model_extinction_TMM}
   &\Sdiag[\sigma_\mathrm{ext},\tilde{\mb}](\mb) = \sigma_\mathrm{ext}(\tilde{\mb}) - \frac{\pi}{k_\mathrm{b}^2} \sum_{i=1}^\Npar \Re\Bigl(\Qcalb(\tilde{\mb})^\top \Bf_i \Xf^i\bigl[\tilde{\mb},\Tcalb^\mathrm{mat}(\tilde{\mb})^{-1}\bigr](m_i)  \Bf_i^\top \Acalb(\tilde{\mb})\Bigr).
\end{align}

\subsection{Directional Scattering}
For a fixed scattering direction $\ab$, the intensity from \cref{eq:directional_scattering} is $s_{\ab}(\mb)=k_\mathrm{b}^{-2}\bigl\|\Ecalb(\ab)^\top\Acalb(\mb)\bigr\|_2^2$. Thus, in \cref{eq:general_obj_func_TMM} we take $L=2$, $\Lcalb=\Ecalb(\ab)$ and $g(\mathbf w)=k_\mathrm{b}^{-2}\|\mathbf w\|_2^2$. Let $\Qcalb(\tilde\mb)\in\CC^{\bar M\times 2}$ solve the adjoint equation $\Tcalb(\tilde\mb)^\top\Qcalb=\Ecalb(\ab)$, so that the particle-level observable update from \cref{eq:zeta_scalar} reads
\begin{align}\label{eq:z_vector}
 \boldsymbol{\zeta}^i[\tilde\mb,\Kcalb](m_i)=\Qcalb(\tilde\mb)^\top\Bf_i\Xf^i[\tilde\mb,\Kcalb](m_i)\Bf_i^\top\Acalb(\tilde\mb)\in\CC^2.
\end{align}
Writing $\tilde{\mathbf w}=\Ecalb(\ab)^\top\Acalb(\tilde\mb)$, \cref{thm:exact_sep_model_TMM} gives the exact separable model
\begin{align}\label{eq:exact_sep_model_directional_scattering_TMM}
 \Sexact[s_{\ab},\tilde\mb](\mb)=s_{\ab}(\tilde\mb)+\frac{1}{k_\mathrm{b}^2}\sum_{i=1}^\Npar\Bigl(\bigl\|\tilde{\mathbf w}+\boldsymbol{\zeta}^i\bigl[\tilde\mb,\Tcalb(\tilde\mb)^{-1}\bigr](m_i)\bigr\|_2^2-\|\tilde{\mathbf w}\|_2^2\Bigr).
\end{align}
The computationally accessible first-order model has the same expression, but uses the block-diagonal choice in each update:
\begin{align}\label{eq:sep_model_directional_scattering_TMM}
 \Sdiag[s_{\ab},\tilde\mb](\mb)=s_{\ab}(\tilde\mb)+\frac{1}{k_\mathrm{b}^2}\sum_{i=1}^\Npar\Bigl(\bigl\|\tilde{\mathbf w}+\boldsymbol{\zeta}^i\bigl[\tilde\mb,\Tcalb^\mathrm{mat}(\tilde\mb)^{-1}\bigr](m_i)\bigr\|_2^2-\|\tilde{\mathbf w}\|_2^2\Bigr).
\end{align}
Both models retain the quadratic structure of~$s_\ab$.

\subsection{Monitor Scattering}\label{section:TMM_SGP_monitor_scattering}
The monitor intensity, integrated over $\tilde n$ scattering directions with weights $\omega_t$, is a quadratic function of the state:
\begin{align}\label{eq:g_monitor_scat}
 s_\Mcal(\mb)=\frac{1}{k_\mathrm{b}^2}\sum_{t=1}^{\tilde n}\omega_t\bigl\|\Ecalb(\ab^t)^\top\Acalb(\mb)\bigr\|_2^2.
\end{align}
Equivalently, with $\Fcalb=\sum_t\omega_t\Ecalb(\ab^t)\Ecalb(\ab^t)^\dagger$ (see \cref{eq:matrix_F}),
\begin{align}\label{eq:monitor_scat_function_TMM}
 s_\Mcal(\mb)=k_\mathrm{b}^{-2}\Acalb(\mb)^\top\Fcalb\overline{\Acalb(\mb)}.
\end{align}
Direct application of the general model to all $2\tilde n$ linear observables would need $2\tilde n$ adjoint solves. Instead, as in \cref{rem:general_functions}, we first linearize this quadratic form with respect to the state at $\Acalb(\tilde\mb)$:
\begin{align}\label{eq:linearized_monitor_scattering_function}
 \tilde s_\Mcal[\tilde\mb](\mb)=-s_\Mcal(\tilde\mb)+\frac{2}{k_\mathrm{b}^2}\Re\bigl(\Acalb(\tilde\mb)^\dagger\overline{\Fcalb}\Acalb(\mb)\bigr).
\end{align}
The linearized function agrees with $s_\Mcal$ in value and first derivative at $\tilde\mb$. It has one linear state observable, so a single adjoint solve is sufficient: $\Qcalb(\tilde\mb)$ solves $\Tcalb(\tilde\mb)^\top\Qcalb=\Fcalb\overline{\Acalb(\tilde\mb)}$. The exact separable model of this linearized function is
\begin{align}\label{eq:ex_sep_model_monitor_scat}
 \Sexact\bigl[\tilde s_\Mcal[\tilde\mb],\tilde\mb\bigr](\mb)=s_\Mcal(\tilde\mb)+\frac{2}{k_\mathrm{b}^2}\sum_{i=1}^\Npar\Re\Bigl(\Qcalb(\tilde\mb)^\top\Bf_i\Xf^i\bigl[\tilde\mb,\Tcalb(\tilde\mb)^{-1}\bigr](m_i)\Bf_i^\top\Acalb(\tilde\mb)\Bigr).
\end{align}
In particular, the single-particle value needed for colour optimization is
\begin{align}\label{eq:single_element_monitor_scat}
 \tilde s_\Mcal[\tilde\mb]\bigl(\tilde\mb+(m_i-\tilde m_i)\eb^i\bigr)=s_\Mcal(\tilde\mb)+\frac{2}{k_\mathrm{b}^2}\Re\Bigl(\Qcalb(\tilde\mb)^\top\Bf_i\Xf^i\bigl[\tilde\mb,\Tcalb(\tilde\mb)^{-1}\bigr](m_i)\Bf_i^\top\Acalb(\tilde\mb)\Bigr).
\end{align}
Replacing the auxiliary inverse $\Tcalb(\tilde\mb)^{-1}$ by $\Tcalb^\mathrm{mat}(\tilde\mb)^{-1}$ gives $\Sdiag\bigl[\tilde s_\Mcal[\tilde\mb],\tilde\mb\bigr]$, a first-order model of both the linearized and the original monitor intensity. The linearized intensity need not be nonnegative away from the reference point, see \cref{rem:negative_reflectance}.

\subsection{Colour Functions}
Let $\Phi_1\colon\mathcal D\to\RR^3$ map the material assignment to its CIE\nbhy $\Lcolour\acolour\bcolour$ values, computed from the reflectance spectrum as in \cref{section:Lab_colour}. Furthermore, let $\Phi_2\colon\RR^3\to\RR$ denote a scalar colour objective differentiable near the reference colour, such as chroma away from the achromatic axis. Then the objective is the composition
\begin{align}\label{eq:general_col_function}
 \Phi(\mb)=\Phi_2\bigl(\Phi_1(\mb)\bigr)
\end{align}
with $\Phi_1(\mb)=\bigl(\Lcolour(\mb),\acolour(\mb),\bcolour(\mb)\bigr)$. The CIE\nbhy $\mathit{XYZ}$ values depend linearly on the quadrature-sampled reflectance. The subsequent CIE\nbhy $\Lcolour\acolour\bcolour$ conversion and $\Phi_2$ are generally nonlinear. Our construction preserves the latter nonlinearities in the surrogate.

\begin{remark}[Differentiability of the Colour Property Functions]\label{rem:colour_differentiability}
The CIE\nbhy $\Lcolour\acolour\bcolour$ conversion $f$ in \cref{eq:Lab_colour_f} is continuously differentiable: its two branches agree in value and derivative at the transition $\varepsilon=(6/29)^3$. The chroma $\Ccolour=\sqrt{(\acolour)^2+(\bcolour)^2}$ is not differentiable at the achromatic point $\acolour=\bcolour=0$, but is differentiable whenever $\Ccolour>0$. For optimization towards zero chroma, the smooth objective $(\Ccolour)^2$ can be used instead. We assume $\Phi_2$ is continuously differentiable in a neighbourhood of the colour being approximated.
\end{remark}

We now display the suppressed wavelength dependence explicitly: $\Tf^i(z,\lambda)$ and its inverse depend on the refractive index at $\lambda$, and $\Tcalb(\mb,\lambda)$ and $\Acalb(\mb,\lambda)$ depend on the entire particle assembly at that wavelength. The coupling, incident-field and far-field operators likewise depend on $\lambda$, and we denote the corresponding monitor matrix by $\Fcalb_\lambda$. For unpolarized illumination, the reflectance is the mean of two orthogonal incident-polarization results, evaluated using the same system matrix but different right-hand sides. Below we display one polarization. In the unpolarized case, the two terms are averaged.

As discussed in \cref{section:TMM_SGP_monitor_scattering} for the monitor scattering intensity, the separable models of the original reflectance would require $2\tilde{n}$ adjoint problems for each wavelength. Therefore, we replace the monitor scattering intensity $s_\Mcal(\cdot,\lambda)$ by its linearization $\tilde{s}_\Mcal[\tilde{\mb}](\cdot,\lambda)$ from \cref{eq:linearized_monitor_scattering_function} in the approximation point~$\tilde{\mb}$ at each sampled wavelength $\lambda_q$, $q=1,\ldots,Q$. This leads to the linearized reflectance
\begin{align}\label{eq:colour_reflectance}
    \tilde{r}_\Mcal[\tilde{\mb}](\mb,\lambda) = \frac{1}{A_\mathrm{p}(\Omega,\db)}\tilde{s}_\Mcal[\tilde{\mb}](\mb,\lambda), \qquad \lambda\in\{\lambda_1,\ldots,\lambda_Q\}.
\end{align}
Since the material optimization does not change the geometry of the particulate system, $A_\mathrm{p}(\Omega,\db)$ is a constant of the optimization problem. Evaluating the spectral integrals in \cref{eq:CIE_XYZ_1,eq:CIE_XYZ_2,eq:CIE_XYZ_3} with the fixed quadrature weights $w_q$ (and corresponding quadrature normalization $\bar y$) yields the CIE\nbhy $\mathit{XYZ}$ values
\begin{align}
    \tilde{X}[\tilde\mb](\mb) &=\frac{1}{\bar{y}} \sum_{q=1}^Q w_q\,x(\lambda_q)\,\tilde{r}_\Mcal[\tilde{\mb}](\mb,\lambda_q)\, I(\lambda_q),\label{eq:sep_model_col_1}\\
    \tilde{Y}[\tilde\mb](\mb) &=\frac{1}{\bar{y}} \sum_{q=1}^Q w_q\,y(\lambda_q)\,\tilde{r}_\Mcal[\tilde{\mb}](\mb,\lambda_q)\, I(\lambda_q),\label{eq:sep_model_col_2}\\
    \tilde{Z}[\tilde\mb](\mb) &=\frac{1}{\bar{y}} \sum_{q=1}^Q w_q\,z(\lambda_q)\,\tilde{r}_\Mcal[\tilde{\mb}](\mb,\lambda_q)\, I(\lambda_q).\label{eq:sep_model_col_3}
\end{align}
Moreover, let $\tilde{L}^\ast[\tilde\mb](\mb)$, $\tilde{a}^\ast[\tilde\mb](\mb)$ and $\tilde{b}^\ast[\tilde\mb](\mb)$ denote the corresponding CIE\nbhy $\Lcolour\acolour\bcolour$ values according to \cref{eq:CIE_Lab_1,eq:CIE_Lab_2,eq:CIE_Lab_3}. Then, the surrogate $\tilde{\Phi}[\tilde{\mb}]\colon\Uadmat\to\RR$ of the colour property function~$\Phi$ is given by
\begin{align}
    \tilde{\Phi}[\tilde{\mb}](\mb) = \Phi_2\bigl(\tilde{L}^\ast[\tilde\mb](\mb),\tilde{a}^\ast[\tilde\mb](\mb),\tilde{b}^\ast[\tilde\mb](\mb)\bigr).
\end{align}

Note that it would also be possible to linearize the entire colour property function~$\Phi$ as described in \cref{rem:general_functions}, which also leads to one adjoint problem per wavelength. However, the nonlinear colour transformation, i.e., the functions~$f$ and $\Phi_2$, is then lost in the model. In contrast, linearizing only the monitor scattering intensity leads to the same computational effort, but retains the nonlinear structure of the colour transformation.

According to \cref{rem:general_functions}, the linearized function $\tilde{s}_\Mcal[\tilde{\mb}](\cdot,\lambda)$ reproduces both the value and the first derivative of $s_\Mcal(\cdot,\lambda)$ in the point of approximation~$\tilde{\mb}$ at every sampled wavelength. Since the CIE\nbhy $\mathit{XYZ}$ values depend linearly on the reflectance spectrum, and since the function~$f$ is continuously differentiable (see \cref{rem:colour_differentiability}), the chain rule yields
\begin{align}
    \tilde{\Phi}[\tilde{\mb}](\tilde{\mb}) = \Phi(\tilde{\mb}) \qquad\text{and}\qquad D\tilde{\Phi}[\tilde{\mb}](\tilde{\mb}) = D\Phi(\tilde{\mb}),
\end{align}
provided that $\Phi_2$ is differentiable in a neighbourhood of the reference colour. The calculation uses a finite sum, so no interchange of differentiation with a spectral integral is required. Consequently, every first-order separable approximation of $\tilde{\Phi}[\tilde{\mb}]$ in~$\tilde{\mb}$ is also a first-order separable approximation of~$\Phi$ in~$\tilde{\mb}$.

The surrogate $\tilde{\Phi}[\tilde{\mb}]$ depends on the states $\Acalb(\mb,\lambda_q)$ at all $Q$ sampled wavelengths. It is of the form of \cref{eq:general_obj_func_TMM_multiple} with $Q$ state equations, $L_k=1$, and $\Lcalb_k=\Fcalb_{\lambda_k}\overline{\Acalb(\tilde{\mb},\lambda_k)}$ for $k=1,\ldots,Q$. For unpolarized illumination, the number of states is $2Q$ (one per wavelength and polarization). Here, the function~$g$ comprises the affine functions of the linearized monitor scattering intensities (cf. \cref{section:TMM_SGP_monitor_scattering}), the evaluation of \cref{eq:sep_model_col_1,eq:sep_model_col_2,eq:sep_model_col_3}, the colour transformation from \cref{eq:CIE_Lab_1,eq:CIE_Lab_2,eq:CIE_Lab_3} and the function~$\Phi_2$. Hence, \cref{cor:several_state_equations} can be applied. For readability, the following expressions use the generic symbol $\lambda$, always evaluated at a sampled node $\lambda_q$. The adjoint problems from \cref{eq:adjoint_equation_multiple} then read
\begin{align}\label{eq:colour_adjoint}
    \Tcalb(\tilde{\mb},\lambda)^\top\Qcalb = \Fcalb_\lambda\overline{\Acalb(\tilde{\mb},\lambda)},
\end{align}
and their unique solutions are denoted by $\Qcalb(\tilde{\mb},\lambda)$. Consequently, one state and one adjoint problem have to be solved per sampled wavelength and polarization when constructing the surrogate (while any shared system matrices can be factorized once).

According to \cref{eq:single_element_monitor_scat}, applied for each wavelength, the linearized reflectance for a change of the material in the $i$\nbhy th design element only is given by
\begin{align}\label{eq:sep_model_col_exact}
    & \tilde{r}_\Mcal[\tilde{\mb}]\bigl(\tilde{\mb} + (m_i - \tilde m_i)\eb^i,\lambda\bigr)= \frac{1}{A_\mathrm{p}(\Omega,\db)}\biggl(s_\Mcal(\tilde{\mb}, \lambda) + \frac{2}{k_\mathrm{b}(\lambda)^2}\Re\Bigl(\Qcalb(\tilde{\mb}, \lambda)^\top\Bf_i\Xf^i_\lambda\bigl[\tilde{\mb},\Tcalb(\tilde{\mb}, \lambda)^{-1}\bigr](m_i) \Bf_i^\top\Acalb(\tilde{\mb}, \lambda)\Bigr)\biggr),
\end{align}
where $s_{\Mcal}(\cdot,\lambda)$ is the original, non-linearized monitor scattering intensity from \cref{eq:monitor_scat_function_TMM} for the wavelength~$\lambda$. Inserting this expression into \cref{eq:sep_model_col_1,eq:sep_model_col_2,eq:sep_model_col_3} yields the values $\tilde{X}[\tilde\mb]\bigl(\tilde{\mb} + (m_i - \tilde m_i)\eb^i\bigr)$, $\tilde{Y}[\tilde\mb]\bigl(\tilde{\mb} + (m_i - \tilde m_i)\eb^i\bigr)$ and $\tilde{Z}[\tilde\mb]\bigl(\tilde{\mb} + (m_i - \tilde m_i)\eb^i\bigr)$, and therefore the corresponding CIE\nbhy $\Lcolour\acolour\bcolour$ values. Then, Statement~\ref{cor:several_state_equations_exact} of \cref{cor:several_state_equations} yields the exact separable model
\begin{align}\label{eq:ex_sep_model_q}
    &\Sexact\bigl[\tilde{\Phi}[\tilde{\mb}],\tilde\mb\bigr](\mb) \\
    &= \Phi(\tilde \mb)(1-\Npar) + \sum_{i=1}^\Npar \Phi_2\Bigl(\tilde{L}^\ast[\tilde\mb]\bigl(\tilde{\mb} + (m_i - \tilde m_i)\eb^i\bigr),\tilde{a}^\ast[\tilde\mb]\bigl(\tilde{\mb} + (m_i - \tilde m_i)\eb^i\bigr),\tilde{b}^\ast[\tilde\mb]\bigl(\tilde{\mb} + (m_i - \tilde m_i)\eb^i\bigr)\Bigr)
\end{align}
of the surrogate~$\tilde{\Phi}[\tilde{\mb}]$, where we used that $\tilde{\Phi}[\tilde{\mb}](\tilde{\mb}) = \Phi(\tilde{\mb})$. Note that this model is exact with respect to the surrogate $\tilde{\Phi}[\tilde{\mb}]$. With respect to the original colour property function~$\Phi$, it remains an approximation, which is first-order correct in the point of approximation~$\tilde\mb$.

Similarly, a corresponding first-order separable approximation of the colour property function~$\Phi$ is obtained by Statement~\ref{cor:several_state_equations_diag} of \cref{cor:several_state_equations}. Therefore, the inverse system matrix $\Tcalb(\tilde{\mb},\lambda)^{-1}$ in \cref{eq:sep_model_col_exact} is replaced by the block-diagonal approximation $\Tcalb^\mathrm{mat}(\tilde{\mb},\lambda)^{-1}$ from \cref{def:diagonal_approximation_TMM}, which yields
\begin{align}
    & \tilde{r}_\Mcal[\tilde{\mb}]\bigl(\tilde{\mb} + (m_i - \tilde m_i)\eb^i,\lambda\bigr) \approx \frac{1}{A_\mathrm{p}(\Omega,\db)}\biggl(s_{\Mcal}(\tilde\mb,\lambda) +\frac{2}{k_\mathrm{b}(\lambda)^2}\Re\Bigl(\Qcalb(\tilde{\mb}, \lambda)^\top \Bf_i \Xf^i_\lambda\bigl[\tilde{\mb},\Tcalb^\mathrm{mat}(\tilde{\mb},\lambda)^{-1}\bigr](m_i) \Bf_i^\top \Acalb(\tilde{\mb}, \lambda)\Bigr)\biggr)\label{eq:sep_model_col_diag}
\end{align}
with (cf. \cref{eq:X_matrix_diag_closed})
\begin{align}
    \Xf^i_\lambda\bigl[\tilde{\mb},\Tcalb^\mathrm{mat}(\tilde{\mb},\lambda)^{-1}\bigr](m_i) = \Tf^i(\tilde m_i,\lambda)^{-1} \bigl(\Tf^i(m_i,\lambda)-\Tf^i(\tilde m_i,\lambda)\bigr) \Tf^i(\tilde m_i,\lambda)^{-1}.
\end{align}
Inserting the right-hand side of \cref{eq:sep_model_col_diag} into \cref{eq:sep_model_col_1,eq:sep_model_col_2,eq:sep_model_col_3}, computing the corresponding CIE\nbhy $\Lcolour\acolour\bcolour$ values and inserting them into \cref{eq:ex_sep_model_q} provides the desired first-order separable approximation, which we denote by $\Sdiag\bigl[\tilde{\Phi}[\tilde{\mb}],\tilde{\mb}\bigr]$. As before, $\Sdiag\bigl[\tilde{\Phi}[\tilde{\mb}],\tilde{\mb}\bigr]$ is also a first-order separable approximation of the original colour property function~$\Phi$ in~$\tilde{\mb}$.
Note that this model requires no inverse of the system matrices~$\Tcalb(\tilde{\mb},\lambda)$.

\begin{remark}[Colour Functions Based on Other Optical Properties]
    In the case of colour functions based on optical properties other than the (linearized) monitor scattering intensity, the corresponding separable models can be obtained by employing the models derived for that optical property for the reflectance term in \cref{eq:sep_model_col_1,eq:sep_model_col_2,eq:sep_model_col_3}, and therefore for $\tilde{L}^\ast[\tilde\mb]\bigl(\tilde{\mb} + (m_i - \tilde m_i)\eb^i\bigr)$, $\tilde{a}^\ast[\tilde\mb]\bigl(\tilde{\mb} + (m_i - \tilde m_i)\eb^i\bigr)$ and $\tilde{b}^\ast[\tilde\mb]\bigl(\tilde{\mb} + (m_i - \tilde m_i)\eb^i\bigr)$, respectively.
\end{remark}

The above derivations can then be employed to express an exact separable model as well as a first-order separable approximation of different colour properties. As an example, the exact separable model of the surrogate $\tilde{C}^\ast[\tilde{\mb}]$ of the chroma function defined in \cref{eq:chroma_function}, based on the linearized monitor scattering intensity in the \ac{tmm} setting, reads
\begin{align}\label{eq:chroma_sep_model}
    &\Sexact\bigl[\tilde{C}^\ast[\tilde{\mb}],\tilde \mb\bigr](\mb) = \Ccolour(\tilde\mb)(1-\Npar) + \sum_{i=1}^{\Npar} \sqrt{\tilde{a}^\ast[\tilde\mb]\bigl(\tilde{\mb} + (m_i - \tilde m_i)\eb^i\bigr)^2 + \tilde{b}^\ast[\tilde\mb]\bigl(\tilde{\mb} + (m_i - \tilde m_i)\eb^i\bigr)^2}
\end{align}
where the tilde colour coordinates are obtained from \cref{eq:sep_model_col_1,eq:sep_model_col_2,eq:sep_model_col_3} and the standard CIE\nbhy $\Lcolour\acolour\bcolour$ transformation \cref{eq:CIE_Lab_1,eq:CIE_Lab_2,eq:CIE_Lab_3}. Inserting the block-diagonal reflectance update from \cref{eq:sep_model_col_diag}, defines the computationally accessible $\Sdiag\bigl[\tilde C^\ast[\tilde\mb],\tilde\mb\bigr]$ used in \cref{section:opt_photonic_crystal}.

\begin{remark}[Admissibility of the Linearized Reflectance]\label{rem:negative_reflectance}
    As noted in \cref{section:TMM_SGP_monitor_scattering}, the linearized monitor scattering intensity~$\tilde{s}_\Mcal[\tilde\mb]$, and hence the linearized reflectance~$\tilde{r}_\Mcal[\tilde{\mb}]$ from \cref{eq:colour_reflectance}, is not sign constrained. Away from the point of approximation, the model may therefore produce negative values of $\tilde{r}_\Mcal[\tilde{\mb}](\cdot,\lambda)$ and, consequently, negative tristimulus values. Since the function~$f$ from \cref{eq:Lab_colour_f} is defined on~$\RR$, the separable models remain well defined in this case. However, the resulting triple can no longer be interpreted as a colour. Note that this only affects the separable model and not the original colour property function~$\Phi$, which is evaluated in the inner loop of \cref{alg:SGP_alg} in order to accept or reject a design update. Moreover, the globalization strategy of \cref{alg:SGP_alg} increases the parameter~$\tau$ of the proximal point term $\tau\Vert\mb^n-\mb\Vert_p^p$ until the objective function decreases sufficiently. Hence, the design update is restricted to a neighbourhood of the point of approximation, in which the linearization is accurate.
\end{remark}
\section{Numerical Experiments}\label{section:numerical_experiments}
In the following examples, we consider optimization problems of the form
\begin{align}\label{eq:general_obj_func_examples}
    \min_{\mb\in\Uadmat}\ & J(\mb),
\end{align}
where $J\colon\MatSpace^\Npar\to\RR$, $\mb\mapsto J(\mb)$ is a real-valued objective function which requires the solution of the electromagnetic scattering problem described in detail in \cref{section:general_electromagnetic_scattering_problem}. We solve this problem by the \ac{tmm} introduced in \cref{section:TMM}. This means that the governing state equation of {\crefname{equation}{Problem}{Problems}\cref{eq:general_obj_func_examples}} with material assignment $\mb\in\MatSpace^\Npar$ is given by
\begin{align}
    \Tcalb(\mb)\Acalb(\mb) = \Pcalb_\mathrm{inc},
\end{align}
see \cref{eq:TMM_system,eq:general_relaxed_mat_problem_TMM}. Furthermore, the admissible set of materials~$\Uadmat$ is given by \cref{eq:material_admissible_set}, which yields $\normalfont \Uadmat = (\mathcal{V}_\mathrm{ad})^\Npar$ with the admissible set of a single design element~$\normalfont \mathcal{V}_\mathrm{ad}$ from \cref{eq:admissible_set_single_element}. We use the \ac{tmm}\nbhy \ac{sgp} framework discussed in \cref{section:material_optimization} to solve the above optimization problem.

The numerical experiments considered in this work are structured as follows: In \cref{section:model_validation}, we minimize the extinction cross-section of a pyramid of \num{35}~particles and validate the accuracy of the employed separable model. In \cref{section:opt_photonic_crystal}, we then optimize the colour of a photonic crystal consisting of \num{627}~particles, which involves a more complex objective function and a larger design domain.

In the following, the material property of each design element, i.e., of each particle of the object, is characterized by the complex refractive index $\nbar\in\CC$. Furthermore, let $\mathbf{1}_\Npar \in \RR^\Npar$ be a vector of size~$\Npar$ with all entries equal to \num{1}, i.e., $\mathbf{1}_\Npar \coloneqq \sum_{i=1}^\Npar \eb^i$, where $\eb^i$ is the $i$\nbhy th unit vector in $\RR^\Npar$. In our experiments, we always consider a plane electromagnetic wave as an incident field.

Constructing a surrogate requires one coupled state solve and one adjoint solve per wavelength and polarization for the linearized observables used below. Further evaluations of the fully coupled objective may be needed in the globalization loop of \cref{alg:SGP_alg}. Each scalar subproblem is minimized over a fixed grid of $\rho\in[0,1]$ with step size~\num{e-3}, i.e., the separable model is evaluated \num{1001}~times per design element. Furthermore, we choose the following algorithmic parameters. The initial proximal point parameter is set to $\tau = \num{e-8}$, which is relatively small. Consequently, the algorithm is allowed to perform large steps with respect to the design updates in the beginning. The parameter~$p$ for the proximal point term used for the solution of the \ac{sgp} subproblem (see step 4 of \cref{alg:SGP_alg}) is chosen to be $p=2$. Additionally, the update parameters for the globalization strategy are set to $\overline{\vartheta} = \num{10}$ and $\underline{\vartheta} = \frac{1}{1.2}\overline{\vartheta}^{-1}$.

The \ac{tmm}\nbhy \ac{sgp} approach was implemented in MATLAB ver. R2025a, where we applied the open-source project CELES by \citet[]{egel2017} to solve the \ac{tmm} system.

\subsection{Validation of Models}\label{section:model_validation}
\begin{figure}[!htb]
    \centering
    \begin{tikzpicture}[baseline=-2pt]
        \draw[line width=1pt,->] (0,0) -- (0.45,0.19) node[anchor=west] {$x_1$};
        \draw[line width=1pt,->] (0,0) -- (-0.36,0.28) node[anchor=east] {$x_2$};
        \draw[line width=1pt,->] (0,0) -- (0,0.55) node[anchor=south] {$x_3$};
    \end{tikzpicture}\hspace{0.8em}%
    \includegraphics[trim={95pt 40pt 32pt 81pt}, clip, width=0.25\textwidth]{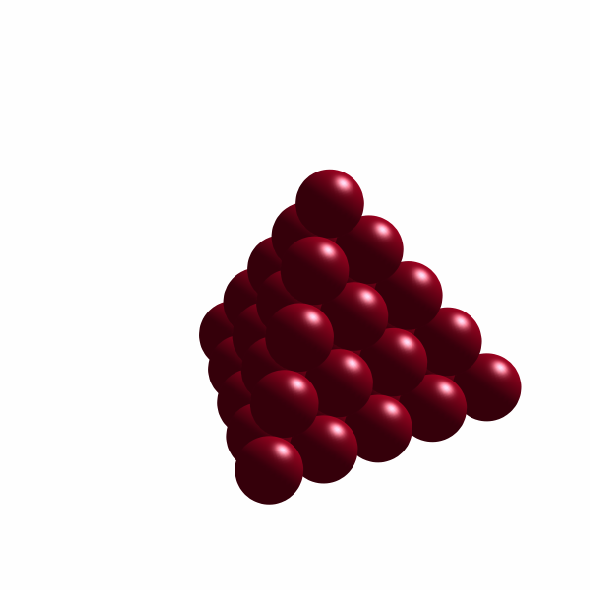}
    \hspace{25pt}
    \begin{tikzpicture}[baseline=-2pt]
        \draw[line width=1pt,->] (0,0) -- (0.5,0) node[anchor=north] {$x_1$};
        \draw[line width=1pt,->] (0,0) -- (0,0.5) node[anchor=east] {$x_2$};
    \end{tikzpicture}\hspace{1em}%
    \includegraphics[trim={36pt 47pt 25pt 38pt}, clip, width=0.27\textwidth]{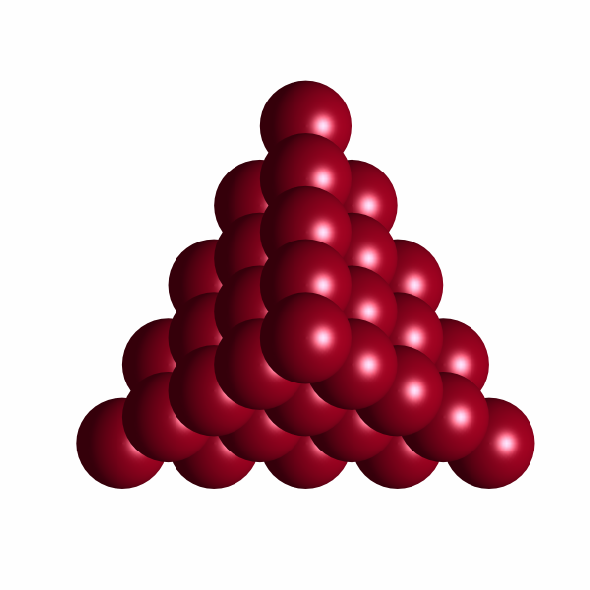}
    \caption{Triangular pyramid of \num{35} particles employed in the optimization example in \cref{section:model_validation}. The right plot shows a top view of the object.}
    \label{fig:pyramid_of_particles}
\end{figure}
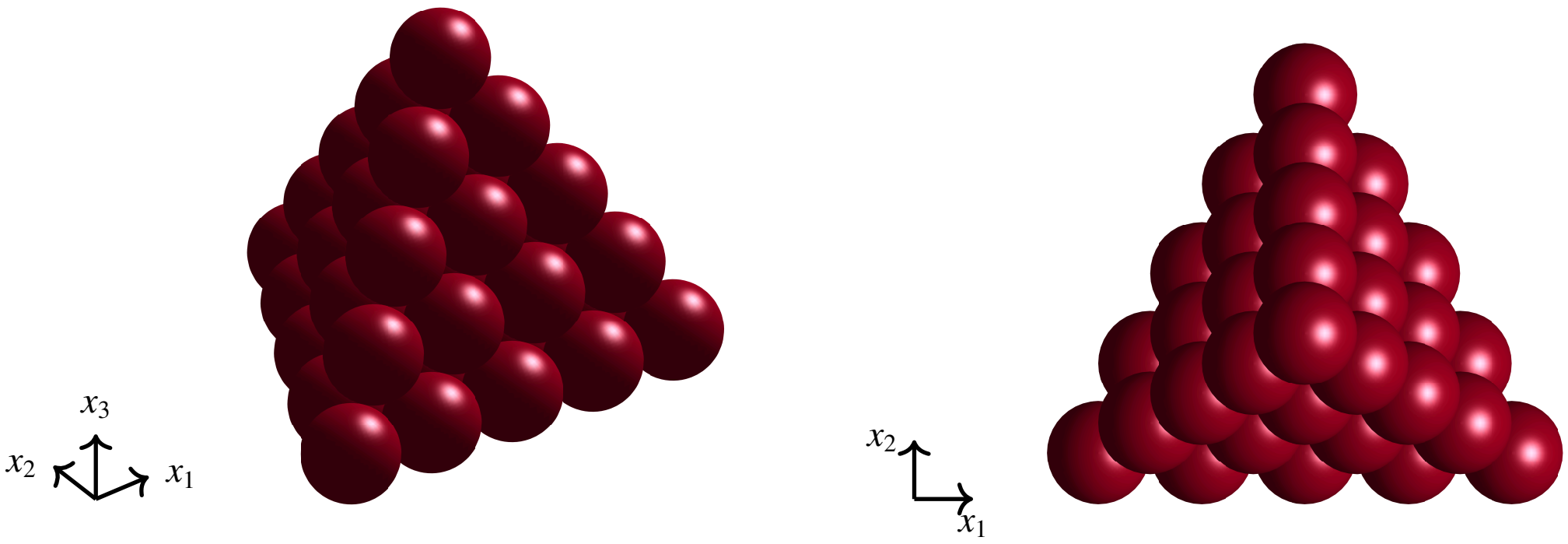

Let us first focus on the validation of the separable models derived in \cref{section:TMM_SGP_models_general}. In the following, we consider a triangular pyramid (tetrahedron) of particles with a total of $\Npar = \num{35}$ nanoparticles, as shown in \cref{fig:pyramid_of_particles}. The particles are in direct physical contact with each other, but they do not overlap. For this system, we optimize the extinction cross-section and validate the accuracy of the employed separable model based on the final design.

The diameter of each particle is given by $d_\mathrm{par} = \qty{250}{\nm}$. The object is illuminated by an $x_1$\nbhy polarized plane wave propagating in $x_3$\nbhy direction with a wavelength of $\lambda = \qty{600}{\nm}$. Based on the results of some convergence studies, the truncation order~$N_\mathrm{t}^i$ required within \ac{tmm} is set to $N_\mathrm{t}^i = 5$ for each particle $i\in\{1,\ldots,\Npar\}$. The particulate system is embedded in a vacuum with $\nbar_\mathrm{b} = \num{1}$, and the material catalogue from \cref{eq:material_catalogue} comprises the two materials $\nbar_1 = 1+1\ci$ and $\nbar_2 = 2$. Together with the linear parameterization from \cref{eq:linear_parameterization}, the admissible set~$\Uadmat$ used for optimization is therefore given by
\begin{align}\label{eq:tmm_admissible_set}
    \Uadmat = \bigl\{(1-\rho)\nbar_1 + \rho \nbar_2 \ \vert \ \rho \in [0,1]\bigr\}^{\Npar}.
\end{align}
The initial material distribution is defined by a design in which each particle of the system is assigned the second material, i.e., we have $\mb^0 \coloneqq \nbar_2 \mathbf{1}_{\Npar}$.

Let us assume that the function $S[J,\tilde\mb](\cdot)\colon\CC^\Npar\to\RR$ represents a separable model of the objective function~$J$ from {\crefname{equation}{Problem}{Problems}\cref{eq:general_obj_func_examples}} at the point of approximation $\tilde\mb\in\CC^\Npar$. Unless stated otherwise, $S$ denotes the computationally accessible first-order separable approximation~$\Sdiag$ derived in \cref{section:TMM_SGP_models_general}. In the following, we will validate this model by modifying the material property of an individual design element $i\in\{1,\ldots,\Npar\}$. To compare the original objective function with the separable model, both the objective function~$J$ and the separable model~$S$ are evaluated at the same approximation point. Subsequently, the material data associated with the $i$\nbhy th design element~$m_i$ is changed. Then, we use the relative model error~$E_\mathrm{r}^i[\tilde\mb](\cdot)\colon\CC\to\RR$ for the $i$\nbhy th design element defined by
\begin{align}\label{eq:rel_model_error}
    E_\mathrm{r}^i[\tilde\mb](m_i) \coloneqq \frac{\Bigl\vert J\bigl(\tilde\mb + (m_i - \tilde m_i) \eb^i\bigr) - S[J,\tilde\mb]\bigl(\tilde\mb + (m_i - \tilde m_i) \eb^i\bigr)\Bigr\vert}{\max_{u\in \mathcal{V}_\mathrm{ad}} J\bigl(\tilde\mb + (u - \tilde m_i) \eb^i\bigr) - \min_{u\in \mathcal{V}_\mathrm{ad}} J\bigl(\tilde\mb + (u - \tilde m_i) \eb^i\bigr)}
\end{align}
to validate the separable model~$S$ for all $\normalfont m_i \in \mathcal{V}_\mathrm{ad}$ at the approximation point~$\tilde\mb$. While $S[J,\tilde\mb]\bigl(\tilde\mb + (m_i - \tilde m_i) \eb^i\bigr)$ can be evaluated cheaply for all $m_i$, each evaluation of $J\bigl(\tilde\mb + (m_i - \tilde m_i) \eb^i\bigr)$ requires the solution of the fully coupled state equation.

\begin{figure}[!htb]
    \centering
    \definecolor{mycolor1}{rgb}{0.00000,0.44700,0.74100}%
\begin{tikzpicture}

\begin{axis}[%
width=0.38\textwidth,
height=0.38\textwidth/1.25,
at={(0\textwidth,0\textwidth)},
scale only axis,
xmin=0,
xmax=14,
table/x expr={\thisrowno{0}-1},
ymin=1.3,
ymax=1.9,
ytick={1.3,1.4,...,1.9},
xlabel={iterations},
ylabel={$\sigma_\mathrm{ext} \ /\ \unit{\mu\square\m}$},
axis background/.style={fill=white},
xmajorgrids,
ymajorgrids,
xminorgrids,
yminorgrids,
]

\addplot [color=mycolor1, line width=1.2pt]
  table[row sep=crcr]{%
1	1.89659583568573\\
2	1.41729092597961\\
3	1.34531772136688\\
4	1.32361733913422\\
5	1.31459021568298\\
6	1.31079375743866\\
7	1.30963432788849\\
8	1.3091059923172\\
9	1.30895936489105\\
10	1.30889391899109\\
11	1.30888020992279\\
12	1.30887293815613\\
13	1.30887186527252\\
14	1.30887186527252\\
};

\end{axis}
\end{tikzpicture}%
    \caption{Convergence history of the objective function for the extinction minimization example given in Problem~\eqref{eq:tmm_pyramid_extinction_min_example}.}
    \label{fig:tmm_pyramid_extinction_min_conv_hist}
\end{figure}
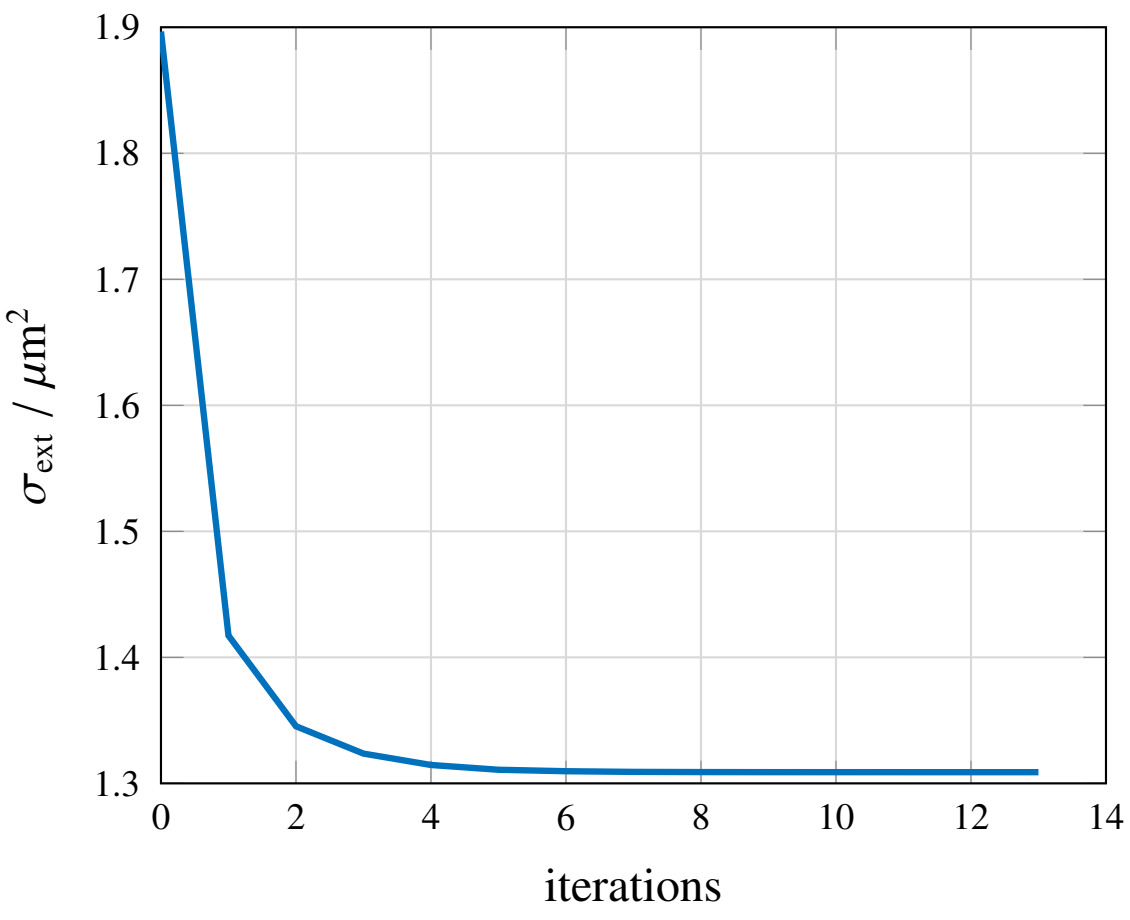

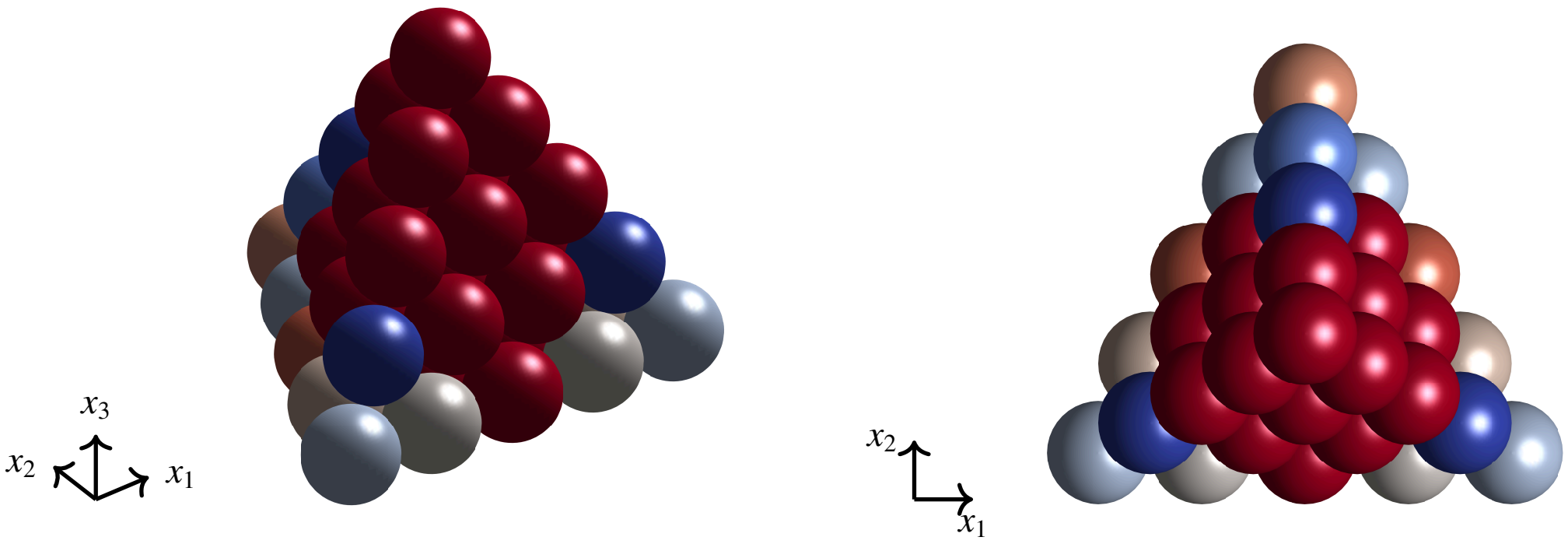
\begin{figure}[!htb]
    \centering
    \begin{tikzpicture}[baseline=-2pt]
        \draw[line width=1pt,->] (0,0) -- (0.45,0.19) node[anchor=west] {$x_1$};
        \draw[line width=1pt,->] (0,0) -- (-0.36,0.28) node[anchor=east] {$x_2$};
        \draw[line width=1pt,->] (0,0) -- (0,0.55) node[anchor=south] {$x_3$};
    \end{tikzpicture}\hspace{0.8em}%
    \includegraphics[trim={95pt 40pt 32pt 81pt}, clip, width=0.25\textwidth]{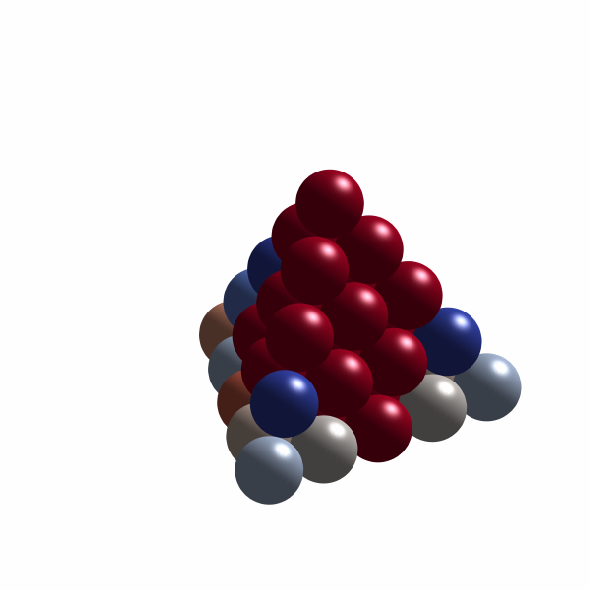}
    \hspace{25pt}
    \begin{tikzpicture}[baseline=-2pt]
        \draw[line width=1pt,->] (0,0) -- (0.5,0) node[anchor=north] {$x_1$};
        \draw[line width=1pt,->] (0,0) -- (0,0.5) node[anchor=east] {$x_2$};
    \end{tikzpicture}\hspace{1em}%
    \includegraphics[trim={36pt 47pt 25pt 38pt}, clip, width=0.27\textwidth]{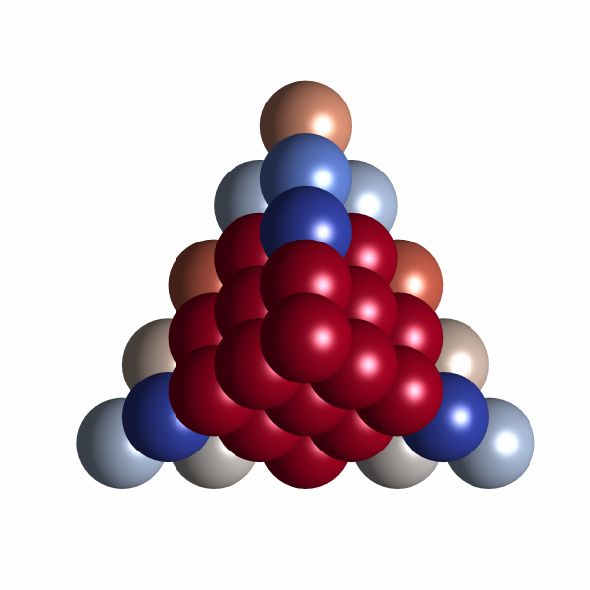}
    \caption{Final material distribution for the extinction minimization example given in Problem~\eqref{eq:tmm_pyramid_extinction_min_example}. The right plot shows a top view of the object. The colour code is the following: blue corresponds to material~1, red to material~2, and grey colours are used to indicate mixed material properties.}
    \label{fig:pyramid_of_particles_design}
\end{figure}

For the considered system of nanoparticles, we optimize the extinction cross-section given by \cref{def:extinction_TMM}, resulting in the optimization problem
\begin{align}\label{eq:tmm_pyramid_extinction_min_example}
\begin{split}
    \min_{\mb\in\Uadmat}\ & \sigma_\mathrm{ext}(\mb) \\
    \text{s.t.}\ & \Tcalb(\mb) \Acalb(\mb) = \Pcalb_\mathrm{inc}
\end{split}
\end{align}
with admissible set~$\Uadmat$ defined by \cref{eq:tmm_admissible_set}. Note that the governing equation in {\crefname{equation}{Problem}{Problems}\cref{eq:tmm_pyramid_extinction_min_example}} describes the \ac{tmm} state equation, which is defined in \cref{section:TMM}. The extinction cross-section computed for the initial material distribution $\mb^0 = \nbar_2 \mathbf{1}_{\Npar}$, i.e., for pure non-absorbing material, is given by $\sigma_\mathrm{ext}(\mb^0) = \qty[round-mode=places,round-precision=4]{1.89659583568573}{\mu \square \m}$. The aforementioned optimization problem was solved using the \ac{tmm}\nbhy \ac{sgp} approach, resulting in the convergence history depicted in \cref{fig:tmm_pyramid_extinction_min_conv_hist} and the optimized design visualized in \cref{fig:pyramid_of_particles_design}. Furthermore, let us denote the final material distribution by $\mb^*$.

The extinction cross-section decreases to $\sigma_\mathrm{ext}(\mb^*) = \qty[round-mode=places,round-precision=4]{1.30887186527252}{\mu \square \m}$ in \num{13}~iterations, which corresponds to a reduction of \qty{31.0}{\percent} with respect to the initial design. Moreover, after iteration~5 the objective function value is already within \qty{0.15}{\percent} of its final value, indicating that the subsequent design changes are relatively minor.

The optimized material distribution, visualized in \cref{fig:pyramid_of_particles_design}, exhibits a pronounced structure: the interior of the pyramid remains at the non-absorbing material~$\nbar_2$, whereas the absorbing material~$\nbar_1$ and mixtures of the two materials are assigned exclusively to particles located on the outer boundary of the object. In total, \num{3}~particles are assigned the absorbing material~$\nbar_1$, \num{20}~particles the pure material~$\nbar_2$, and the remaining \num{12}~particles retain an intermediate greyness parameter $\rho\in(0,1)$. Note that the real part of~$\nbar_1$ coincides with the refractive index of the background medium, so that the resulting shell reduces the refractive index contrast at the boundary of the object. Furthermore, the final design is symmetric with respect to the $x_2x_3$\nbhy plane, which reflects the symmetry of the object together with the direction of incidence and the polarization of the illuminating plane wave. Consequently, particles located at mirrored positions are characterized by the same material property. In the following, let us verify the accuracy of the employed separable approximation.

\paragraph{Accuracy of Separable Model:}
In order to verify the accuracy of our first-order separable model, given for the extinction cross-section in \cref{eq:sep_model_extinction_TMM}, we compare it with the true objective function for the change of the material property of a single design element. Therefore, we gradually modify the material property through the greyness parameter~$\rho$. We perform this comparison for all \num{35}~design elements based on the final design~$\mb^*$. For the true objective function, namely the extinction cross-section defined in \cref{def:extinction_TMM}, the \ac{tmm} system from \cref{eq:TMM_system} is solved for each $\rho \in [0,1]$. Here, we employ an equidistant grid of \num{201}~values of~$\rho$ per design element, so that this study requires a total of \num{7035}~solutions of the \ac{tmm} state system. In contrast, in order to compute the corresponding first-order separable approximation, the system is evaluated only once at the point of approximation given by $\tilde\mb = \mb^*$.

These fully coupled evaluations also allow a numerical test of the coordinate-wise minimality condition \cref{eq:separable_local_minimum} from \cref{def:coordinatewise_minimum} on the 201\nbhy point validation grid. In contrast to the \ac{sgp} subproblems, which minimize the surrogates on the 1001\nbhy point grid, this test uses only fully coupled objective computations. For every particle, none of the \num{201} tested material values improves $J(\mb^*)$, and the minima of the separable and the fully coupled one-particle sections coincide on the sampled candidates. Since the grid covers the entire admissible interval, this test goes beyond first-order stationarity. As only finitely many candidates are examined, it does not prove \cref{eq:separable_local_minimum} on the continuous interval. This property was reached with the first-order approximation~$\Sdiag$ and the proximal globalization of \cref{alg:SGP_alg}, i.e., without the reduction strategy of \cref{rem:coordinatewise_globalization}.

\begin{figure}[!htb]
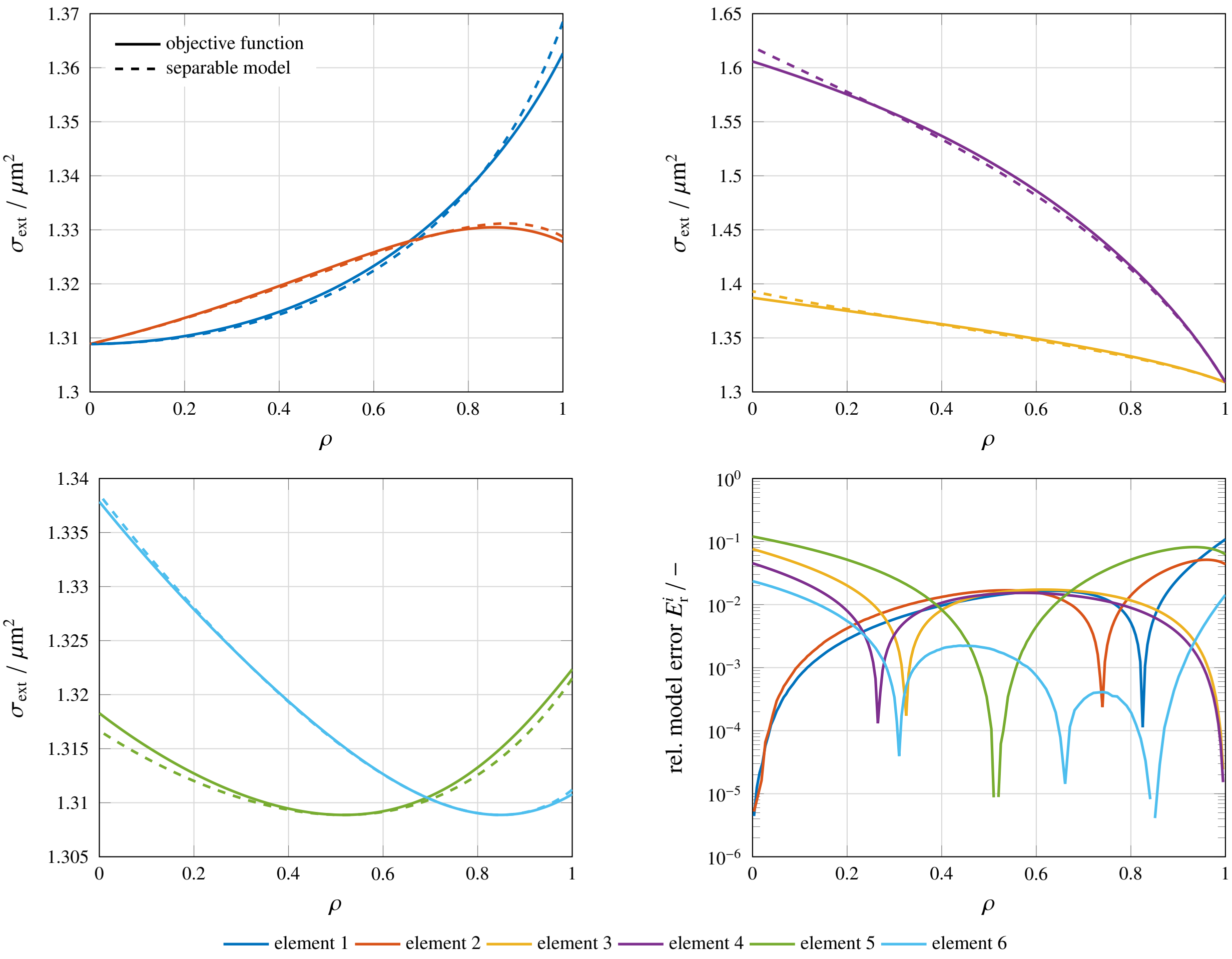

    \centering
    \input{pics/tmm_extinction_min/pyramid_opt_rho0}
    \hfill
    \input{pics/tmm_extinction_min/pyramid_opt_rho1}
    \par
    \input{pics/tmm_extinction_min/pyramid_opt_rho_grey}
    \hfill
    \input{pics/tmm_extinction_min/pyramid_opt_error_all}
    \par
    \ref*{leg:pyramid_elements}
    \caption{Accuracy of the first-order separable model based on the final design for the extinction minimization example given in Problem~\eqref{eq:tmm_pyramid_extinction_min_example}, evaluated for six design elements. These design elements are chosen such that, in each case, two of them are characterized by $\rho = 0$, $\rho = 1$, and $\rho \in (0,1)$, respectively, in the optimized material distribution. Elements~1 and~2 are characterized by $\rho = 0$, elements~3 and~4 by $\rho = 1$, and elements~5 and~6 by $\rho = \num{0.516}$ and $\rho = \num{0.848}$, respectively. The colours identify the design element. The first three plots show the true objective functions (solid curves) together with the corresponding first-order separable models (dashed curves), where there are two elements per plot. The last plot collects the associated relative model errors~$E_\mathrm{r}^i$ as defined in \cref{eq:rel_model_error}.}
    \label{fig:tmm_pyramid_model_comparison_opt}
\end{figure}

To exemplarily visualize the true objective function and the corresponding separable model, we select two design elements that are characterized within the optimized material distribution by $\rho = 0$, two by $\rho = 1$, and another two by $\rho \in (0,1)$, respectively. In the following, these are referred to as elements~1 to~6. The results are presented in \cref{fig:tmm_pyramid_model_comparison_opt}. There, the true objective functions are given by the solid curves and the corresponding separable approximations by the dashed curves, where the colour identifies the design element. The last plot collects the associated relative model errors~$E_\mathrm{r}^i$, as defined in \cref{eq:rel_model_error}.

As can be seen in \cref{fig:tmm_pyramid_model_comparison_opt}, the true objective functions can be highly nonconvex. A purely local linear gradient model cannot generally represent these nonconvex variations over the full admissible interval. For design elements whose minimum is found for $\rho \in (0,1)$, such as elements~5 and~6, the local linear model can be nearly constant and thus misses variations away from the reference point. Nevertheless, our first-order separable model provides an approximation that is not only accurate in proximity to the approximation point, but also highly accurate over the entire design interval. This can also be observed from the relative model error~$E_\mathrm{r}^i$ visualized in the last plot of \cref{fig:tmm_pyramid_model_comparison_opt}. In all cases, the separable model yields an extremely precise approximation of the physical objective function. The particle-level construction also avoids the aggregation of several dipole-based design coordinates required in the particle-homogeneous \ac{dda}\nbhy \ac{sgp} setting. That aggregation partly sacrifices the separability of the underlying \ac{dda} model, as discussed by \citet[Section~6.3]{Nees2022}. For a related ten-particle line, \citet[Section~5.4.1, Figure~5.55]{nees2025Thesis} reports a direct comparison of \ac{tmm}\nbhy \ac{sgp} and particle-based \ac{dda}\nbhy \ac{sgp} surrogates with their respective fully coupled objectives. In those numerical experiments the \ac{tmm}\nbhy based surrogate exhibits markedly smaller relative model errors, in parts of the tested material intervals by approximately one order of magnitude.

\subsection{Colour Optimization of a Photonic Crystal}\label{section:opt_photonic_crystal}
\begin{figure}[!htb]
    \centering
    \begin{tikzpicture}
        \node[anchor=south west,inner sep=0] (img) {\includegraphics[trim={105 145 109 321}, clip, width=0.55\textwidth]{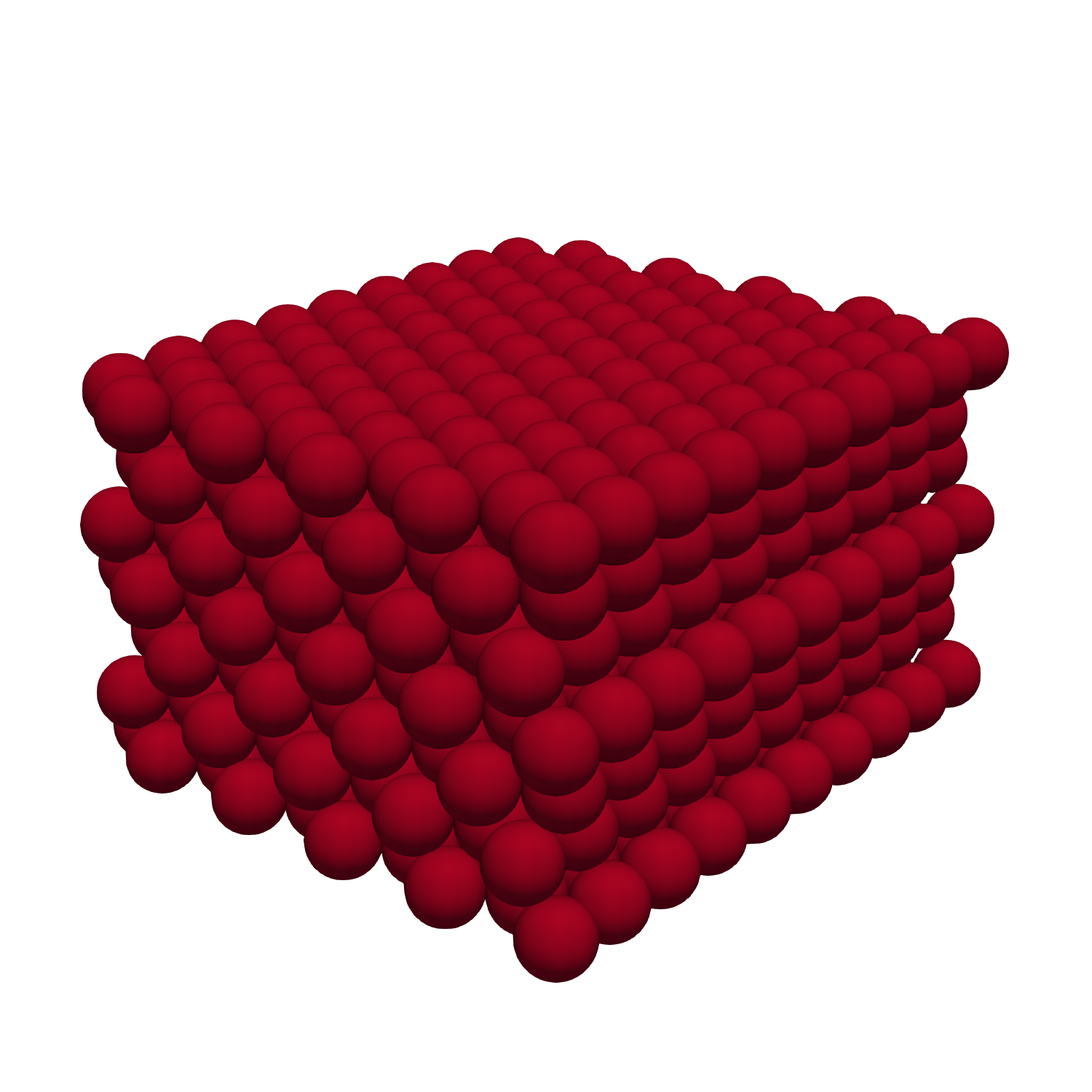}};
        \coordinate (o) at ([shift={(0.9,0.35)}]img.south west);
        \draw[line width=1pt,->] (o) -- ++(0.42,0.19) node[anchor=west] {$x_1$};
        \draw[line width=1pt,->] (o) -- ++(-0.48,0.16) node[anchor=east] {$x_2$};
        \draw[line width=1pt,->] (o) -- ++(0,0.6) node[anchor=south] {$x_3$};
    \end{tikzpicture}
    \caption{Photonic crystal composed of \num{7}~layers of closely packed spherical particles employed in the optimization examples presented in \cref{section:opt_photonic_crystal}.}
    \label{fig:photonic_crystal}
\end{figure}

In contrast to the particulate system considered in \cref{section:model_validation}, we now focus on a photonic crystal, i.e., a system of equal-sized spherical nanoparticles arranged in a hexagonal pattern. We consider \num{7}~layers of closely packed particles forming a cropped face\nbhy centred cubic (FCC) crystal with about $10\times 10$ particles per layer, which results in a total number of $\Npar = \num{627}$~particles. The object is visualized in \cref{fig:photonic_crystal}. The diameter of all particles in the photonic crystal is denoted by $d_\mathrm{par}$ and will be varied in the following, and the truncation order required within \ac{tmm} is set to $N_\mathrm{t}^i = 3$ for each particle $i\in\{1,\ldots,\Npar\}$. The object of interest is illuminated by an unpolarized incident plane wave which propagates in the $-x_3$\nbhy direction, i.e., in \cref{fig:photonic_crystal} the light propagates from top to bottom. Consequently, the state and the adjoint equation have to be solved for both polarizations in every iteration of \cref{alg:SGP_alg}.

\begin{figure}[!htb]
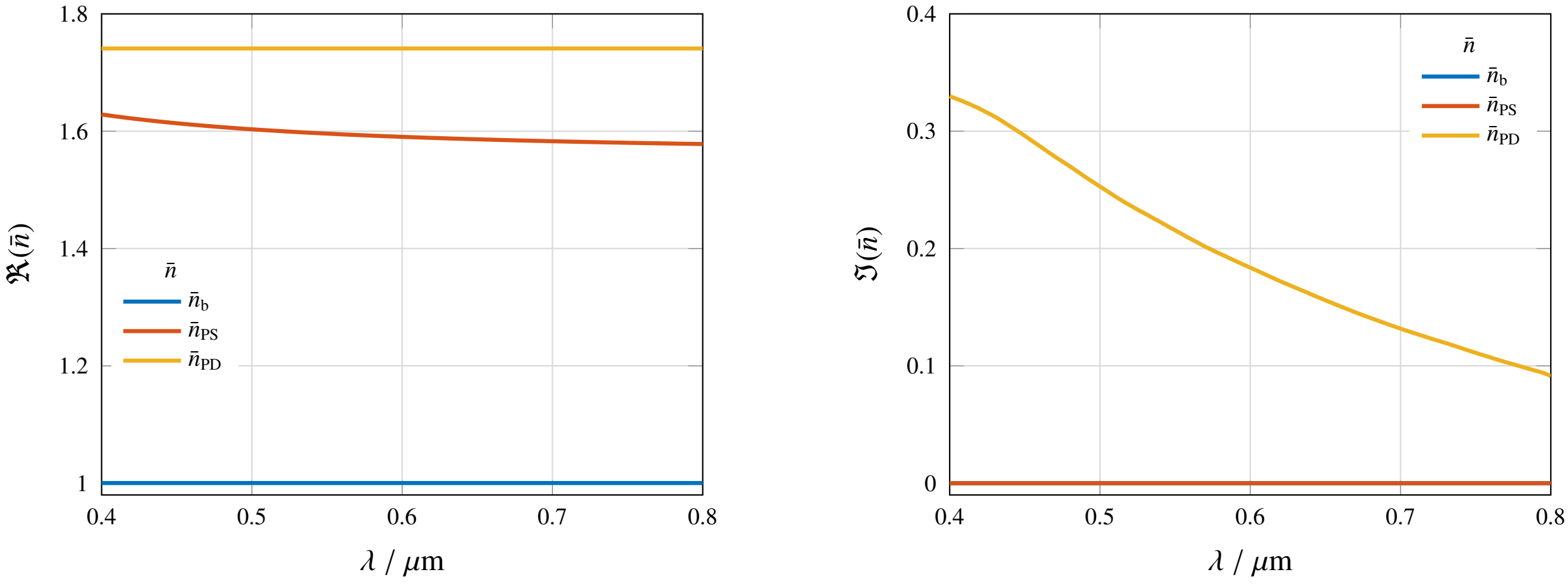

    \centering
    \input{pics/photonic_crystal/photonic_crystal_real_ref_ind}
    \hfill
    \input{pics/photonic_crystal/photonic_crystal_imag_ref_ind}
    \caption{Wavelength\nbhy dependent complex refractive indices of air (blue curve), polystyrene (red curve), and polydopamine (yellow curve). Each index is split into its real (left plot) and imaginary (right plot) parts, which are used in the optimization examples of \cref{section:opt_photonic_crystal}. The blue curve in the right plot is hidden by the red curve.}
    \label{fig:photonic_crystal_ref_index}
\end{figure}

The photonic crystal is embedded in air with a real-valued refractive index~$\nbar_\mathrm{b}$, the value of which is given by \citet[]{Ciddor1996}. Furthermore, we will focus on two materials: polystyrene (PS) with a complex refractive index~$\nbar_\mathrm{PS}$ \citep[]{sultanova2009}, and polydopamine (PD) with a complex refractive index~$\nbar_\mathrm{PD}$ \citep[]{xiao2015}. The two prescribed optical endpoint functions represent non-absorbing polystyrene and a phenomenological absorbing PS--PD composition used in the cited model. Intermediate indices are obtained by the numerical interpolation in \cref{eq:linear_parameterization}. The wavelength\nbhy dependent complex refractive indices of the two admissible materials are given by $\nbar_1 = \nbar_\mathrm{PS}$ for polystyrene and $\nbar_2 = \nbar_\mathrm{PS} + \ci \Im(\nbar_\mathrm{PD})$ for the polystyrene\nbhy polydopamine composite. A detailed explanation of the materials employed can be found in \citet[]{carina2024,bleyer2025}. The complex refractive indices of the materials are plotted in \cref{fig:photonic_crystal_ref_index} as functions of the wavelength. Together with the linear parameterization from \cref{eq:linear_parameterization}, the corresponding admissible set~$\Uadmat$ is therefore given by
\begin{align}\label{eq:photonic_crystal_admissible_set}
    \Uadmat = \bigl\{(1-\rho)\nbar_1 + \rho \nbar_2 \ \vert \ \rho \in [0,1]\bigr\}^{\Npar}.
\end{align}

We start the following optimizations with pure polystyrene, so that the whole object has non-absorbing material properties. Consequently, the initial material distribution is given by $\mb^0 \coloneqq \nbar_1 \mathbf{1}_{\Npar}$. This implies that the material property of each design element within the considered system corresponds to the first material, characterized by the complex refractive index~$\nbar_1$. For the particle system and the colour objective functions considered in this section, the colours attained along the optimization path are sufficiently saturated, so that the chroma function is differentiable (cf. \cref{rem:colour_differentiability}).

Absorbing materials have been studied as a means of improving structural colour \citep{xiao2015,Bai2019,Kawamura2016,carina2024,bleyer2025}. Here, we ask where absorbing material should be placed within a prescribed assembly, which leads to a large, wavelength-dependent and nonconvex material optimization problem.

In the following, we maximize the chroma function~$\Ccolour$ defined by \cref{eq:chroma_function} and computed using the reflectance~$r_\Mcal$ given by \cref{eq:reflectance} for the aforementioned photonic crystal. The reflectance is computed by measuring the light captured by a monitor~$\Mcal$ with position vector of its centre $\db_\Mcal = (0,0,1)^\top$ and opening angle $\theta_\Mcal = \pi$. Hence, the monitor collects all light reflected from the structure (cf. \cref{rem:reflection_transmission}). Since we consider colour property functions, we sample the visible spectral interval $\Lambda = [\qty{400}{\nm},\qty{800}{\nm}]$ at $Q=81$ wavelengths $\lambda_q=\bigl(400+5(q-1)\bigr)\,\unit{\nm}$, $q=1,\ldots,81$, and use the corresponding fixed quadrature for the colour objective. The optimization problem is given as follows:
\begin{align}\label{eq:chroma_opt_problem}
\begin{split}
    \min_{\mb\in\Uadmat}\ & J(\mb) \\
    \text{s.t.}\ & \Tcalb(\mb,\lambda) \Acalb(\mb,\lambda) = \Pcalb_\mathrm{inc}(\lambda)
\end{split}
\end{align}
for all sampled wavelengths $\lambda=\lambda_q$, $q=1,\ldots,Q$, and both polarizations of the incident field (suppressed in the notation), with $J(\mb) \coloneqq -\Ccolour(\mb)$ and admissible set of materials~$\Uadmat$ given by \cref{eq:photonic_crystal_admissible_set}. Note that the particle arrangement shown in \cref{fig:photonic_crystal} is symmetric with respect to the $x_2x_3$\nbhy plane. Since the unpolarized incident field propagates along the $-x_3$\nbhy axis and the monitor~$\Mcal$ is centred on this axis, {\crefname{equation}{Problem}{Problems}\cref{eq:chroma_opt_problem}} inherits this symmetry. The selected wavelength spacing of \qty{5}{\nm} is consistent with the sampling of the CIE~2006 colour matching data used in \cref{section:Lab_colour}. This requires \num{81}~state evaluations per reflectance spectrum and polarization. For a fresh separable model, the two polarizations and their state/adjoint equations therefore entail $81\times2\times2=324$ coupled solves. This count excludes additional full objective evaluations made when testing design updates in the globalization loop.

The monitor scattering intensity from \cref{eq:monitor_scat_integral} is approximated by the Romberg quadrature rule given in \cref{eq:directional_scat_Romberg}, where the monitor is discretized into $\tilde{n} = \num{561}$~scattering directions. Since the computation of a separable model of the objective function, i.e., of the chroma function~$\Ccolour$, requires the computation of a separable model of the reflectance, the solution of $2\tilde{n}$ adjoint problems per wavelength and incident polarization would be necessary in every major \ac{sgp} iteration. This would be far too expensive in terms of computation time. Therefore, we replace the monitor scattering intensity by its linearization from \cref{eq:linearized_monitor_scattering_function}, which yields the linearized reflectance from \cref{eq:colour_reflectance}. Based on that function, we employ the first-order separable approximation $\Sdiag\bigl[\tilde{C}^\ast[\tilde{\mb}],\tilde\mb\bigr]$ of the chroma function derived from \cref{eq:chroma_sep_model}. Hence, the negative of this model serves as the separable model of the objective function~$J$ from {\crefname{equation}{Problem}{Problems}\cref{eq:chroma_opt_problem}}.

In order to compare the initial and the optimized colours, let us additionally introduce the hue of a colour. The hue, denoted by $\hcolour$, is given by the angle measured from the positive $\acolour$\nbhy axis to the respective colour projected onto the $\acolour\bcolour$\nbhy plane, i.e.,
\begin{align}\label{eq:hue_function}
    \hcolour \coloneqq \atantwo(\bcolour,\acolour) \in (-\pi,\pi],
\end{align}
where $\atantwo$ denotes the two\nbhy argument arctangent. Note that the tuple~$(\Lcolour,\acolour,\bcolour)$ represents the Cartesian coordinates of the CIE\nbhy $\Lcolour\acolour\bcolour$ colour space, while the tuple~$(\Lcolour,\Ccolour,\hcolour)$ corresponds to the polar coordinates, respectively.

In the following, we consider photonic crystals with the particle diameters $d_\mathrm{par}\in\{\qty{200}{\nm},\allowbreak\qty{225}{\nm},\allowbreak\qty{245}{\nm},\allowbreak\qty{260}{\nm}\}$. For the initial design $\mb^0$, i.e., for pure polystyrene, these crystals exhibit a blue ($d_\mathrm{par} = \qty{200}{\nm}$), green ($d_\mathrm{par} = \qty{225}{\nm}$), yellow to beige ($d_\mathrm{par} = \qty{245}{\nm}$), and red ($d_\mathrm{par} = \qty{260}{\nm}$) colour, respectively. Thus, a broad range of colours can be obtained by changing the particle diameter~$d_\mathrm{par}$ alone. However, the chroma of the initial designs is relatively low, resulting in colours that are not particularly saturated or vivid (see \cref{tab:chroma_opt_particle}). For a vivid and more saturated colour, the ratio of the peak intensity to the incoherent scattering, i.e., the scattering that occurs at lower wavelengths in the spectrum, is decisive. Previous experiments on related material systems motivate the material choices and their optical interpretation \citep{carina2024,bleyer2025}. The present optimized particle-wise designs are assessed numerically.

For all four settings, we solve {\crefname{equation}{Problem}{Problems}\cref{eq:chroma_opt_problem}} with one design variable per particle, i.e., with $\Npar = \num{627}$~degrees of freedom. In the published particle-homogeneous \ac{dda}\nbhy \ac{sgp} formulation, several dipole variables are collapsed into a single particle-wise design coordinate. This partly sacrifices separability of the dipole-level surrogate \citep[Section~6.3]{Nees2022}.

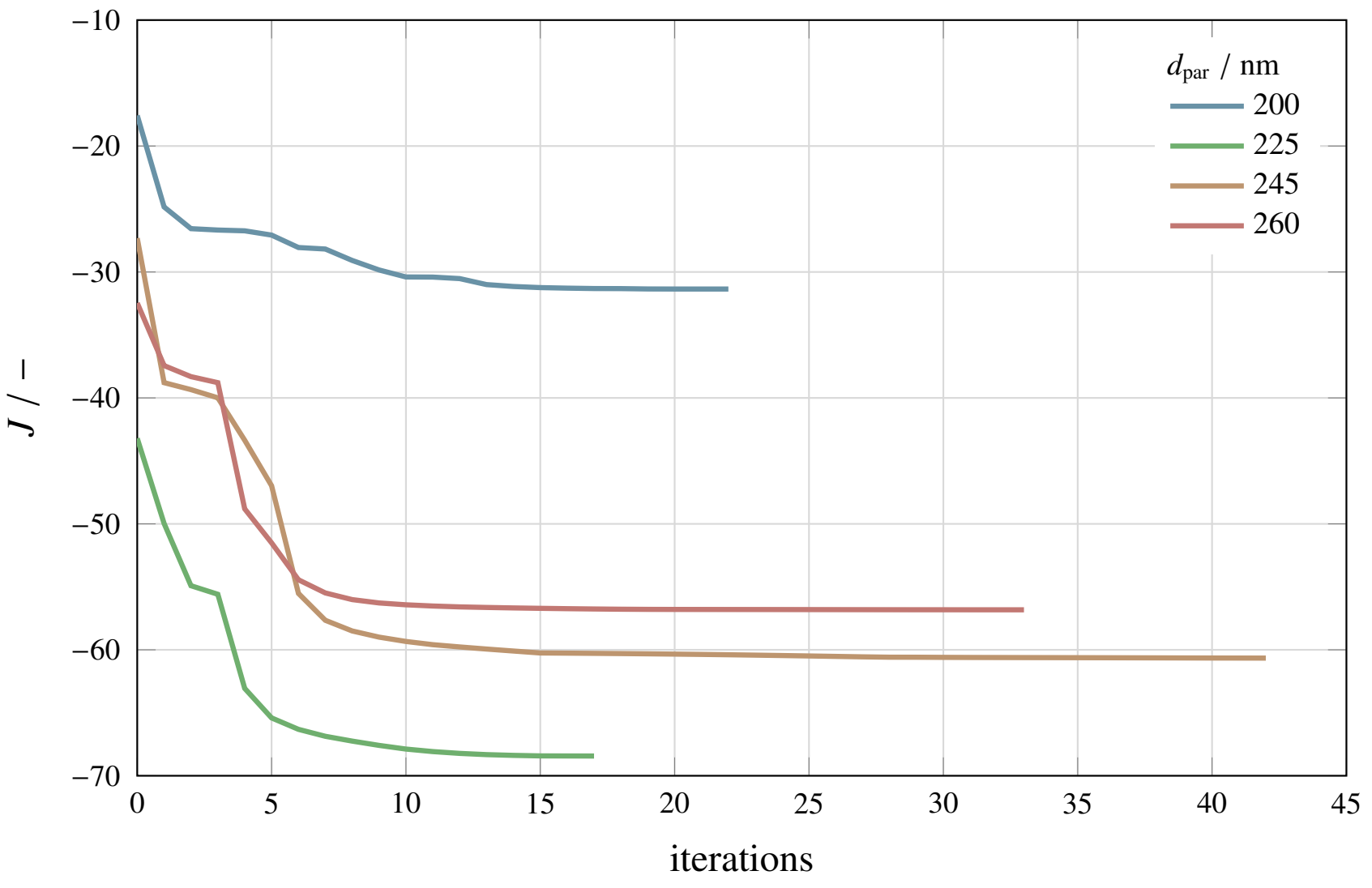
\begin{figure}[!htb]
    \centering
    \definecolor{mycolor6}{rgb}{0.41312,0.57079,0.65213}%
\definecolor{mycolor7}{rgb}{0.43507,0.68561,0.43129}%
\definecolor{mycolor8}{rgb}{0.73924,0.58401,0.43717}%
\definecolor{mycolor9}{rgb}{0.76241,0.46939,0.45250}%
\begin{tikzpicture}

\begin{axis}[%
width=0.6\textwidth,
height=0.6\textwidth/1.6,
at={(0\textwidth,0\textwidth)},
scale only axis,
xmin=0,
xmax=45,
table/x expr={\thisrowno{0}-1},
xtick={0,5,...,45},
xlabel={iterations},
ymin=-70,
ymax=-10,
ytick={-70,-60,...,-10},
ylabel={$J \ /\ -$},
axis background/.style={fill=white},
xmajorgrids,
ymajorgrids,
legend style={nodes={scale=0.8, transform shape}, legend cell align=left, align=center, draw=none},
]

\addlegendimage{white, fill=white}
\addlegendentry{\hspace{-25pt}$d_\mathrm{par} \ /\ \unit{\nm}$}

\addplot [color=mycolor6, line width=1.2pt]
  table[row sep=crcr]{%
1	-17.5714474\\
2	-24.8370209\\
3	-26.5616150\\
4	-26.6713791\\
5	-26.7330017\\
6	-27.0715866\\
7	-28.0588856\\
8	-28.1781464\\
9	-29.0844040\\
10	-29.8406925\\
11	-30.3967037\\
12	-30.4064655\\
13	-30.5258827\\
14	-31.0018311\\
15	-31.1519032\\
16	-31.2411995\\
17	-31.2841434\\
18	-31.3154202\\
19	-31.3213501\\
20	-31.3487034\\
21	-31.3538971\\
22	-31.3541775\\
23	-31.3541889\\
};
\addlegendentry{200}

\addplot [color=mycolor7, line width=1.2pt]
  table[row sep=crcr]{%
1	-43.2091980\\
2	-49.9686127\\
3	-54.9098663\\
4	-55.5912514\\
5	-63.0645256\\
6	-65.4060669\\
7	-66.3142853\\
8	-66.8657227\\
9	-67.2476807\\
10	-67.5888062\\
11	-67.8802185\\
12	-68.0799942\\
13	-68.2239990\\
14	-68.3221970\\
15	-68.3814087\\
16	-68.4223480\\
17	-68.4281082\\
18	-68.4326935\\
};
\addlegendentry{225}

\addplot [color=mycolor8, line width=1.2pt]
  table[row sep=crcr]{%
1	-27.3155575\\
2	-38.7887802\\
3	-39.3411522\\
4	-39.9927292\\
5	-43.3620300\\
6	-46.9701920\\
7	-55.5225487\\
8	-57.6555099\\
9	-58.5030098\\
10	-58.9930840\\
11	-59.3359871\\
12	-59.5918655\\
13	-59.7718124\\
14	-59.9391022\\
15	-60.0962296\\
16	-60.2507629\\
17	-60.2649689\\
18	-60.2807617\\
19	-60.2981911\\
20	-60.3180580\\
21	-60.3407631\\
22	-60.3652840\\
23	-60.3929214\\
24	-60.4227638\\
25	-60.4545898\\
26	-60.4875107\\
27	-60.5224838\\
28	-60.5581627\\
29	-60.5862541\\
30	-60.5889626\\
31	-60.5998917\\
32	-60.6072083\\
33	-60.6115646\\
34	-60.6164360\\
35	-60.6203766\\
36	-60.6248970\\
37	-60.6294899\\
38	-60.6347122\\
39	-60.6400986\\
40	-60.6461182\\
41	-60.6518021\\
42	-60.6548004\\
43	-60.6578369\\
};
\addlegendentry{245}

\addplot [color=mycolor9, line width=1.2pt]
  table[row sep=crcr]{%
1	-32.4620552\\
2	-37.4332428\\
3	-38.3101425\\
4	-38.7892265\\
5	-48.8088722\\
6	-51.5079308\\
7	-54.4397278\\
8	-55.4849472\\
9	-56.0133858\\
10	-56.2792320\\
11	-56.4259605\\
12	-56.5226936\\
13	-56.5922699\\
14	-56.6367035\\
15	-56.6744652\\
16	-56.7063599\\
17	-56.7341499\\
18	-56.7578201\\
19	-56.7772102\\
20	-56.7889709\\
21	-56.7944450\\
22	-56.7969398\\
23	-56.7991219\\
24	-56.8011742\\
25	-56.8038445\\
26	-56.8065491\\
27	-56.8096046\\
28	-56.8130112\\
29	-56.8160019\\
30	-56.8181763\\
31	-56.8197403\\
32	-56.8220367\\
33	-56.8229408\\
34	-56.8234558\\
};
\addlegendentry{260}

\end{axis}
\end{tikzpicture}%
    \caption{Convergence histories of the objective function for the chroma maximization example given in Problem~\eqref{eq:chroma_opt_problem} for the photonic crystal structure shown in \cref{fig:photonic_crystal} with different particle sizes $d_\mathrm{par} \in \{\qty{200}{\nm},\qty{225}{\nm},\qty{245}{\nm},\qty{260}{\nm}\}$.}
    \label{fig:chroma_opt_conv_hist}
\end{figure}

\begin{figure}[!htb]
    \centering
    \input{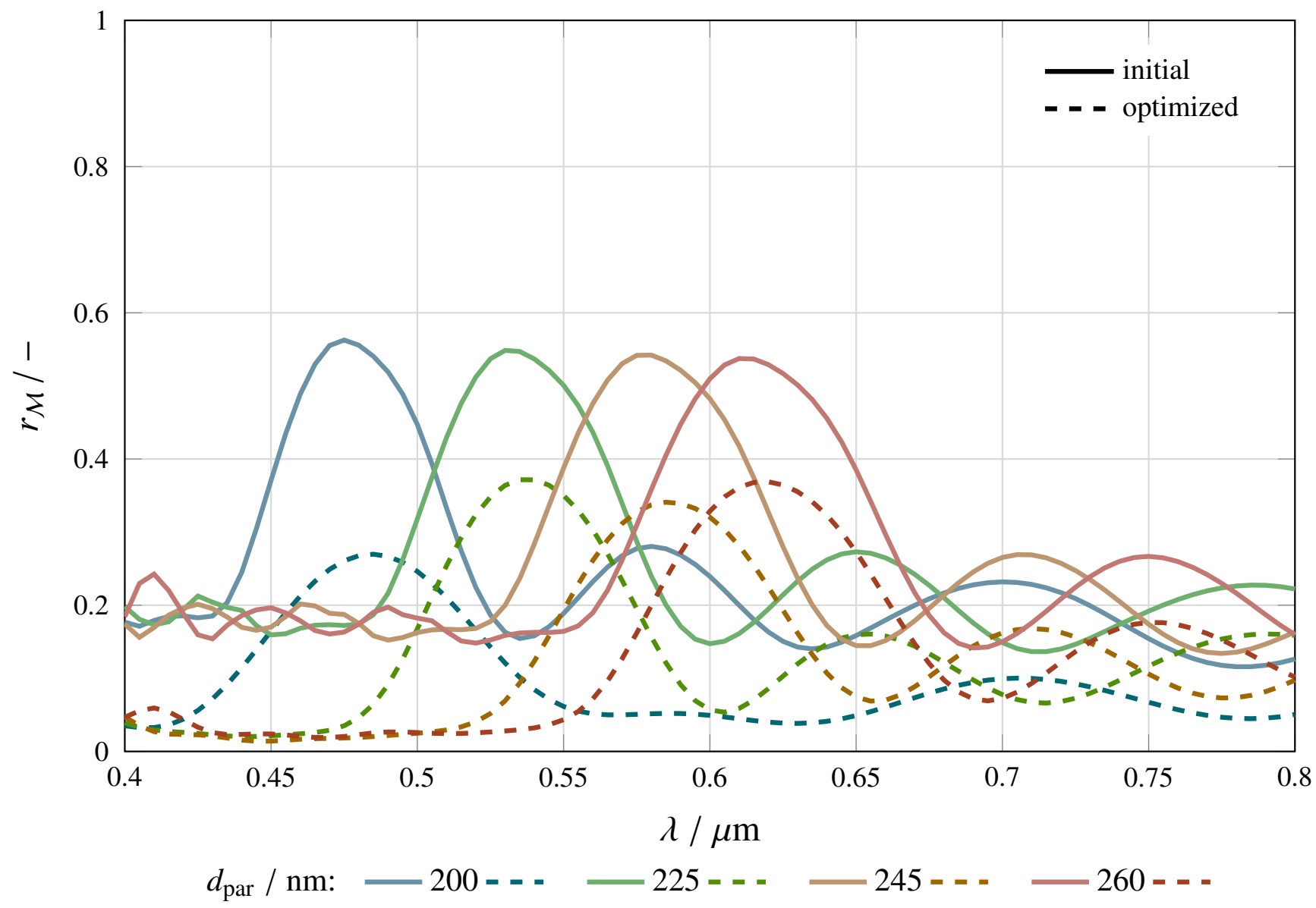}
    \caption{Initial and final reflectance spectra for the chroma maximization example given in Problem~\eqref{eq:chroma_opt_problem} for the photonic crystal structure shown in \cref{fig:photonic_crystal} with different particle sizes $d_\mathrm{par} \in \{\qty{200}{\nm},\qty{225}{\nm},\qty{245}{\nm},\qty{260}{\nm}\}$. Solid curves correspond to the initial reflectance spectra, and dashed curves to the reflectance spectra of the optimized designs. The colours of the curves correspond to the colours calculated from the respective CIE\nbhy $\Lcolour\acolour\bcolour$ values. Below the plot, the initial colour (left) and the optimized colour (right) are given for each particle diameter.}
    \label{fig:chroma_opt_spectrum}
\end{figure}

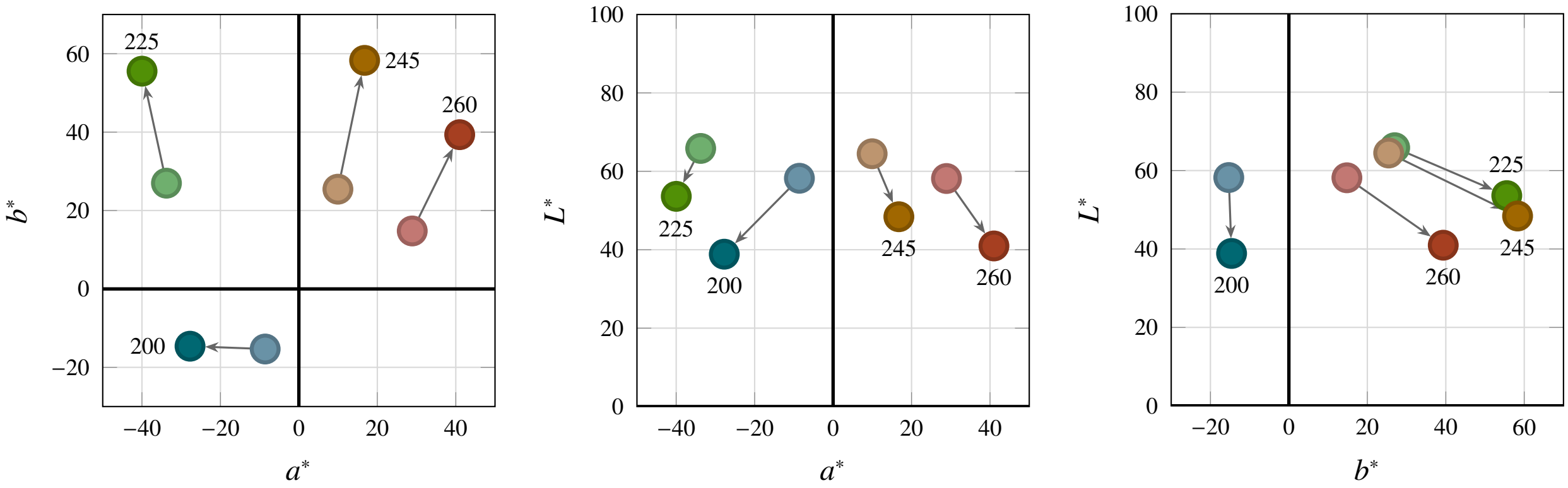
\begin{figure}[!htb]
    \centering
    \begin{tikzpicture}

\begin{axis}[%
height=0.25\textwidth,
axis equal image,
at={(0\textwidth,0\textwidth)},
scale only axis,
xmin=-50,
xmax=50,
xlabel={$\acolour$},
ymin=-30,
ymax=70,
ylabel={$\bcolour$},
axis background/.style={fill=white},
xmajorgrids,
ymajorgrids,
]
\addplot [color=black, line width=1.0pt, forget plot]
  table[row sep=crcr]{%
-80	0\\
100	0\\
};
\addplot [color=black, line width=1.0pt, forget plot]
  table[row sep=crcr]{%
0	-80\\
0	100\\
};

\draw[-{Stealth[length=4pt]}, line width=0.6pt, black!60, shorten <=4.5pt, shorten >=4.5pt] (axis cs:-8.5924,-15.3273) -- (axis cs:-27.7429,-14.6089);
\draw[-{Stealth[length=4pt]}, line width=0.6pt, black!60, shorten <=4.5pt, shorten >=4.5pt] (axis cs:-33.7568,26.9725) -- (axis cs:-39.9762,55.5422);
\draw[-{Stealth[length=4pt]}, line width=0.6pt, black!60, shorten <=4.5pt, shorten >=4.5pt] (axis cs:9.9510,25.4385) -- (axis cs:16.7614,58.2960);
\draw[-{Stealth[length=4pt]}, line width=0.6pt, black!60, shorten <=4.5pt, shorten >=4.5pt] (axis cs:28.8996,14.7850) -- (axis cs:40.9965,39.3471);

\addplot[
        scatter,%
        scatter/@pre marker code/.code={%
            \edef\temp{\noexpand\definecolor{mapped color}{rgb}{\pgfplotspointmeta}}%
            \temp
            \scope[draw=mapped color!80!black,fill=mapped color]},%
        scatter/@post marker code/.code={\endscope},%
        only marks,
        mark=*,
        point meta={TeX code symbolic={\edef\pgfplotspointmeta{\thisrow{RED},\thisrow{GREEN},\thisrow{BLUE}}}},
        mark size = 4pt,
]
table{%
x   y   RED   GREEN   BLUE
-8.5924	-15.3273	0.41312	0.57079	0.65213\\
-27.7429	-14.6089	0.00000	0.40804	0.44874\\
-33.7568	26.9725	0.43507	0.68561	0.43129\\
-39.9762	55.5422	0.31790	0.56318	0.02319\\
9.9510	25.4385	0.73924	0.58401	0.43717\\
16.7614	58.2960	0.62940	0.40209	0.00000\\
28.8996	14.7850	0.76241	0.46939	0.45250\\
40.9965	39.3471	0.65162	0.24876	0.12895\\
};

\node[font=\scriptsize, inner sep=0, anchor=east] at (axis cs:-33.5,-14.6) {200};
\node[font=\scriptsize, inner sep=0, anchor=south] at (axis cs:-40.0,60.5) {225};
\node[font=\scriptsize, inner sep=0, anchor=west] at (axis cs:21.5,58.3) {245};
\node[font=\scriptsize, inner sep=0, anchor=south] at (axis cs:41.0,44.5) {260};

\end{axis}
\end{tikzpicture}%
\hfill
    \begin{tikzpicture}

\begin{axis}[%
height=0.25\textwidth,
axis equal image,
at={(0\textwidth,0\textwidth)},
scale only axis,
xmin=-50,
xmax=50,
xlabel={$\acolour$},
ymin=0,
ymax=100,
ytick={0,20,...,100},
ylabel={$\Lcolour$},
axis background/.style={fill=white},
xmajorgrids,
ymajorgrids,
]
\addplot [color=black, line width=1.0pt, forget plot]
  table[row sep=crcr]{%
-80	0\\
100	0\\
};
\addplot [color=black, line width=1.0pt, forget plot]
  table[row sep=crcr]{%
0	-80\\
0	100\\
};

\draw[-{Stealth[length=4pt]}, line width=0.6pt, black!60, shorten <=4.5pt, shorten >=4.5pt] (axis cs:-8.5924,58.2246) -- (axis cs:-27.7429,38.8579);
\draw[-{Stealth[length=4pt]}, line width=0.6pt, black!60, shorten <=4.5pt, shorten >=4.5pt] (axis cs:-33.7568,65.8370) -- (axis cs:-39.9762,53.6143);
\draw[-{Stealth[length=4pt]}, line width=0.6pt, black!60, shorten <=4.5pt, shorten >=4.5pt] (axis cs:9.9510,64.4807) -- (axis cs:16.7614,48.4159);
\draw[-{Stealth[length=4pt]}, line width=0.6pt, black!60, shorten <=4.5pt, shorten >=4.5pt] (axis cs:28.8996,58.1636) -- (axis cs:40.9965,40.9638);

\addplot[
        scatter,%
        scatter/@pre marker code/.code={%
            \edef\temp{\noexpand\definecolor{mapped color}{rgb}{\pgfplotspointmeta}}%
            \temp
            \scope[draw=mapped color!80!black,fill=mapped color]},%
        scatter/@post marker code/.code={\endscope},%
        only marks,
        mark=*,
        point meta={TeX code symbolic={\edef\pgfplotspointmeta{\thisrow{RED},\thisrow{GREEN},\thisrow{BLUE}}}},
        mark size = 4pt,
]
table{%
x   y   RED   GREEN   BLUE
-8.5924	58.2246	0.41312	0.57079	0.65213\\
-27.7429	38.8579	0.00000	0.40804	0.44874\\
-33.7568	65.8370	0.43507	0.68561	0.43129\\
-39.9762	53.6143	0.31790	0.56318	0.02319\\
9.9510	64.4807	0.73924	0.58401	0.43717\\
16.7614	48.4159	0.62940	0.40209	0.00000\\
28.8996	58.1636	0.76241	0.46939	0.45250\\
40.9965	40.9638	0.65162	0.24876	0.12895\\
};

\node[font=\scriptsize, inner sep=0, anchor=north] at (axis cs:-27.7,33.5) {200};
\node[font=\scriptsize, inner sep=0, anchor=north] at (axis cs:-40.0,48.0) {225};
\node[font=\scriptsize, inner sep=0, anchor=north] at (axis cs:16.8,43.0) {245};
\node[font=\scriptsize, inner sep=0, anchor=north] at (axis cs:41.0,35.5) {260};

\end{axis}
\end{tikzpicture}%
\hfill
    \begin{tikzpicture}

\begin{axis}[%
height=0.25\textwidth,
axis equal image,
at={(0\textwidth,0\textwidth)},
scale only axis,
xmin=-30,
xmax=70,
xlabel={$\bcolour$},
ymin=0,
ymax=100,
ytick={0,20,...,100},
ylabel={$\Lcolour$},
axis background/.style={fill=white},
xmajorgrids,
ymajorgrids,
]
\addplot [color=black, line width=1.0pt, forget plot]
  table[row sep=crcr]{%
-80	0\\
100	0\\
};
\addplot [color=black, line width=1.0pt, forget plot]
  table[row sep=crcr]{%
0	-80\\
0	100\\
};

\draw[-{Stealth[length=4pt]}, line width=0.6pt, black!60, shorten <=4.5pt, shorten >=4.5pt] (axis cs:-15.3273,58.2246) -- (axis cs:-14.6089,38.8579);
\draw[-{Stealth[length=4pt]}, line width=0.6pt, black!60, shorten <=4.5pt, shorten >=4.5pt] (axis cs:26.9725,65.8370) -- (axis cs:55.5422,53.6143);
\draw[-{Stealth[length=4pt]}, line width=0.6pt, black!60, shorten <=4.5pt, shorten >=4.5pt] (axis cs:25.4385,64.4807) -- (axis cs:58.2960,48.4159);
\draw[-{Stealth[length=4pt]}, line width=0.6pt, black!60, shorten <=4.5pt, shorten >=4.5pt] (axis cs:14.7850,58.1636) -- (axis cs:39.3471,40.9638);

\addplot[
        scatter,%
        scatter/@pre marker code/.code={%
            \edef\temp{\noexpand\definecolor{mapped color}{rgb}{\pgfplotspointmeta}}%
            \temp
            \scope[draw=mapped color!80!black,fill=mapped color]},%
        scatter/@post marker code/.code={\endscope},%
        only marks,
        mark=*,
        point meta={TeX code symbolic={\edef\pgfplotspointmeta{\thisrow{RED},\thisrow{GREEN},\thisrow{BLUE}}}},
        mark size = 4pt,
]
table{%
x   y   RED   GREEN   BLUE
-15.3273	58.2246	0.41312	0.57079	0.65213\\
-14.6089	38.8579	0.00000	0.40804	0.44874\\
26.9725	65.8370	0.43507	0.68561	0.43129\\
55.5422	53.6143	0.31790	0.56318	0.02319\\
25.4385	64.4807	0.73924	0.58401	0.43717\\
58.2960	48.4159	0.62940	0.40209	0.00000\\
14.7850	58.1636	0.76241	0.46939	0.45250\\
39.3471	40.9638	0.65162	0.24876	0.12895\\
};

\node[font=\scriptsize, inner sep=0, anchor=north] at (axis cs:-14.6,33.5) {200};
\node[font=\scriptsize, inner sep=0, anchor=south] at (axis cs:55.5,59.0) {225};
\node[font=\scriptsize, inner sep=0, anchor=north] at (axis cs:58.3,43.0) {245};
\node[font=\scriptsize, inner sep=0, anchor=north] at (axis cs:39.3,35.5) {260};

\end{axis}
\end{tikzpicture}%
    \caption{CIE\nbhy $\Lcolour\acolour\bcolour$ colour plots for the chroma maximization example given in Problem~\eqref{eq:chroma_opt_problem} for the photonic crystal structure shown in \cref{fig:photonic_crystal} with different particle sizes $d_\mathrm{par} \in \{\qty{200}{\nm},\qty{225}{\nm},\qty{245}{\nm},\qty{260}{\nm}\}$. The CIE\nbhy $\Lcolour\acolour\bcolour$ values are calculated from the initial and optimized reflectance spectra shown in \cref{fig:chroma_opt_spectrum} and are projected onto the $\acolour\bcolour$\nbhy, $\acolour\Lcolour$\nbhy, and $\bcolour\Lcolour$\nbhy plane, respectively. The inner point colours represent the colours calculated from the CIE\nbhy $\Lcolour\acolour\bcolour$ values. The less saturated points correspond to the initial designs, whereas the more vivid ones to the optimized designs. The arrows point from the initial to the optimized colour, and the labels denote the particle diameter $d_\mathrm{par}$ in \unit{\nm}.}
    \label{fig:chroma_opt_Lab}
\end{figure}

\begin{figure}[!htbp]
    \centering
    \begin{subfigure}[b]{\textwidth}
        \centering
        \begin{tikzpicture}[baseline=(img.center)]
            \node[anchor=south west,inner sep=0] (img) {\includegraphics[trim={105 145 109 321}, clip, height=0.245\textwidth]{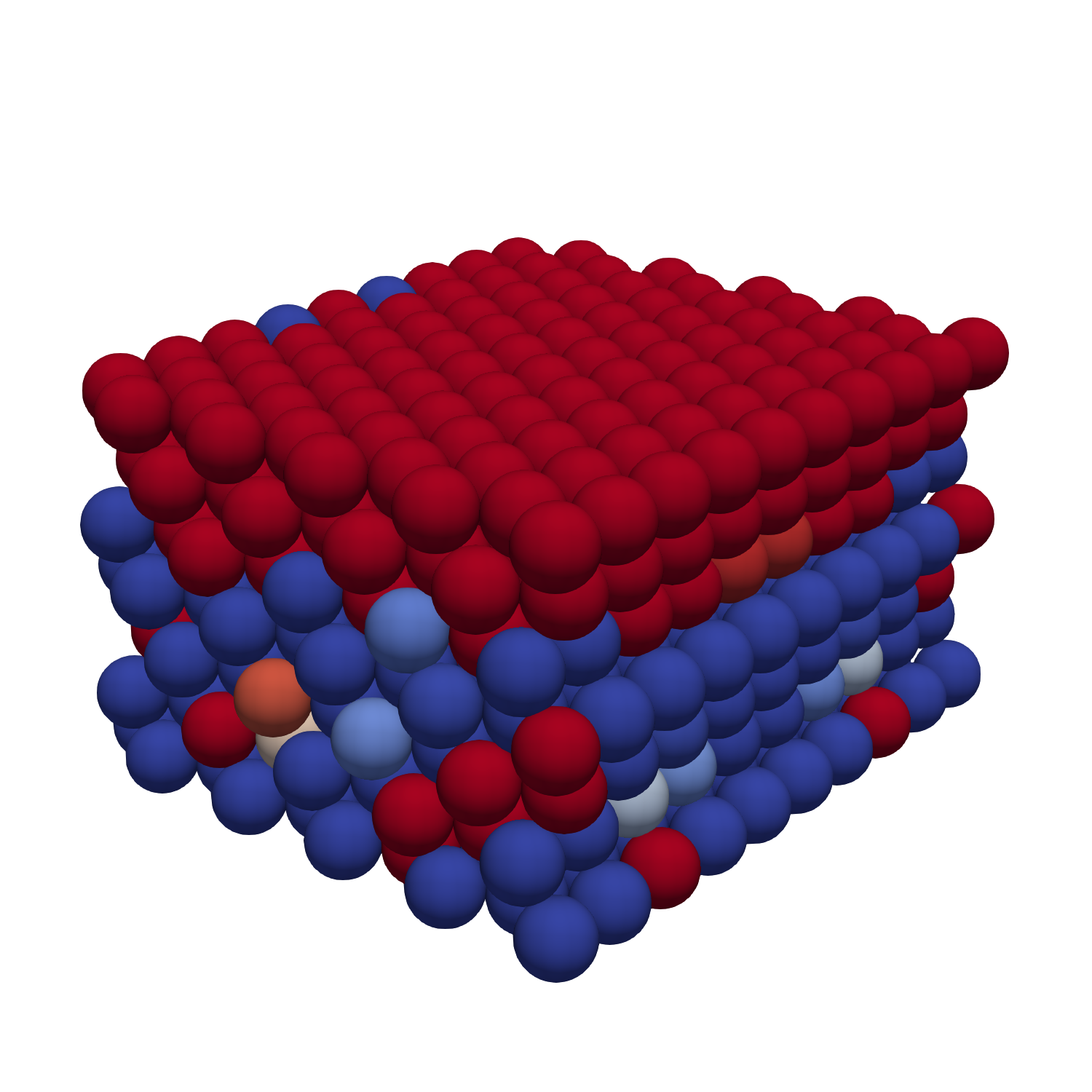}};
            \coordinate (o) at ([shift={(0.55,0.21)}]img.south west);
            \draw[line width=1pt,->] (o) -- ++(0.32,0.14) node[anchor=west] {$x_1$};
            \draw[line width=1pt,->] (o) -- ++(-0.36,0.12) node[anchor=east] {$x_2$};
            \draw[line width=1pt,->] (o) -- ++(0,0.45) node[anchor=south] {$x_3$};
        \end{tikzpicture}\hspace{2.5em}%
        \begin{tikzpicture}[baseline=(img.center)]
            \node[anchor=south west,inner sep=0] (img) {\includegraphics[trim={152 348 167 347}, clip, height=0.2\textwidth]{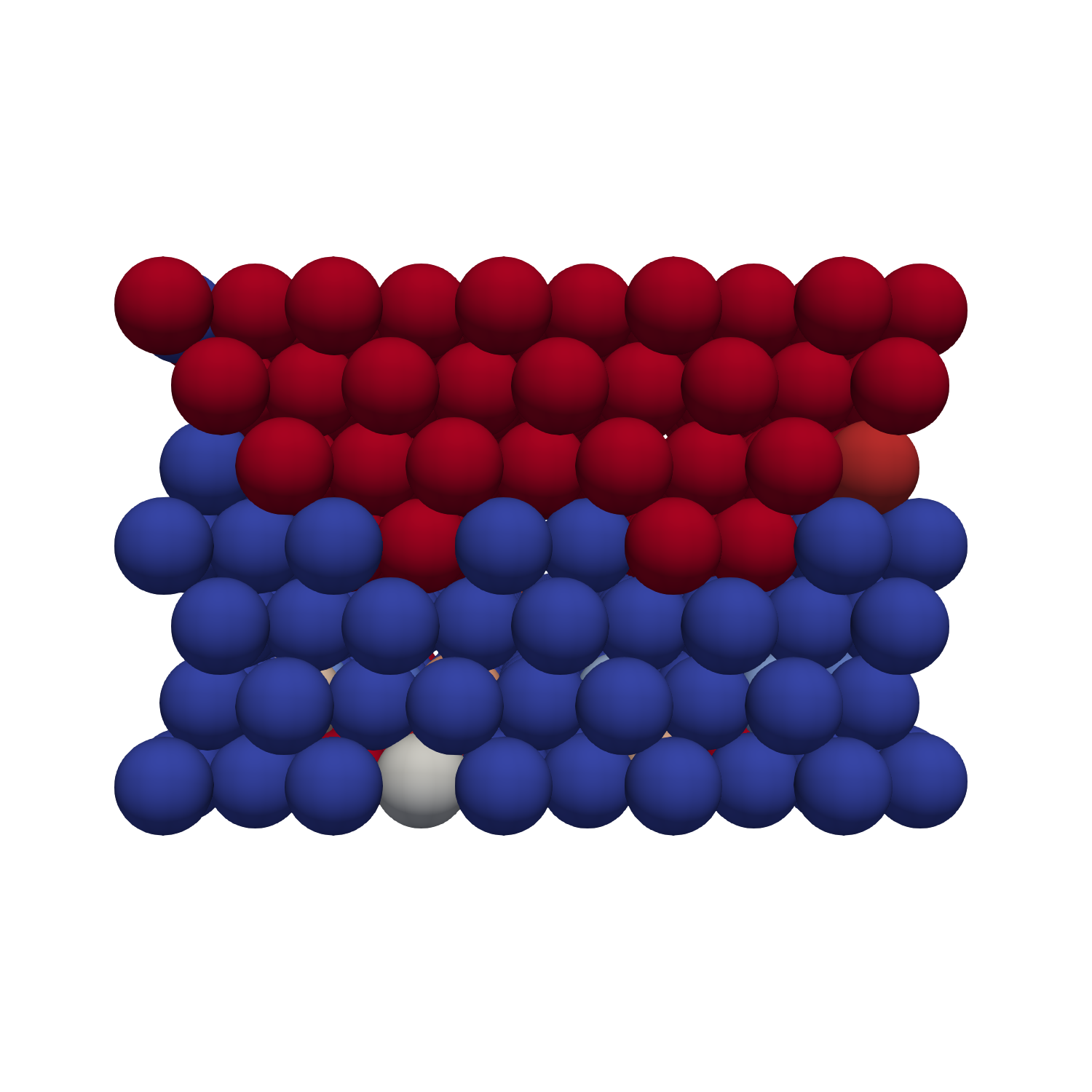}};
            \coordinate (o) at ([shift={(-0.25,0.05)}]img.south west);
            \draw[line width=1pt,->] (o) -- ++(-0.45,0) node[anchor=east] {$x_2$};
            \draw[line width=1pt,->] (o) -- ++(0,0.45) node[anchor=south] {$x_3$};
        \end{tikzpicture}
        \caption{Final design for $d_\mathrm{par} = \qty{200}{\nm}$.}
    \end{subfigure}
    \par\medskip
    \begin{subfigure}[b]{\textwidth}
        \centering
        \begin{tikzpicture}[baseline=(img.center)]
            \node[anchor=south west,inner sep=0] (img) {\includegraphics[trim={105 145 109 321}, clip, height=0.245\textwidth]{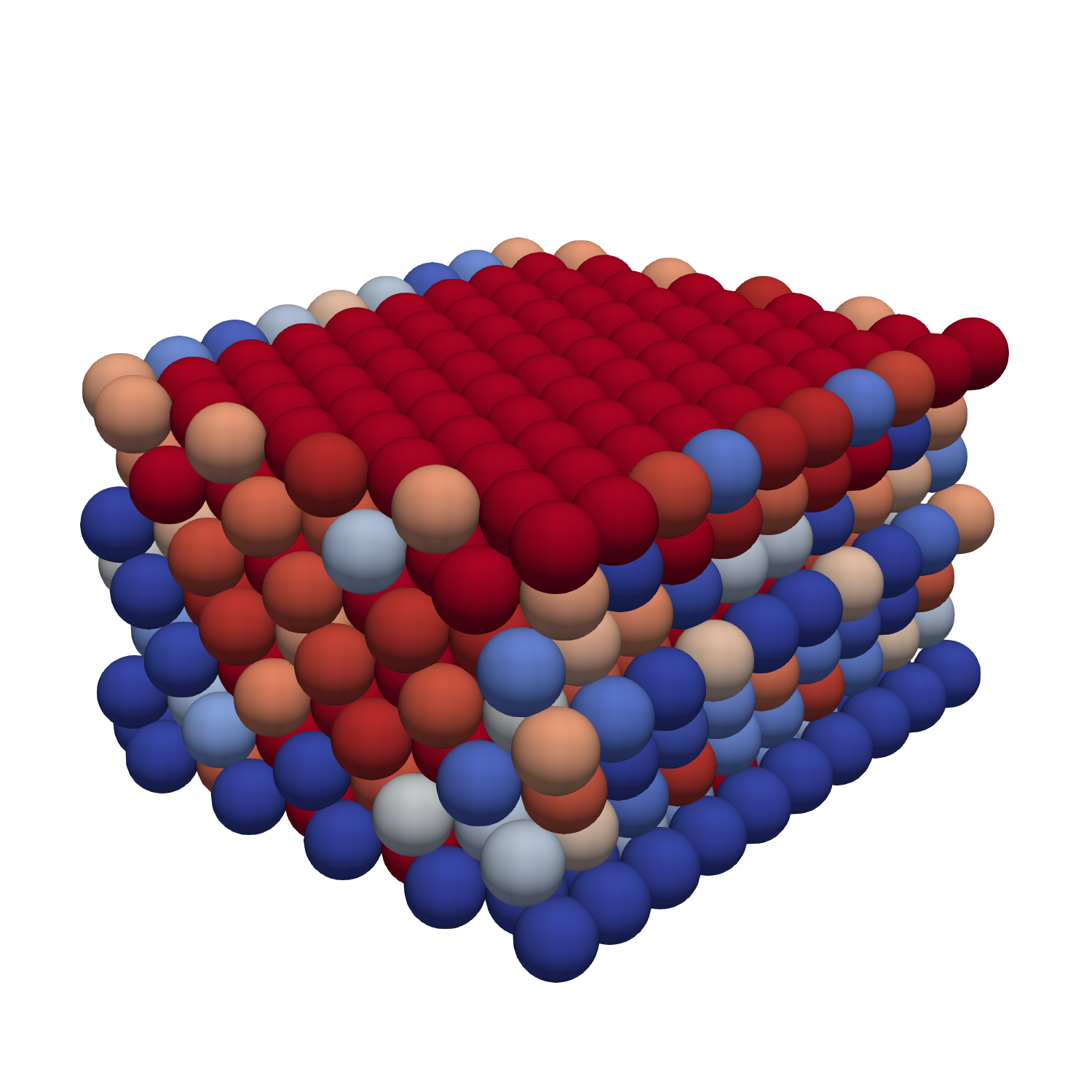}};
            \coordinate (o) at ([shift={(0.55,0.21)}]img.south west);
            \draw[line width=1pt,->] (o) -- ++(0.32,0.14) node[anchor=west] {$x_1$};
            \draw[line width=1pt,->] (o) -- ++(-0.36,0.12) node[anchor=east] {$x_2$};
            \draw[line width=1pt,->] (o) -- ++(0,0.45) node[anchor=south] {$x_3$};
        \end{tikzpicture}\hspace{2.5em}%
        \begin{tikzpicture}[baseline=(img.center)]
            \node[anchor=south west,inner sep=0] (img) {\includegraphics[trim={152 348 167 347}, clip, height=0.2\textwidth]{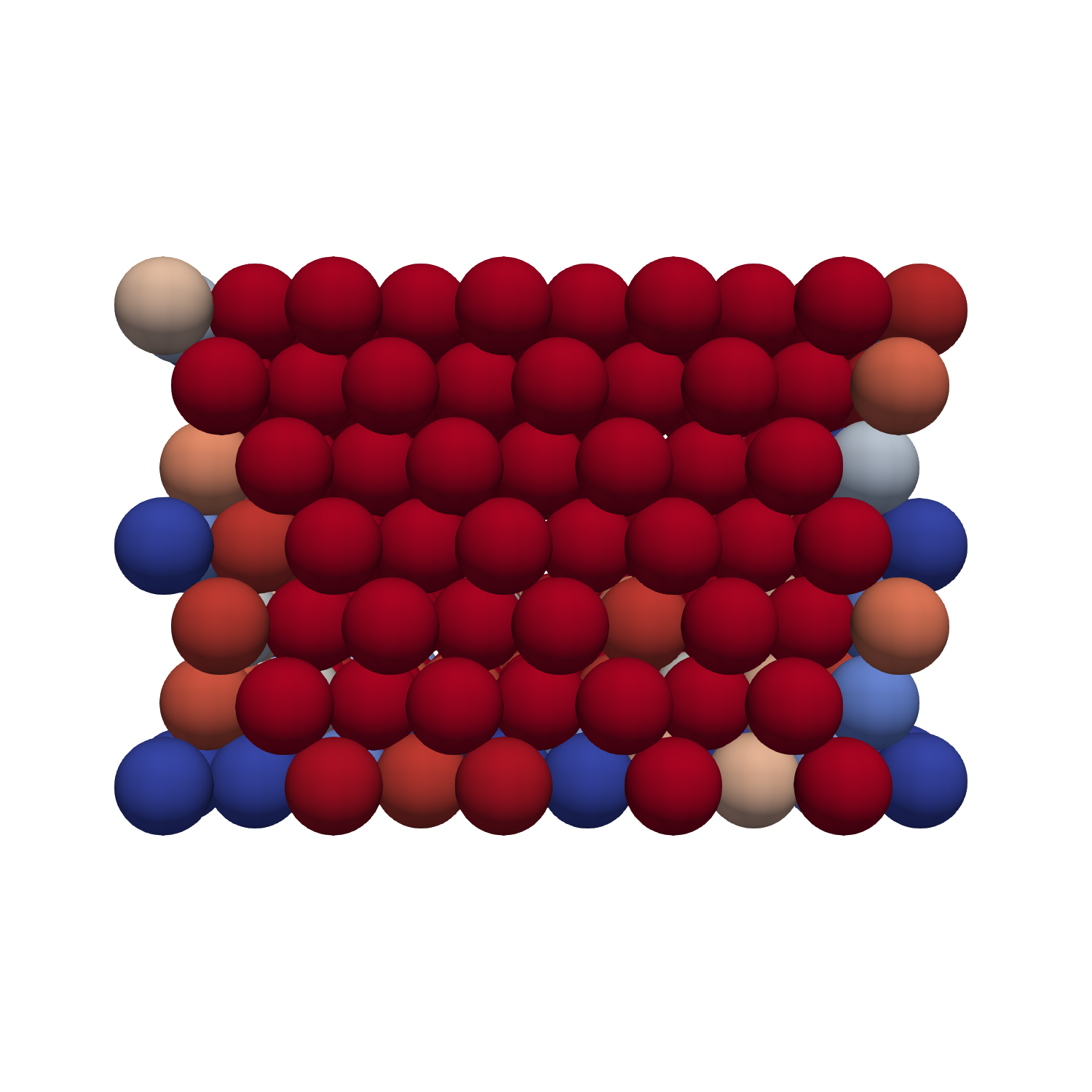}};
            \coordinate (o) at ([shift={(-0.25,0.05)}]img.south west);
            \draw[line width=1pt,->] (o) -- ++(-0.45,0) node[anchor=east] {$x_2$};
            \draw[line width=1pt,->] (o) -- ++(0,0.45) node[anchor=south] {$x_3$};
        \end{tikzpicture}
        \caption{Final design for $d_\mathrm{par} = \qty{225}{\nm}$.}
    \end{subfigure}
    \par\medskip
    \begin{subfigure}[b]{\textwidth}
        \centering
        \begin{tikzpicture}[baseline=(img.center)]
            \node[anchor=south west,inner sep=0] (img) {\includegraphics[trim={105 145 109 321}, clip, height=0.245\textwidth]{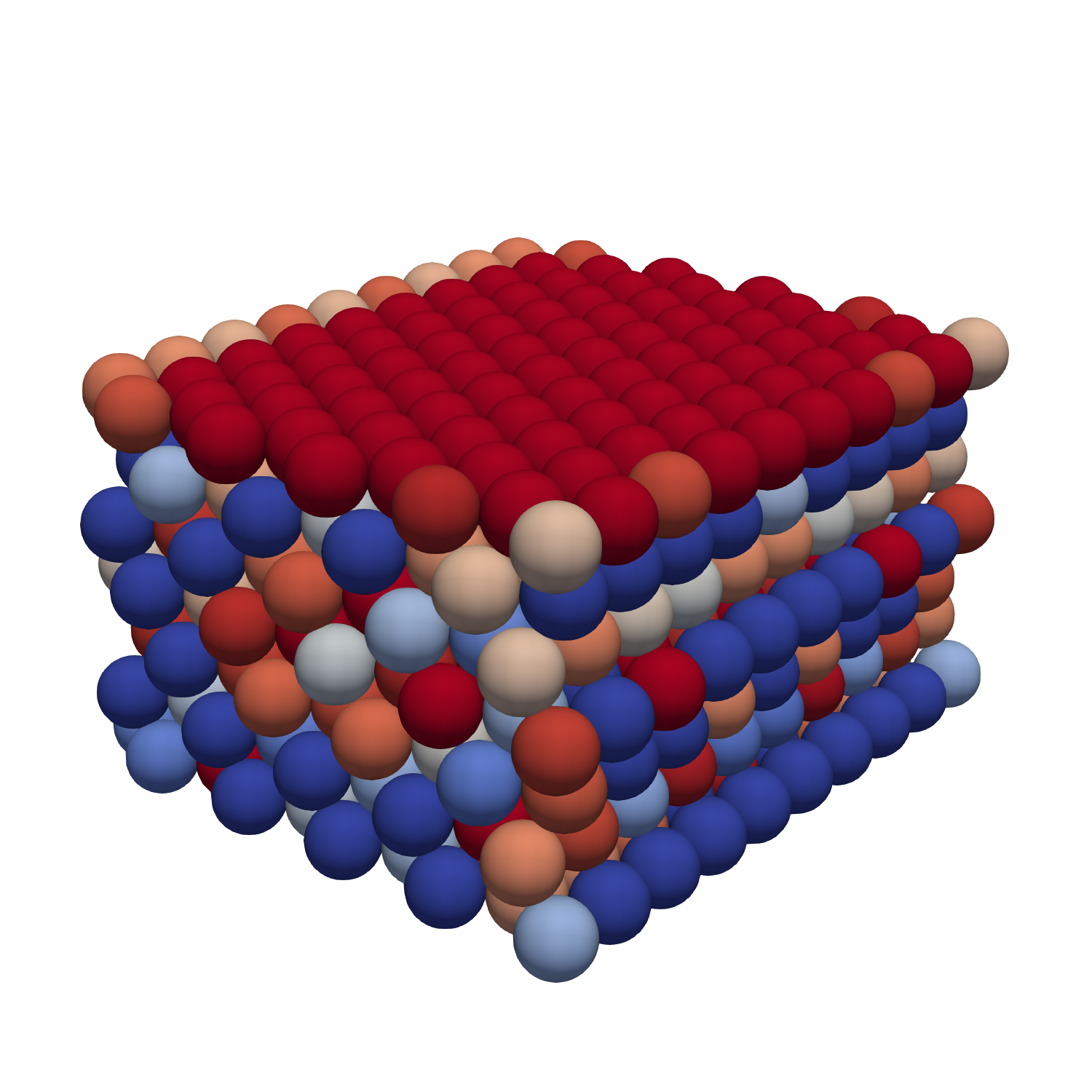}};
            \coordinate (o) at ([shift={(0.55,0.21)}]img.south west);
            \draw[line width=1pt,->] (o) -- ++(0.32,0.14) node[anchor=west] {$x_1$};
            \draw[line width=1pt,->] (o) -- ++(-0.36,0.12) node[anchor=east] {$x_2$};
            \draw[line width=1pt,->] (o) -- ++(0,0.45) node[anchor=south] {$x_3$};
        \end{tikzpicture}\hspace{2.5em}%
        \begin{tikzpicture}[baseline=(img.center)]
            \node[anchor=south west,inner sep=0] (img) {\includegraphics[trim={152 348 167 347}, clip, height=0.2\textwidth]{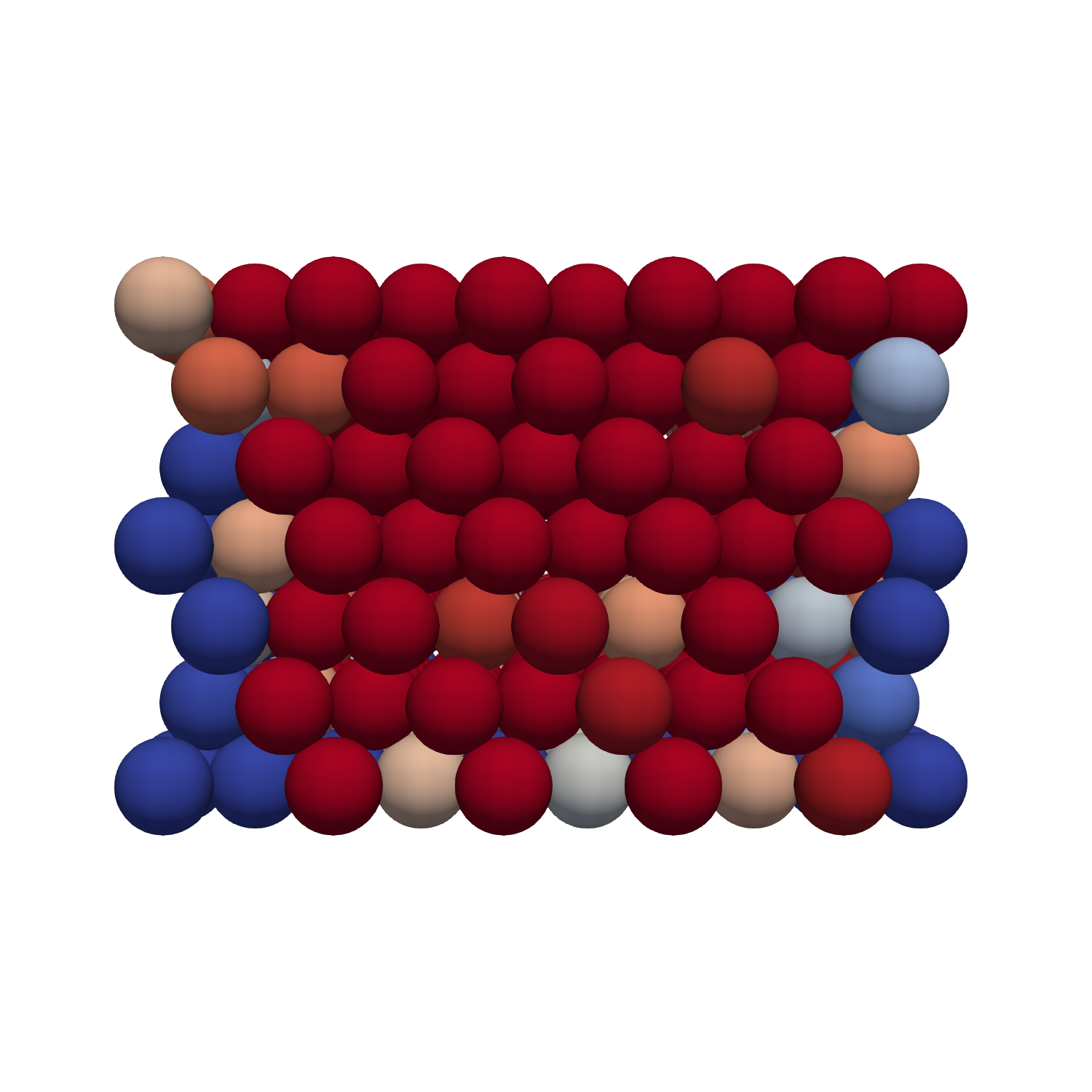}};
            \coordinate (o) at ([shift={(-0.25,0.05)}]img.south west);
            \draw[line width=1pt,->] (o) -- ++(-0.45,0) node[anchor=east] {$x_2$};
            \draw[line width=1pt,->] (o) -- ++(0,0.45) node[anchor=south] {$x_3$};
        \end{tikzpicture}
        \caption{Final design for $d_\mathrm{par} = \qty{245}{\nm}$.}
    \end{subfigure}
    \par\medskip
    \begin{subfigure}[b]{\textwidth}
        \centering
        \begin{tikzpicture}[baseline=(img.center)]
            \node[anchor=south west,inner sep=0] (img) {\includegraphics[trim={105 145 109 321}, clip, height=0.245\textwidth]{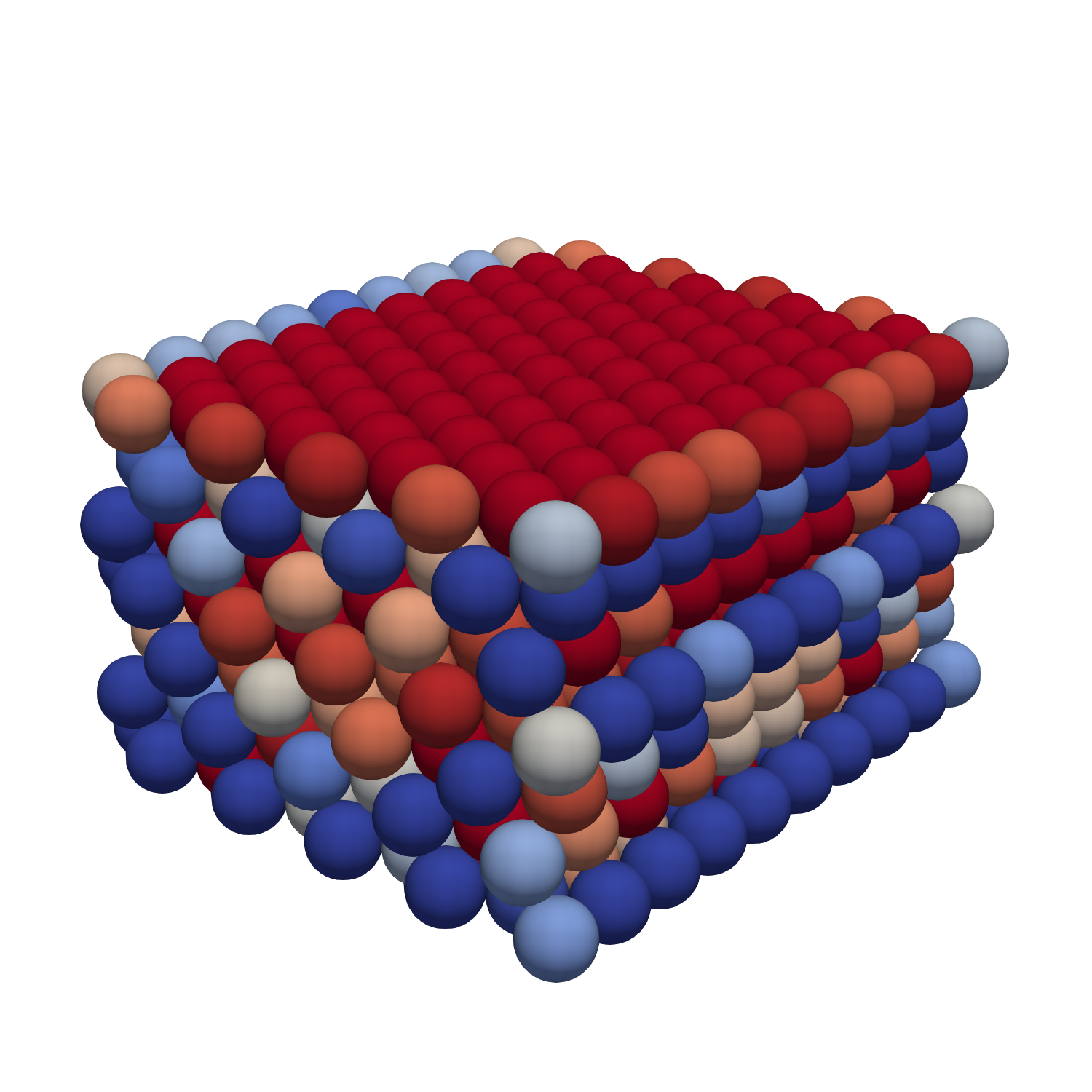}};
            \coordinate (o) at ([shift={(0.55,0.21)}]img.south west);
            \draw[line width=1pt,->] (o) -- ++(0.32,0.14) node[anchor=west] {$x_1$};
            \draw[line width=1pt,->] (o) -- ++(-0.36,0.12) node[anchor=east] {$x_2$};
            \draw[line width=1pt,->] (o) -- ++(0,0.45) node[anchor=south] {$x_3$};
        \end{tikzpicture}\hspace{2.5em}%
        \begin{tikzpicture}[baseline=(img.center)]
            \node[anchor=south west,inner sep=0] (img) {\includegraphics[trim={152 348 167 347}, clip, height=0.2\textwidth]{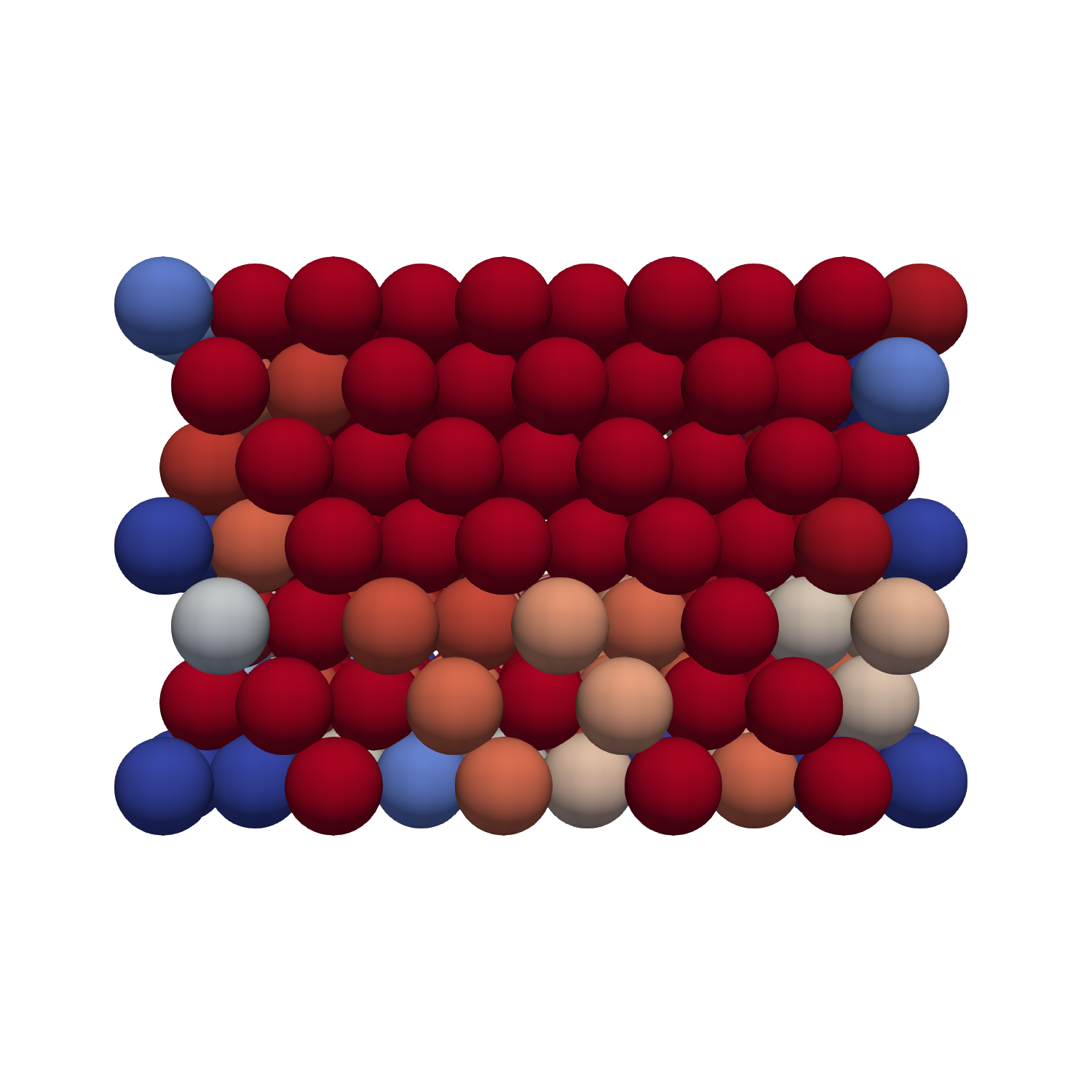}};
            \coordinate (o) at ([shift={(-0.25,0.05)}]img.south west);
            \draw[line width=1pt,->] (o) -- ++(-0.45,0) node[anchor=east] {$x_2$};
            \draw[line width=1pt,->] (o) -- ++(0,0.45) node[anchor=south] {$x_3$};
        \end{tikzpicture}
        \caption{Final design for $d_\mathrm{par} = \qty{260}{\nm}$.}
    \end{subfigure}
    \caption{Final material distribution for the chroma maximization example given in Problem~\eqref{eq:chroma_opt_problem} for the photonic crystal structure shown in \cref{fig:photonic_crystal} with different particle sizes $d_\mathrm{par} \in \{\qty{200}{\nm},\qty{225}{\nm},\qty{245}{\nm},\qty{260}{\nm}\}$. The representation on the right corresponds to a slice through the $x_2x_3$\nbhy plane of the final design. The colour code is the following: red corresponds to material~1, blue to material~2, and grey colours are used to indicate mixed material properties.}
    \label{fig:chroma_opt_design}
\end{figure}

Applying the presented \ac{tmm}\nbhy \ac{sgp} approach to the aforementioned optimization problem yields the following results. The initial and final CIE\nbhy $\Lcolour\acolour\bcolour$ colour values, together with the respective chroma value~$\Ccolour$ and the corresponding hue~$\hcolour$ of the colours are provided in \cref{tab:chroma_opt_particle}. The convergence histories and the initial and optimized reflectance spectra are presented in \cref{fig:chroma_opt_conv_hist,fig:chroma_opt_spectrum}, respectively. Additionally, \cref{fig:chroma_opt_Lab} illustrates the initial and optimized colours in the CIE\nbhy $\Lcolour\acolour\bcolour$ colour space. Furthermore, the final designs for each considered particle diameter are depicted in \cref{fig:chroma_opt_design}.

\begin{table}[!htb]
    \centering
    \sisetup{table-alignment-mode=format,table-number-alignment=center}
    \begin{tabular}{
    c
    c
    S[table-format=2.1,round-mode=places,round-precision=1]
    S[table-format=-2.1,round-mode=places,round-precision=1]
    S[table-format=-2.1,round-mode=places,round-precision=1]
    S[table-format=2.1,round-mode=places,round-precision=1]
    S[table-format=-1.2,round-mode=places,round-precision=2]
    }
        \toprule
        {$d_\mathrm{par}$ / \unit{\nm}} & design & $\Lcolour$ & $\acolour$ & $\bcolour$ & $\Ccolour$ & $\hcolour$ \\
        \midrule
        200 & initial   & 58.224565392893126 & -8.592416907493462 & -15.327299901393655 & 17.5714474 & -2.08173794519588 \\
        200 & optimized & 38.857937220757080 & -27.742851399462616 & -14.608878947500703 & 31.3541889 & -2.65690648719654 \\[1.5mm]
        225 & initial   & 65.836988497570587 & -33.756788292016452 &  26.972466629594159 & 43.2091980 & 2.46744743093146 \\
        225 & optimized & 53.614341283527040 & -39.976176862848931 & 55.542229651415944 & 68.4326935 & 2.19465064602669 \\[1.5mm]
        245 & initial   & 64.480668968943334 &  9.951002430638844 & 25.438499242821777 & 27.3155575 & 1.19791745829450 \\
        245 & optimized & 48.415929523788506 & 16.761448727765416 & 58.296028266367763 & 60.6578369 & 1.29082526364118 \\[1.5mm]
        260 & initial   & 58.163560764841932 & 28.899609595161692 & 14.785043018957555 & 32.4620552 & 0.472884561386898 \\
        260 & optimized & 40.963783159442443 & 40.996451651048041 & 39.347123966520634 & 56.8234558 & 0.764872613647453 \\
        \bottomrule
    \end{tabular}
    \caption{Initial and optimized CIE\nbhy $\Lcolour\acolour\bcolour$ colour values together with their chroma value~$\Ccolour$ for the chroma maximization example given in Problem~\eqref{eq:chroma_opt_problem} for the photonic crystal structure shown in \cref{fig:photonic_crystal} with different particle sizes $d_\mathrm{par} \in \{\qty{200}{\nm},\qty{225}{\nm},\qty{245}{\nm},\qty{260}{\nm}\}$. Additionally, the initial and final hue values~$\hcolour$ in \unit{\radian} of the resulting colours are given.}
    \label{tab:chroma_opt_particle}
\end{table}

The objective decreases for all particle diameters, see the convergence histories in \cref{fig:chroma_opt_conv_hist} and the chroma values in \cref{tab:chroma_opt_particle}. The chroma~$\Ccolour$ increases by \qtyrange{58.4}{122.1}{\percent} with respect to the initial designs within \numrange{17}{42}~iterations. In all cases, the incoherent scattering at lower wavelengths is strongly reduced, whereas the peak intensity decreases considerably less, see the reflectance spectra in \cref{fig:chroma_opt_spectrum}. Hence, the ratio of peak intensity to incoherent scattering increases, resulting in colours that are more saturated and vivid than the initial ones (see \cref{fig:chroma_opt_Lab}).

Furthermore, we can see that the problem is considerably nonlinear with regard to the particle size. This phenomenon can be explained as follows: On the one hand, the complex refractive index of polydopamine, particularly the imaginary part (see right plot in \cref{fig:photonic_crystal_ref_index}), is a function of the considered wavelength. Therefore, the polystyrene\nbhy polydopamine composite absorbs light energy more efficiently at lower wavelengths. As a consequence, the incorporation of absorber has a considerably greater effect on scattering at shorter wavelengths than at longer ones. On the other hand, the sensitivity of the human eye depends on the wavelength, which is reflected in the different weights of the wavelengths in the CIE~2006 colour matching functions \citep{CIE2006}. Consequently, the perception of colour is nonlinear in general. Accordingly, modifications in the scattering intensity at specific wavelengths may influence the colour property to a greater extent than others, thereby giving rise to the observed nonlinear behaviour. This is also reflected in the amount of absorber in the final designs, which does not depend monotonically on the particle diameter. In particular, the mean greyness parameter of the final designs is given by \num{0.50}, \num{0.25}, \num{0.29}, and \num{0.26} for $d_\mathrm{par} = \qty{200}{\nm}$, \qty{225}{\nm}, \qty{245}{\nm}, and \qty{260}{\nm}, respectively.

However, the addition of absorber also produces darker colours, as can be seen from the CIE\nbhy $\Lcolour\acolour\bcolour$ plots in \cref{fig:chroma_opt_Lab} and the corresponding lightness values~$\Lcolour$ in \cref{tab:chroma_opt_particle}. This is expected, since a considerable amount of absorbing material has been incorporated. Moreover, the outer boundary of the CIE\nbhy $\Lcolour\acolour\bcolour$ colour space in \cref{fig:CIElab_colour_space} shows that the maximum distance to the $\Lcolour$\nbhy axis, i.e., maximal chroma, is only attained for lightness values close to \num{60}. Allowing a finer spatial variation of material properties within each particle could offer additional freedom to balance chroma against lightness.

The final designs depicted in \cref{fig:chroma_opt_design} show that, for all considered particle diameters, the absorber is preferably placed at the bottom of the object, i.e., in the layers facing away from the incident light, whereas the top layer consists almost exclusively of PS. This is in accordance with the findings of experimental investigations \citep{carina2024,bleyer2025}. For $d_\mathrm{par} = \qty{200}{\nm}$, the final design is nearly binary and the absorber occupies the lower four of the seven layers. In total, \num{295}~particles are assigned the absorbing material~$\nbar_2$, while only \num{35}~particles retain an intermediate greyness parameter $\rho\in(0,1)$. For the larger particle diameters, considerably less absorber is incorporated. Here, \num{69}, \num{93}, and \num{80}~particles are assigned the absorbing material~$\nbar_2$ for $d_\mathrm{par} = \qty{225}{\nm}$, \qty{245}{\nm}, and \qty{260}{\nm}, respectively, while \num{232}, \num{242}, and \num{264}~particles retain an intermediate greyness parameter. Moreover, for these diameters, the majority of absorber is located at the lateral boundaries and in the bottom layer of the object. Although only about \qty{24}{\percent} of the particles lie on the lateral boundary, they account for \qtyrange{54}{62}{\percent} of the total greyness for these three diameters. This is attributed to the finite structure of the photonic crystal considered here. For a detailed explanation of this phenomenon in a similar setting, we refer to \citet{bleyer2025}. Finally, all final designs are symmetric with respect to the $x_2x_3$\nbhy plane, which reflects the symmetry of {\crefname{equation}{Problem}{Problems}\cref{eq:chroma_opt_problem}}.

Moreover, the resulting hue of the optimized colour can deviate from the initial hue (see \cref{tab:chroma_opt_particle} and \cref{fig:chroma_opt_Lab}). As a consequence of the nonlinearity of the colour matching functions and the manner in which the human eye perceives colour, the optimization process can also result in a different hue, which may be considered an undesired result. Nevertheless, the objective of the optimization problem was to maximize the chroma of the initial colour without any given hue direction. If the hue is to be preserved, an additional tracking term penalizing deviations from a prescribed hue could be incorporated into the objective function. However, that extension is not investigated here.

A particle-based \ac{dda}\nbhy \ac{sgp} optimization of the corresponding \num{627}\nbhy particle photonic crystal for the same four particle diameters is documented by \citet[Section~5.3.3, Figure~5.42]{nees2025Thesis}. These runs require more outer iterations than the present \ac{tmm}\nbhy \ac{sgp} runs in \cref{fig:chroma_opt_conv_hist}. Since forward solvers and surrogate constructions differ, this comparison does not isolate the effect of model accuracy, and computation times are not compared.

\section{Concluding Remarks}\label{section:conclusion}
\subsection{Summary}\label{section:summary}
We constructed exact separable models and computationally accessible first-order separable approximations for material optimization governed by the fully coupled \ac{tmm}. A change of one particle's material induces a low-rank update of the \ac{tmm} system matrix. The Sherman--Morrison--Woodbury identity provides the exact update, while a block-diagonal approximation of an auxiliary inverse produces an inexpensive surrogate without decoupling the forward or adjoint scattering problems. The construction covers general continuously differentiable functions of linear state observables and extends to several wavelengths and polarizations. Linearizing monitor scattering before applying the colour transformation reduces the adjoint workload while retaining first-order correctness and the nonlinear CIE\nbhy $\Lcolour\acolour\bcolour$ transformation.

In the \num{35}\nbhy particle extinction example, the objective decreased by \qty{31.0}{\percent} in \num{13}~iterations. Direct evaluations of the fully coupled objective at \num{201} candidate material values per particle showed high model accuracy and no improving sampled single-particle change for the final design. This finite-grid check goes beyond first-order stationarity but does not prove coordinate-wise minimality on the continuous interval. For the \num{627}\nbhy particle photonic crystal, the chroma increased by \qtyrange{58.4}{122.1}{\percent} over four particle diameters, using \numrange{17}{42}~iterations. The optimized patterns favoured absorber at the lower layers, and sometimes at lateral boundaries, broadly consistent with trends reported in \citet{carina2024,bleyer2025}. For this same geometry and these four particle diameters, the \ac{dda}\nbhy \ac{sgp} convergence histories documented in the dissertation by \citet[Section~5.3.3, Figure~5.42]{nees2025Thesis} require more outer iterations. Computation times were not compared.

These experiments illustrate how first-order correct approximations to exact separable models can retain the nonconvex particle-wise dependence that a local gradient model misses. Together with global minimization of the sampled scalar subproblems, this can permit substantial improvements in few outer \ac{sgp} iterations. The earlier published \ac{dda}\nbhy \ac{sgp} experiments demonstrate the additional possibility of escaping unfavourable stationary points \citep[Section~6.1]{Nees2022}. Intermediate refractive indices are numerical design variables rather than physical mixtures, and exact particle-by-particle placement is not assumed to be realizable. The observed spatial and layer-wise material preferences may nevertheless guide experimental studies.

\subsection{Outlook}\label{section:outlook}
Possible directions for future research are the following.

\paragraph{Multi-Material Optimization and Greyness Penalization:}
The numerical experiments use two prescribed materials and a continuous interpolation between them, sampled on a fine candidate grid when solving each particle-wise subproblem. Extension to larger finite material catalogues is possible within the graph-based \ac{sgp} parameterizations of \citet[Section~2.2]{Nees2022} and the original discrete-material framework of \citet{semmler2018}. One could also introduce greyness penalties to favour the prescribed discrete materials over intermediate interpolation values.

\paragraph{Dispersity and Synthesis Constraints:}
Particle sizes, shapes and compositions in experiments are dispersed rather than identical. Extensions combining \ac{tmm}\nbhy based material optimization with stochastic-gradient techniques \citep{pflug2020,uihlein2023,grieshammer2024a,grieshammer2024b} would therefore be valuable. Synthesis constraints also suggest replacing unrestricted particle-wise assignments with experimentally controllable quantities, such as the total absorber fraction within each layer.

\paragraph{Geometry Optimization:}
Angle-insensitive structural colours may require departures from perfect particle ordering. Joint optimization of material functions and particle positions or orientations is a possible extension. In the \ac{tmm} decomposition, position dependence enters through the coupling blocks and material dependence through the particle T\nbhy matrices. A simultaneous or nested approach would also need to enforce particle non-overlap and account for potentially expensive geometry updates.

\paragraph{General Further Research:}
Further possible extensions are:
\begin{itemize}
    \item With the \ac{tmm}\nbhy \ac{sgp} approach, it is possible to consider much larger photonic crystals in terms of the number of particles in $x_1$\nbhy{} and $x_2$\nbhy directions. This can be helpful to exclude possible boundary effects, which naturally occur due to the relatively small and finite structure considered in our examples, and to be even closer to the experimental setup.
    \item In this work, we considered spherical particles, for which the individual T\nbhy matrices are given by the Mie coefficients. The \ac{tmm}\nbhy \ac{sgp} approach can be extended to non-spherical or layered particles, e.g., core-shell particles, as long as the corresponding T\nbhy matrices are available for all admissible materials and all considered wavelengths.
\end{itemize}

\section*{CRediT authorship contribution statement}
\textbf{Nico Nees:} Conceptualization, Methodology, Software, Formal analysis, Investigation, Validation, Visualization, Writing -- original draft, Writing -- review \& editing.
\textbf{Lukas Pflug:} Methodology, Supervision, Writing -- review \& editing.
\textbf{Michael Stingl:} Conceptualization, Methodology, Supervision, Resources, Funding acquisition, Writing -- review \& editing.

\section*{Declaration of competing interest}
The authors declare that they have no known competing financial interests or personal relationships that could have appeared to influence the work reported in this paper.

\section*{Funding}
The authors acknowledge funding by the Deutsche Forschungsgemeinschaft (DFG, German Research Foundation) -- Project-ID 416229255 -- SFB 1411. The funding source had no involvement in the study design, in the collection, analysis and interpretation of data, in the writing of the report, or in the decision to submit the article for publication.

\section*{Acknowledgements}
The authors gratefully acknowledge the scientific support and HPC resources provided by the Erlangen National High Performance Computing Center (NHR@FAU) of the Friedrich-Alexander-Universität Erlangen-Nürnberg (FAU) under the NHR project b297dc -- Structure and Colour Investigation and Optimization of Particulate Systems (SCIOPS). NHR funding is provided by federal and Bavarian state authorities. NHR@FAU hardware is partially funded by the German Research Foundation (DFG) -- 440719683. Furthermore, the authors thank Florian Prohaska for his support and advice regarding the efficient implementation of the TMM\nbhy SGP optimization interface.

\section*{Data Availability Statement}
Data for this paper, including all data shown in the figures, are available at \url{https://zenodo.org/} (DOI: \url{https://doi.org/10.5281/zenodo.23588900}). 

\section*{Declaration of generative AI and AI-assisted technologies in the manuscript preparation process}
During the preparation of this work the authors used Claude (Anthropic) and ChatGPT (OpenAI) to assist with language editing, including shortening and restructuring individual sections, as well as with checking the equations for correctness. After using these tools, the authors reviewed and edited the content as needed and take full responsibility for the content of the published article.

\appendix
\section{Derivation of the TMM System}\label{appendix:TMM}
In this appendix, we summarize the derivation of the \ac{tmm} system from \cref{section:TMM}, mainly following \citet[]{egel2017,BEUTEL2024109076}. For more details, we refer to these references and to \citet[]{Mishchenko1996,Mackowski1994}. Here, the dependencies on the materials~$\mb$ and positions~$\Rb$ are not stated explicitly.

\subsection{Field Expansions}\label{appendix:TMM_expansions}
The total electromagnetic field $\Ef_{\mathrm{tot}}\colon\RR^3\to\CC^3$ is the superposition of the a-priori known incident field~$\Ef_\mathrm{inc}$, here the plane wave from \cref{section:general_electromagnetic_scattering_problem}, and the scattered fields $\Ef_{\mathrm{sca},i}\colon\RR^3\to\CC^3$ of the individual particles:
\begin{align}
    \Ef_\mathrm{tot}(\xb) = \Ef_\mathrm{inc}(\xb) + \sum_{i=1}^\Npar \Ef_{\mathrm{sca},i}(\xb).
\end{align}
The field exciting particle~$i$ consists of the incident field and the fields scattered by all other particles \citep[Equations~(1) and~(5)]{egel2017}:
\begin{align}\label{eq:TMM_incident_field_particle_i}
    \Ef_{\mathrm{inc},i}(\xb) = \Ef_{\mathrm{inc}}(\xb) + \sum_{\substack{j=1\\j\neq i}}^\Npar \Ef_{\mathrm{sca},j}(\xb).
\end{align}
A plane wave can be decomposed into its transverse electric (TE) and transverse magnetic (TM) modes with respect to a reference plane \citep[Section~2.4]{Hohenester2020}, which are computed separately in \ac{tmm}. Here, the reference plane associated with a direction is the plane spanned by that direction and the polar axis~$\eb^3$. For TE\nbhy polarized light, the electric field is perpendicular to this plane, whereas for TM\nbhy polarized light, the magnetic field is perpendicular to it.

Let $(x,\theta_\xb,\varphi_\xb)$ denote the spherical coordinates of $\xb\in\RR^3$. The vector spherical wave functions $\Psi_{\ell n m}^{(v)}\colon\RR^3\to\CC^3$ are given by \citep[Appendix~A]{egel2017}
\begin{align}
    \Psi_{1 n m}^{(v)}(\xb) &= \frac{1}{\sqrt{2n(n+1)}} \nabla \times \Bigl(\xb\, z_n^{(v)}(k_\mathrm{b} x) P_n^{\lvert m\rvert}(\cos{\theta_\xb}) \exp(\ci m \varphi_\xb)\Bigr),\\
    \Psi_{2 n m}^{(v)}(\xb) &= \frac{1}{k_\mathrm{b}} \nabla \times \Psi_{1 n m}^{(v)}(\xb),
\end{align}
where $\ell=1$ refers to TE and $\ell=2$ to TM polarization. The type~$v$ distinguishes regular waves ($v=1$), with the spherical Bessel function $z_n^{(1)}=j_n$, from outgoing waves ($v=3$), with the spherical Hankel function of the first kind $z_n^{(3)}=h_n^{(1)}$ \citep[Appendix~C.3]{Hohenester2020}. The functions~$P_n^m$ are the normalized associated Legendre functions defined in \ref{appendix:TMM_plane_wave}. With the multi-index $\nu=(\ell,n,m)$ from \cref{section:TMM}, the exciting and the scattered field of particle~$i$ are expanded as
\begin{align}
    \Ef_{\mathrm{inc},i}(\xb) = \sum_{\nu=1}^{M^i} \pf_\nu^i\, \Psi_\nu^{(1)}(\xb - \rb^i)
    \qquad\text{and}\qquad
    \Ef_{\mathrm{sca},i}(\xb) = \sum_{\nu=1}^{M^i} \af_\nu^i\, \Psi_\nu^{(3)}(\xb - \rb^i),
\end{align}
with the exciting field coefficients~$\pf_\nu^i$ and the scattering coefficients~$\af_\nu^i$. The coefficients of each particle are stored in vectors in which the TE and TM coefficients alternate, consistent with the \ac{tmm} data format proposed by \citet[]{asadova2025}, i.e.,
\begin{align}
    \af^i = \bigl({^{\mathrm{TE}}}\af^i_1,{^{\mathrm{TM}}}\af^i_1,\ldots,{^{\mathrm{TE}}}\af^i_{M^i/2},{^{\mathrm{TM}}}\af^i_{M^i/2}\bigr)^\top \in\CC^{M^i},
\end{align}
and analogously for $\pf^i$ and for the incident field coefficients~$\pf^i_\mathrm{inc}$, i.e., the expansion coefficients of $\Ef_\mathrm{inc}$ about~$\rb^i$.

\subsection{\texorpdfstring{T\nbhy Matrix}{T-Matrix} and Translation Operator}\label{appendix:TMM_Tmatrix}
The T\nbhy matrix $\Tf^i\in\CC^{M^i\times M^i}$ relates the exciting field coefficients to the scattering coefficients of particle~$i$ (cf. \citeauthor{egel2017}, \citeyear{egel2017}, Equation~(4); \citeauthor{BEUTEL2024109076}, \citeyear{BEUTEL2024109076}, Equation~(15)):
\begin{align}\label{eq:TMM_Tmatrix_particle_i_equation}
    \af^i = \Tf^i \pf^i.
\end{align}
As a consequence of \cref{eq:TMM_incident_field_particle_i}, the exciting field coefficients are obtained by translating the scattered fields of all other particles to the centre of particle~$i$ \citep[Section~2.2]{egel2017}:
\begin{align}\label{eq:TMM_incident_field_coeff_particle_i}
    \pf^i = \pf_{\mathrm{inc}}^i + \sum_{\substack{j=1\\j\neq i}}^\Npar \Cf^{i\leftarrow j} \af^j,
\end{align}
where the entries of the coupling matrix $\Cf^{i\leftarrow j}=\Cf(\rb^i-\rb^j)\in\CC^{M^i\times M^j}$ are given by
\begin{align}\label{eq:coupling_matrix}
    \Cf^{i\leftarrow j}_{\nu\xi} = \Cf_{\nu\xi}(\rb^i-\rb^j) = C_{\nu\xi}^{(3)}(\rb^i-\rb^j).
\end{align}
Here, the first index~$\nu$ refers to the receiving particle~$i$ and the second index~$\xi$ to the emitting particle~$j$. The translation operator $C_{\nu\xi}^{(v)}\colon\RR^3\to\CC$ is given by (cf. \citeauthor{egel2017}, \citeyear{egel2017}, Appendix~A2; \citeauthor{BEUTEL2024109076}, \citeyear{BEUTEL2024109076}, Appendix~C)
\begin{align}
    C_{\ell n m,\ell' n' m'}^{(v)}(\xb) = \delta_{\ell \ell'} A_{nm,n'm'}^{(v)}(\xb) + (1-\delta_{\ell \ell'}) B_{nm,n'm'}^{(v)}(\xb),
\end{align}
with $\nu = (\ell,n,m)$, $\xi = (\ell',n',m')$, the Kronecker delta~$\delta$, and the translation coefficients $A_{nm,n'm'}^{(v)},B_{nm,n'm'}^{(v)}\colon\RR^3\to\CC$ from \citet[Equations~(A.7) and (A.8)]{egel2017}, where $v=3$ corresponds to the translation of outgoing waves.

\subsection{State Equation}\label{appendix:TMM_state_equation}
Inserting \cref{eq:TMM_incident_field_coeff_particle_i} into \cref{eq:TMM_Tmatrix_particle_i_equation} yields
\begin{align}
    \af^i = \Tf^i \biggl(\pf_{\mathrm{inc}}^i + \sum_{\substack{j=1 \\ j\neq i}}^\Npar \Cf^{i\leftarrow j}\af^j\biggr),
\end{align}
which is equivalent to \citet[Equations~(4) and~(6)]{egel2017} and \citet[Equation~(16)]{BEUTEL2024109076}. Hence, the scattering coefficients of a particle are determined both by the incident field and by the fields scattered by all other particles. Multiplying by the inverse T\nbhy matrix (see \cref{rem:invertibility_Tmatrix}) and separating the unknowns leads to the state equation of particle~$i$ in \cref{eq:TMM_state_equation_i}, and combining these equations for all particles yields the \ac{tmm} system in \cref{eq:TMM_system}.

\subsection{Transformation Operator and Plane Wave Coefficients}\label{appendix:TMM_plane_wave}
Both the incident field coefficients and the far-field quantities in \cref{section:TMM_optical_properties} are expressed by the transformation operator $B_{j,\nu}\colon \RR \to\CC$ corresponding to the $j$\nbhy th polarization mode, $j\in\{1,2\}$, which is defined according to \citet[Equation~(A.13)]{egel2017} by
\begin{equation}\label{eq:transformation_operator}
\begin{aligned}
    B_{j,\nu}(\theta) & = B_{j,\ell n m}(\theta) \\
    & = -\frac{1}{\ci^{n+1}}\frac{1}{\sqrt{2n(n+1)}} (\ci\delta_{j1} + \delta_{j2})\Bigl(\delta_{\ell j} \tau_n^{\lvert m \rvert}(\theta)  + (1-\delta_{\ell j}) m \pi_n^{\lvert m \rvert}(\theta)\Bigr)\\
    & = -\frac{1}{\ci^{n+1}}\frac{1}{\sqrt{2n(n+1)}} \cdot \begin{cases}
        \ci \tau_n^{\lvert m \rvert}(\theta), & \text{if } j=1,\,\ell=1, \\
        \ci m \pi_n^{\lvert m \rvert}(\theta), & \text{if } j=1,\,\ell=2, \\
        m \pi_n^{\lvert m \rvert}(\theta), & \text{if } j=2,\,\ell=1, \\
        \tau_n^{\lvert m \rvert}(\theta), & \text{if } j=2,\,\ell=2, \\
    \end{cases}
\end{aligned}
\end{equation}
with the Kronecker delta~$\delta$ and the spherical functions $\pi_n^m\colon(0,\pi)\to\RR$ and $\tau_n^m\colon\RR\to\RR$ defined for $n\in\NN$ and $m\in\{-n,\ldots,n\}$ by
\begin{align}
    \pi_n^m(\theta) = \frac{P_n^m(\cos{\theta})}{\sin{\theta}}
    \qquad\text{and}\qquad
    \tau_n^m(\theta) = \frac{\partial P_n^m(\cos{\theta})}{\partial \theta}.
\end{align}
Note that the quotient defining $\pi_n^m$ is not defined at the poles $\theta\in\{0,\pi\}$, which correspond to exact forward and backward scattering. Since $P_n^{\lvert m\rvert}(\cos\theta)$ behaves like $\sin^{\lvert m\rvert}\theta$ there, the limits are finite for $\lvert m\rvert = 1$ and zero for $\lvert m\rvert\geq 2$, while for $m=0$ the quotient diverges. The latter case is harmless, because $\pi_n^{\lvert m\rvert}$ enters the transformation operator~$B_{j,\nu}$ only in the product $m\,\pi_n^{\lvert m\rvert}$, which vanishes for $m = 0$. In an implementation, these limits have to be evaluated explicitly instead of the quotient.
In the above formulations, $P_n^m\colon\RR\to\RR$ denote the fully normalized associated Legendre functions defined by (see \citet[Equation~(B.2)]{BEUTEL2024109076} for the associated Legendre polynomials and \url{https://de.mathworks.com/help/matlab/ref/legendre.html#f89-1002493} for the full normalization)
\begin{align}
    P_n^m(x) = \underbrace{(-1)^m\sqrt{\frac{\bigl(n+\frac{1}{2}\bigr)(n-m)!}{(n+m)!}}}_\text{normalization} \underbrace{\frac{(-1)^m}{2^n n!}(1-x^2)^\frac{m}{2} \frac{\mathrm{d}^{n+m}}{\mathrm{d}x^{n+m}}(x^2-1)^n}_\text{associated Legendre polynomials}.
\end{align}
Note that the two factors $(-1)^m$ cancel, i.e., the functions~$P_n^m$ used here do not carry the Condon\nbhy Shortley phase. They are normalized such that $\int_{-1}^{1} \bigl(P_n^m(x)\bigr)^2\,\mathrm{d}x = 1$, which is the normalization underlying the vector spherical wave functions of \citet{egel2017} and, consequently, the T\nbhy matrix convention used throughout this work.

Furthermore, the operator $B_{j,\nu}^\dagger$ is the same as $B_{j,\nu}$ but with all explicit $\ci$ in \cref{eq:transformation_operator} set to~$-\ci$.

With the transformation operator, we can state the incident field expansion coefficients~$\pf_\mathrm{inc}^i$ of the $i$\nbhy th particle for a normalized plane wave illumination. Therefore, let $(d,\theta_\db,\varphi_\db)$ be the spherical coordinates of the incident direction $\db\in\RR^3$, where $d$ describes the radial distance, $\varphi_\db$ the azimuthal angle and $\theta_\db$ the polar angle, and let $j\in\{1,2\}$ define the incident polarization state (TE for $j=1$, TM for $j=2$). Furthermore, let the corresponding radial basis vector $\hat{\rb}_\db\in\RR^3$ be defined by
\begin{align}
    \hat{\rb}_\db = \begin{pmatrix}
        \sin{\theta_\db} \cos{\varphi_\db} \\
        \sin{\theta_\db} \sin{\varphi_\db} \\
        \cos{\theta_\db}
    \end{pmatrix}.
\end{align}
Then, the incident field expansion coefficients are given by \citep[]{egel2017}
\begin{align}\label{eq:incident_coefficients_t}
    \pf_{\mathrm{inc},\nu}^i(\rb^i) & = \pf_{\mathrm{inc},\ell n m}^i(\rb^i) = 4 \exp(-\ci m \varphi_\db) \exp\bigl(\ci k_\mathrm{b} (\hat{\rb}_\db)^\top \rb^i\bigr) B_{j,\ell n m}^\dagger(\theta_\db).
\end{align}

\section{Technical Details for the Single-Particle Update}\label{appendix:single_particle_update}

This appendix records the matrix-algebra steps underlying \cref{lem:single_element_update}. All identities are evaluated at a fixed sampled wavelength $\lambda_q$. In the calculations below, the shorthand $\partial/\partial m_i$ and the partial derivatives of $\Tf^i(m_i)^{-1}$ refer to differentiation with respect to the single complex component $m_{i,q}$ relevant to this state equation. For a direction $h_i\in\MatSpace=\CC^Q$ the corresponding directional derivative is obtained by multiplying that scalar derivative by $h_{i,q}$.

\begin{proof}
    \textit{Statement~\ref{lem:single_element_update_exact}:} The state system given by \cref{eq:TMM_system_SGP} implies that
    \begin{align}
        \Lcalb^\top\Acalb\bigl(\tilde{\mb}+(m_i - \tilde m_i)\eb^i\bigr) = \Lcalb^\top \Tcalb\bigl(\tilde{\mb}+(m_i - \tilde m_i)\eb^i\bigr)^{-1}\Pcalb_\mathrm{inc}.
    \end{align}
    We can split the matrix $\Tcalb\bigl(\tilde{\mb}+(m_i - \tilde m_i)\eb^i\bigr)$ according to \cref{eq:TMM_split_matrix,eq:TMM_diagonal_part} into
    \begin{align}
        \Tcalb\bigl(\tilde{\mb}+(m_i - \tilde m_i)\eb^i\bigr) &= \Tcalb^\mathrm{pos} + \Tcalb^\mathrm{mat}\bigl(\tilde{\mb}+(m_i - \tilde m_i)\eb^i\bigr) \\
        &= \Tcalb^\mathrm{pos} + \sum_{j=1}^\Npar \Bf_j \Tf^j\Bigl((\eb^j)^\top \bigl(\tilde{\mb}+(m_i - \tilde m_i)\eb^i\bigr)\Bigr)^{-1} \Bf_j^\top.
    \end{align}
    Using the fact that
    \begin{align}
        \Tf^j\Bigl((\eb^j)^\top\bigl(\tilde{\mb}+(m_i - \tilde m_i)\eb^i\bigr)\Bigr) = \begin{cases}
            \Tf^j(\tilde m_j), & \text{if } j \neq i,\\
            \Tf^j(m_j), & \text{if } j = i,
        \end{cases}
    \end{align}
    we have
    \begin{align}
        \Tcalb\bigl(\tilde{\mb}+(m_i - \tilde m_i)\eb^i\bigr) &= \Tcalb^\mathrm{pos} + \sum_{\substack{j=1 \\j\neq i}}^\Npar \Bf_j \Tf^j(\tilde m_j)^{-1} \Bf_j^\top + \Bf_i \Tf^i(m_i)^{-1} \Bf_i^\top \\
        &= \Tcalb^\mathrm{pos} + \sum_{j=1}^\Npar \Bf_j \Tf^j(\tilde m_j)^{-1} \Bf_j^\top - \Bf_i \bigl(\Tf^i(\tilde m_i)^{-1} - \Tf^i(m_i)^{-1}\bigr) \Bf_i^\top \\
        &= \Tcalb(\tilde{\mb}) - \Bf_i \Delta\Tf^i(\tilde m_i,m_i) \Bf_i^\top,
    \end{align}
    and therefore
    \begin{align}
        \Lcalb^\top\Acalb\bigl(\tilde{\mb}+(m_i - \tilde m_i)\eb^i\bigr) = \Lcalb^\top \bigl(\Tcalb(\tilde{\mb}) - \Bf_i \Delta\Tf^i(\tilde m_i,m_i) \Bf_i^\top\bigr)^{-1}\Pcalb_\mathrm{inc}.
    \end{align}
    Applying the Sherman--Morrison--Woodbury matrix identity from \cref{eq:SMW_formula} with the matrices $\Ub = -\Bf_i$ and $\Vb = \Delta\Tf^i(\tilde m_i,m_i) \Bf_i^\top$ yields
    \begin{align}
        \Lcalb^\top\Acalb\bigl(\tilde{\mb}+(m_i - \tilde m_i)\eb^i\bigr) &= \Lcalb^\top\Bigl(\Tcalb(\tilde{\mb})^{-1} + \Tcalb(\tilde{\mb})^{-1} \Bf_i \bigl(\one_{M^i} - \Delta\Tf^i(\tilde m_i,m_i) \Bf_i^\top\Tcalb(\tilde{\mb})^{-1}\Bf_i\bigr)^{-1} \Delta\Tf^i(\tilde m_i,m_i) \Bf_i^\top \Tcalb(\tilde{\mb})^{-1}\Bigr)\Pcalb_\mathrm{inc}\\
        &= \Lcalb^\top \Acalb(\tilde{\mb}) + \Lcalb^\top \Tcalb(\tilde{\mb})^{-1} \Bf_i \Xf^i\bigl[\tilde{\mb},\Tcalb(\tilde{\mb})^{-1}\bigr](m_i) \Bf_i^\top \Acalb(\tilde{\mb}).
    \end{align}
    Then, using the adjoint problem defined by \cref{eq:general_adjoint_equation_TMM}, i.e., $\Lcalb^\top \Tcalb(\tilde{\mb})^{-1} = \Qcalb(\tilde\mb)^\top$, and the definition of~$\boldsymbol{\zeta}^i$ in \cref{eq:zeta_scalar}, we get \cref{eq:single_element_update_exact}.

    \textit{Statement~\ref{lem:single_element_update_first_order}:} At first, we have to show that $\boldsymbol{\zeta}^i[\tilde{\mb},\Kcalb](\tilde m_i) = \mathbf{0}$. This holds, since $\Delta\Tf^i(\tilde m_i,\tilde m_i) = \mathbf{0}$ and, therefore, $\Xf^i[\tilde{\mb},\Kcalb](\tilde m_i) = \mathbf{0}$. Next, we have to compute the derivative of the linear functionals. According to \cref{eq:TMM_split_matrix,eq:TMM_diagonal_part}, only the $i$\nbhy th diagonal block of the system matrix depends on~$\tilde m_i$, and therefore, we have
    \begin{align}
        \frac{\partial \Tcalb(\tilde{\mb})}{\partial \tilde m_i} = \Bf_i D_{\tilde m_i}\bigl(\Tf^i(\tilde m_i)^{-1}\bigr) \Bf_i^\top.
    \end{align}
    As discussed before \cref{eq:chain_rule_g}, the identity $\Tcalb(\mb)\Acalb(\mb) = \Pcalb_\mathrm{inc}$ holds for all~$\mb$ in a neighbourhood of~$\tilde{\mb}$, and both factors are complex differentiable in this neighbourhood. Therefore, it is possible to differentiate this identity with respect to~$m_i$ by the product rule. Evaluating the result at $\mb = \tilde{\mb}$, we get
    \begin{align}
        \frac{\partial \Tcalb(\tilde{\mb})}{\partial \tilde m_i}\Acalb(\tilde{\mb}) + \Tcalb(\tilde{\mb})\frac{\partial \Acalb(\tilde{\mb})}{\partial \tilde m_i} = \mathbf{0}
        \quad\Longrightarrow\quad
        \frac{\partial \Acalb(\tilde{\mb})}{\partial \tilde m_i} = - \Tcalb(\tilde{\mb})^{-1} \Bf_i D_{\tilde m_i} \bigl(\Tf^i(\tilde m_i)^{-1}\bigr) \Bf_i^\top \Acalb(\tilde{\mb}).
    \end{align}
    Note that the same result is obtained by applying the implicit function theorem to the function $(\mb,\yb)\mapsto\Tcalb(\mb)\yb-\Pcalb_\mathrm{inc}$. Hence, together with $\Lcalb^\top \Tcalb(\tilde{\mb})^{-1} = \Qcalb(\tilde\mb)^\top$, we get
    \begin{align}
        \frac{\partial\bigl(\Lcalb^\top\Acalb(\tilde{\mb})\bigr)}{\partial\tilde m_i} = -\Qcalb(\tilde{\mb})^\top\Bf_i D_{\tilde m_i}\bigl(\Tf^i(\tilde m_i)^{-1}\bigr)\Bf_i^\top\Acalb(\tilde{\mb}).
    \end{align}
    Finally, we have to compute the derivative of~$\boldsymbol{\zeta}^i[\tilde{\mb},\Kcalb](m_i)$ with respect to~$m_i$. For the sake of readability, we abbreviate the matrix that is inverted in \cref{eq:X_matrix} by
    \begin{align}\label{eq:X_abbreviation}
        \Yf^i(m_i) = \one_{M^i} -\Delta\Tf^i(\tilde m_i,m_i) \Bf_i^\top \Kcalb \Bf_i \in \CC^{M^i\times M^i},
    \end{align}
    which is invertible by assumption, so that $\Xf^i[\tilde{\mb},\Kcalb](m_i) = \Yf^i(m_i)^{-1} \Delta\Tf^i(\tilde m_i,m_i)$.
    For the derivative of the inverse we use the identity $\frac{\partial (\Yf^i(m_i))^{-1}}{\partial m_i} = -(\Yf^i(m_i))^{-1}\frac{\partial \Yf^i(m_i)}{\partial m_i}(\Yf^i(m_i))^{-1}$ together with $\frac{\partial \Yf^i(m_i)}{\partial m_i} = -\frac{\partial \Delta\Tf^i(\tilde m_i,m_i)}{\partial m_i}\Bf_i^\top \Kcalb \Bf_i$, which gives
    \begin{align}
        \frac{\partial \Yf^i(m_i)^{-1}}{\partial m_i} = \Yf^i(m_i)^{-1} \frac{\partial \Delta\Tf^i(\tilde m_i,m_i)}{\partial m_i} \Bf_i^\top \Kcalb \Bf_i  \Yf^i(m_i)^{-1}.
    \end{align}
    With the product rule, this yields
    \begin{align}
        & \frac{\partial \Xf^i[\tilde{\mb},\Kcalb](m_i)}{\partial m_i} = \Yf^i(m_i)^{-1} \frac{\partial \Delta\Tf^i(\tilde m_i,m_i)}{\partial m_i} + \Yf^i(m_i)^{-1} \frac{\partial \Delta\Tf^i(\tilde m_i,m_i)}{\partial m_i} \Bf_i^\top \Kcalb \Bf_i  \Yf^i(m_i)^{-1} \Delta\Tf^i(\tilde m_i,m_i).
    \end{align}
    Furthermore, it holds
    \begin{align}
        \frac{\partial \Delta\Tf^i(\tilde m_i,m_i)}{\partial m_i} = \frac{\partial}{\partial m_i} \bigl(\Tf^i(\tilde m_i)^{-1} - \Tf^i(m_i)^{-1}\bigr) = -\frac{\partial \bigl(\Tf^i(m_i)^{-1}\bigr)}{\partial m_i} = -D_{m_i} \bigl(\Tf^i(m_i)^{-1}\bigr).
    \end{align}
    Hence, we have
    \begin{align}
        \frac{\partial \boldsymbol{\zeta}^i[\tilde{\mb},\Kcalb](m_i)}{\partial m_i} &= \Qcalb(\tilde{\mb})^\top \Bf_i \frac{\partial \Xf^i[\tilde{\mb},\Kcalb](m_i)}{\partial m_i} \Bf_i^\top \Acalb(\tilde{\mb}) \\
        &= -\Qcalb(\tilde{\mb})^\top \Bf_i \Yf^i(m_i)^{-1} D_{m_i} \bigl(\Tf^i(m_i)^{-1}\bigr) \bigl(\one_{M^i} + \Bf_i^\top \Kcalb \Bf_i\Xf^i[\tilde{\mb},\Kcalb](m_i)\bigr) \Bf_i^\top \Acalb(\tilde{\mb}).
    \end{align}
    Evaluating at the point of approximation, we have $\Delta\Tf^i(\tilde m_i,\tilde m_i) = \mathbf{0}$ and hence $\Yf^i(\tilde m_i) = \one_{M^i}$ as well as $\Xf^i[\tilde{\mb},\Kcalb](\tilde m_i) = \mathbf{0}$. From this, we get that
    \begin{align}
        \frac{\partial \boldsymbol{\zeta}^i[\tilde{\mb},\Kcalb](m_i)}{\partial m_i}\biggr\vert_{m_i=\tilde m_i} = -\Qcalb(\tilde{\mb})^\top \Bf_i D_{\tilde m_i} \bigl(\Tf^i(\tilde m_i)^{-1}\bigr) \Bf_i^\top \Acalb(\tilde{\mb}).
    \end{align}
    With that, we can conclude that both derivative formulas coincide, independently of the choice of~$\Kcalb$, which finishes the proof.
\end{proof}

\section{Coordinate-wise Optimality under Successive Reduction}\label{appendix:coordinatewise_globalization}
This appendix concerns a hypothetical globalization of simultaneous updates based on exact separable models, without a proximal term. It is not a convergence result for the numerical algorithm used in \cref{section:numerical_experiments}.

\begin{proposition}[Successive reduction of exact-model updates]\label{prop:coordinatewise_reduction}
Let $D=D_1\times\cdots\times D_N$, where each $D_i$ is a nonempty compact subset of a finite-dimensional real vector space, and let $J\colon D\to\RR$ be continuous. At $x\in D$, suppose the exact coordinate-wise problems are solved globally, yielding
\begin{align}\label{eq:best_coordinate_improvement}
    t_i\in\argmin_{t\in D_i}J\bigl(x+(t-x_i)e_i\bigr),\qquad
    \delta_i(x)=J(x)-J\bigl(x+(t_i-x_i)e_i\bigr)\geq0,\qquad
    \delta(x)=\max_i\delta_i(x).
\end{align}
If $\delta(x)=0$, stop. Otherwise first test the simultaneous update of all $t_i\ne x_i$. If it fails
\begin{align}\label{eq:coordinatewise_decrease_test}
    J(x^+)\leq J(x)-\eta\delta(x),\qquad 0<\eta\leq1,
\end{align}
retain the better half of the proposed coordinate changes according to $\delta_i(x)$, including a maximizer of $\delta_i$, and repeat the test. Continue until the test succeeds. Then every accumulation point of the resulting iteration is a coordinate-wise minimum of $J$. If all $D_i$ are finite, the procedure terminates at such a point after finitely many accepted iterations.
\end{proposition}
\begin{proof}
The ranking and successive reduction preserve a coordinate $i^*$ for which $\delta_{i^*}(x)=\delta(x)>0$. At the latest, when only this change remains, the exact model and the definition of $t_{i^*}$ give $J(x^+)=J(x)-\delta(x)$, so \cref{eq:coordinatewise_decrease_test} holds. Hence the inner reduction terminates and every accepted step satisfies that test. Since $J$ is bounded below on compact $D$,
\begin{align}
    \eta\sum_{k=0}^{n-1}\delta(x^k)\leq J(x^0)-J(x^n)\leq J(x^0)-\min_D J,
\end{align}
which implies $\delta(x^k)\to0$. Compactness of each $D_i$ and continuity of $J$ imply continuity of $\delta_i$ and $\delta$ (minimum of a jointly continuous function over a fixed compact set). Thus any accumulation point $x^*$ satisfies $\delta(x^*)=0$, precisely the coordinate-wise minimality condition. If $D$ is finite, every accepted step strictly decreases $J$ because $\delta(x)>0$. No design can be revisited, and termination is therefore finite.
\end{proof}

\clearpage
\newrefcontext[sorting=nyt]
\printbibliography

\end{document}